\documentclass[11pt,reqno]{amsart}
\usepackage[a4paper,vmargin=30mm,hmargin=28mm]{geometry}
\usepackage[T1]{fontenc}
\usepackage{lmodern, soul,comment}
\usepackage{microtype}
\usepackage{mathtools,amssymb,mathrsfs}
\usepackage{stmaryrd}
\usepackage{enumitem}
\usepackage{xcolor}
\usepackage{hyperref}
\usepackage{mathtools}
\mathtoolsset{showonlyrefs}
\hypersetup{
  colorlinks=true,
  linkcolor=black,
  citecolor=black,
  urlcolor=black,
  pdfauthor={Shiping Cao, Haiping Jiao, Hua Qiu},
  pdftitle={Homogenization of p-resistance forms on nested fractal graphs}
}

\numberwithin{equation}{section}
\newtheorem{theorem}{Theorem}[section]
\newtheorem{lemma}[theorem]{Lemma}
\newtheorem{proposition}[theorem]{Proposition}
\newtheorem{corollary}[theorem]{Corollary}
\theoremstyle{definition}
\newtheorem{definition}[theorem]{Definition}
\newtheorem{assumption}[theorem]{Assumption}
\theoremstyle{remark}
\newtheorem{remark}[theorem]{Remark}

\newcommand{\R}{\mathbb R}
\newcommand{\Tr}{\operatorname{Tr}}

\newcommand{\cE}{\mathcal E}
\newcommand{\cR}{\mathcal R}
\newcommand{\cQ}{\mathcal Q}
\newcommand{\cS}{\mathcal S}
\newcommand{\cF}{\mathcal F}
\newcommand{\cK}{\mathcal K}
\newcommand{\cW}{\mathcal W}
\newcommand{\cD}{\mathcal D}
\newcommand{\cV}{\mathcal V}

\setlist[enumerate]{leftmargin=2.2em,itemsep=.25em,topsep=.45em}
\setlist[itemize]{leftmargin=1.8em,itemsep=.25em,topsep=.45em}
\allowdisplaybreaks[2]
\title{Stochastic homogenization of $p$-resistance forms on nested fractals}

\author{Shiping Cao}
\address{Department of Mathematics, The Chinese University of Hong Kong, Shatin, Hong Kong}
\email{spcao@math.cuhk.edu.hk}
\thanks{ The research of SC is partially  supported by the General Research Fund grant (project CUHK14304626) from the Hong Kong Research Grant Council,
	and by a direct grant for research (project 4053774) from the Chinese University of Hong Kong.}

\author{Haiping Jiao}
\address{School of Mathematics, Nanjing University, Nanjing, 210093, P. R. China.}
\email{hpjiao@smail.nju.edu.cn}
\date{}

\author{Hua Qiu}
\address{School of Mathematics, Nanjing University, Nanjing, 210093, P. R. China.}
\thanks{The research of HQ was supported by the National Natural Science Foundation of China under its General Program (Grant No.~12471087) and Key Program (Grant No.~12531004).}
\email{huaqiu@nju.edu.cn}


\begin{document}

\begin{abstract}
We study the stochastic homogenization of $p$-energies on infinite nested fractals. We assign i.i.d. random $p$-resistance forms $D_{0,\alpha}$ to level-$0$ cells of an infinite pre-nested fractal, and study the scaling limit of the resulting random $p$-resistance form $\cD^{(0)}$. Under a moment condition on the reciprocal of $D_{0,\alpha}$, we prove that the sequence of rescaled and renormalized random energies is almost surely Mosco convergent to an equally weighted self-similar $p$-resistance form. As an application, when $p=2$, we prove the almost sure weak convergence of the associated rescaled continuous-time random walk. In particular, we allow arbitrarily small random conductances and impose no condition on the upper bound of random conductances.
\end{abstract}

\maketitle
{\small\setlength{\baselineskip}{10pt}\tableofcontents}

\section{Introduction}
Stochastic homogenization is the phenomenon in which random media exhibit deterministic behavior at large scales.  A well-known and important class of examples is provided by random conductance models. In this paper, we study the stochastic homogenization of $p$-resistance forms on infinite nested fractals,  for $1<p<\infty$. We assign independent and identically distributed (i.i.d.) $p$-resistance forms to the cells of an infinite pre-nested fractal and prove the almost sure convergence of the resulting random energy forms after scaling.

In the 1980s, Künnemann \cite{Kunnemann1983} and Kozlov \cite{Kozlov1985} initiated the study of random conductance models on $\mathbb{Z}^d$. Assuming that the conductances are stationary, ergodic, and uniformly elliptic, they proved the weak convergence of the associated random walks to Brownian motion under diffusive rescaling. This phenomenon has since been widely studied and established under weaker assumptions. In the i.i.d. setting, Barlow and Deuschel \cite{BarlowDeuschel2010} removed the assumption of a uniform upper bound on the conductances while retaining a uniform positive lower bound. Later, for general stationary ergodic environments, Andres, Deuschel and Slowik \cite{AndresDeuschelSlowik2015} replaced both bounds with moment conditions, and their result was refined in \cite{BS} under weaker moment assumptions for $d\geq 3$. We also point out that homogenization has been proved for degenerate random conductance models, including those on the infinite cluster of supercritical percolation (see \cite{AndresBarlowDeuschelHambly2013,BergerBiskup2007,MathieuPiatnitski2007,SidoraviciusSznitman2004}).\medskip 

However, the fractal structure of a random medium can give rise to random walks whose large-scale behavior differs from that of Brownian motion. In 1986, Kesten \cite{Kesten1986} proved that simple random walk on the incipient infinite cluster of two-dimensional critical bond percolation is subdiffusive: for some $\varepsilon>0$, the family $\{n^{-1/2+\varepsilon}X_n\}_{n\geq1}$
is tight under the annealed law, where $X_n$ denotes the position of the walk at time $n$ and $X_0=0$. 

\begin{figure}[htp]
\includegraphics[width=3.1cm]{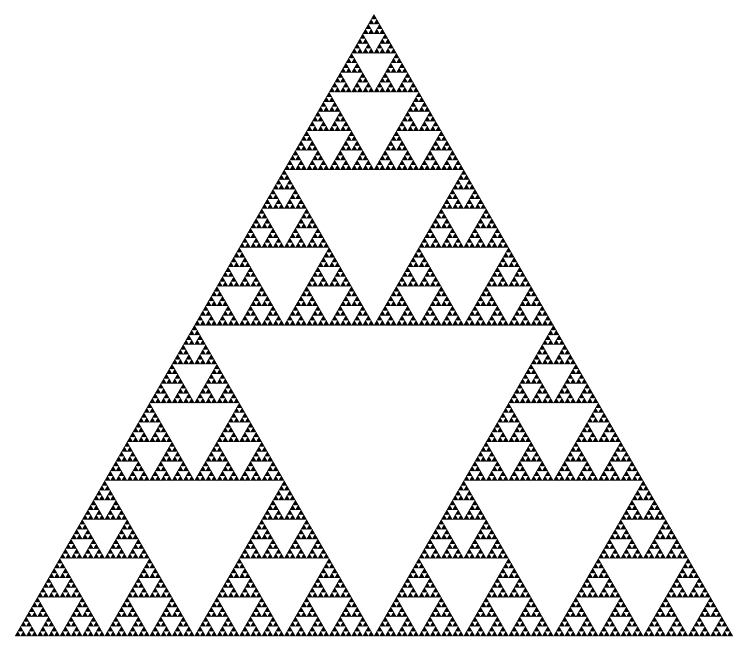}\quad
\includegraphics[width=3cm]{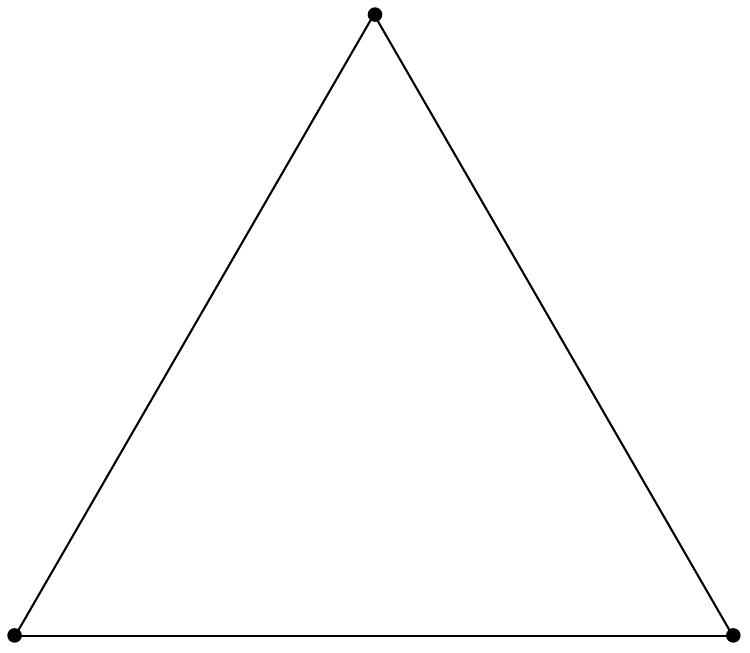}\quad
\includegraphics[width=3cm]{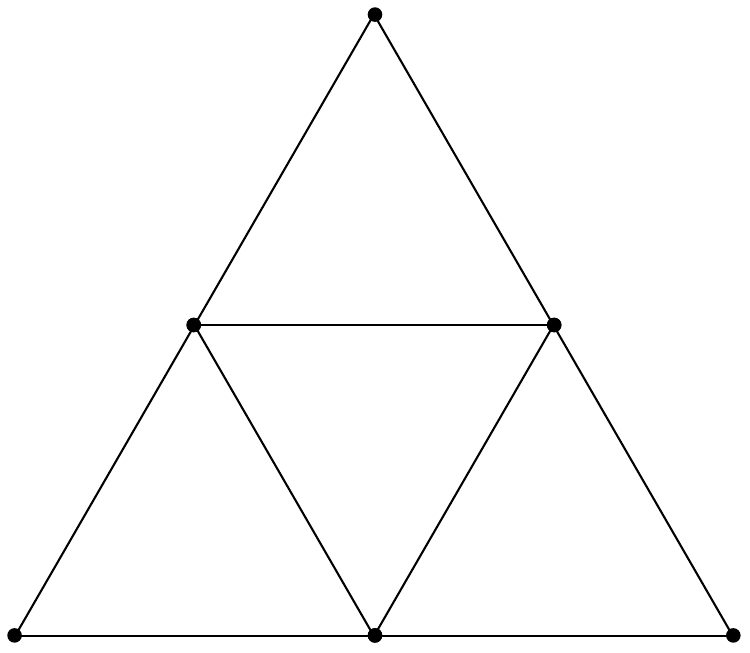}\quad
\includegraphics[width=3cm]{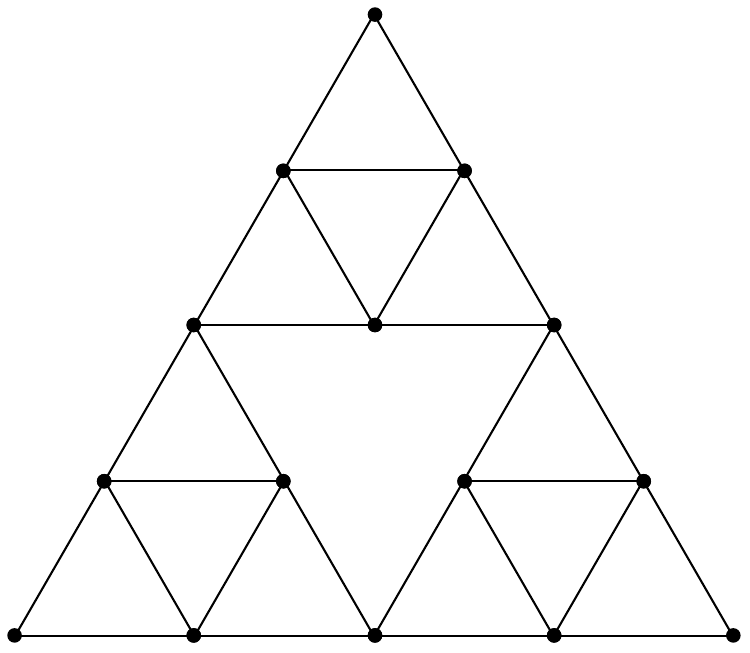}
\caption{The Sierpi\'nski gasket and its graph approximations (also known as pre-Sierpi\'nski gaskets)}\label{fig1}
\end{figure}

Such phenomena motivated the study of diffusion processes on self-similar fractals. In 1987 and 1988, Kusuoka \cite{Kusuoka1987}, Goldstein \cite{Goldstein1987}, and Barlow and Perkins \cite{BarlowPerkins1988} independently constructed Brownian motion on the Sierpi\'nski gasket as a scaling limit of random walks on approximating graphs (see Figure~\ref{fig1}). The construction was later extended probabilistically to nested fractals by Lindstr{\o}m \cite{Lindstrom1990} and to planar and higher-dimensional generalized Sierpi\'nski carpets by Barlow and Bass \cite{BarlowBass1989,BarlowBass1999}.

Note that these limit processes can be viewed as diffusion processes associated with limit energies of discrete electrical networks, and the idea has led to further generalizations. In particular, Kigami \cite{Kigami1993} used vertex graph approximations to construct Dirichlet forms on post-critically finite (p.c.f.) self-similar sets, including nested fractals. Kusuoka and Zhou \cite{KusuokaZhou} used cell graph approximations to construct Dirichlet forms on a class of self-similar sets, including the planar Sierpi\'nski carpet.
For background on Dirichlet forms and on diffusion processes on fractals, see \cite{CF,FukushimaOshimaTakeda2011} and \cite{Barlow1998,Kigami2001}, respectively. 

 We also point out that $p$-energies are natural, non-trivial extensions of Dirichlet forms and have attracted growing interest in analysis on fractals \cite{KajinoShimizuSurvey,Kigami2020}. Very recently, Kajino and Shimizu  \cite{KajinoShimizu2026} developed a general framework for $p$-resistance forms and studied their properties. On concrete fractals, earlier in 2004, Herman, Peirone and Strichartz \cite{HermanPeironeStrichartz2004} constructed $p$-energies on Sierpi\'nski gasket type fractals, and recently this construction was extended by Cao, Gu and Qiu \cite{CaoGuQiu2022} to p.c.f. self-similar sets using techniques of Kusuoka and Zhou \cite{KusuokaZhou}, and to a class of conductively homogeneous compact metric spaces by Kigami \cite{Kigami2023}. Later, Shimizu~\cite{Shimizu2024} constructed self-similar $p$-energies on the standard planar Sierpi\'nski carpet, and Murugan and Shimizu~\cite{MuruganShimizu2025} gave a construction for all $1<p<\infty$. These works provide useful tools for the study of the fractal nature of porous media. 
\medskip

Self-similar sets provide simple and natural models for studying homogenization in porous media. 

In 1993, Kozlov \cite{Kozlov1993} initiated the study of homogenization phenomena on fractal graphs. He showed that certain anisotropic energies on the pre-Sierpi\'nski gaskets (each consisting of the vertices of all triangular cells at a fixed scale) converge, after renormalization, to a constant multiple of the standard energy. Similar conclusions extend to many nested fractals by combining convergence results for the renormalization map \cite{Peirone2000} with uniqueness results for self-similar energies up to a constant factor \cite{Sabot1997}.

\begin{figure}[htp]
\includegraphics[width=3.8cm]{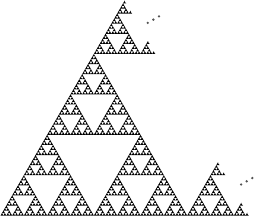}\qquad
\includegraphics[width=3.8cm]{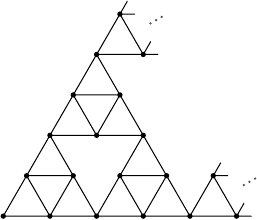}
\caption{The infinite Sierpi\'nski gasket and level-$0$  infinite pre-Sierpi\'nski gasket}\label{fig2}
\end{figure}

Later, Kumagai and Kusuoka \cite{KumagaiKusuoka1996} studied stochastic homogenization on infinite pre-nested fractals (each consisting of the boundary vertices of all cells at a fixed scale). They assigned i.i.d. conductance configurations to the cells and proved the weak convergence of the associated rescaled random walks to a constant time change of Brownian motion on the corresponding infinite fractal under a uniform two-sided bound assumption on the cell forms and an additional condition on the renormalization map. See Figure \ref{fig2} for an illustration of the infinite Sierpi\'nski gasket and the associated level-$0$ infinite pre-Sierpi\'nski gasket.

This result was later refined by Croydon, Hambly and Kumagai \cite{CroydonHamblyKumagai2017}, who replaced the uniform upper bound with a second moment condition on the renormalized conductances.
However, their argument still relies on a uniform lower bound on conductances. So far, in contrast to the classical results \cite{AndresDeuschelSlowik2015} and \cite{BS} on $\mathbb{Z}^d$, the only stochastic homogenization theorem on a nested fractal allowing  arbitrarily small conductances has been established for Vicsek sets \cite{HamblyMetz1998}, where the argument relies on their tree-like structures. \medskip

In this paper, we prove stochastic homogenization for $p$-resistance forms on nested fractals with i.i.d. random $p$-resistance forms assigned to cells, assuming only a moment condition on the reciprocal of cell energies.

More precisely, let $K$ be a nested fractal with essential boundary $V_0$, and let $\cK$ be the associated infinite nested fractal. Let $\cV_0$ be the associated level-$0$ infinite pre-nested fractal. We assign i.i.d. random $p$-resistance forms $D_{0,\alpha}$ to the level-$0$ cells of $\cV_0$, where $\alpha$ is the address of the corresponding cell, and denote the resulting random form on $\cV_0$ by $\cD^{(0)}$. With respect to a non-degenerate fixed point $E_0$ of the energy renormalization map $\cR$, we define the reciprocal ratio of a $p$-resistance form $E$ on $V_0$ by
\[
\llbracket E\rrbracket
:=\sup_{u\in\mathbb R^{V_0}\setminus\mathbb R\mathbf1}
\frac{E_0(u)}{E(u)}.
\]
Under the assumption
\begin{equation}\label{eq:intro-moments}
	\mathbb E\llbracket D_{0,\alpha}\rrbracket^{1/(p-1)}<\infty,
\end{equation}
we prove that the renormalized energy forms $\cD^{(n)}$ on the level-$n$ infinite pre-nested fractal $\cV_n$ (each obtained by rescaling $\cD^{(0)}$ on $\cV_0$) are almost surely Mosco convergent to a deterministic self-similar form $\cE$ on $\cK$ (see Theorem \ref{thm:mosco-homogenization}). To prove this, we establish a $p$-version of Peirone's convergence theorem \cite{Peirone2000}, showing that, for every $p$-resistance form $E$ on $V_0$, $\cR^n(E)$ converges to a fixed point of $\cR$ (Theorem~\ref{thm:peirone}).

When $p=2$, we further show that, for almost every environment, the associated rescaled random walks converge weakly to the diffusion associated with $\cE$. See Theorem~\ref{thm:process-homogenization}. \medskip 

Note that we allow arbitrarily small conductances, and require no condition concerning the upper bound of conductances. Our results greatly extend earlier theorems \cite{KumagaiKusuoka1996,HamblyMetz1998,CroydonHamblyKumagai2017} even in the $p=2$ case, and our moment condition is optimal even for the unit interval.\medskip

The rest of the paper is organized as follows.

In Section 2, we introduce precise definitions of nested fractals, $p$-resistance forms and the random environments. At the end of this section, we state our main theorems concerning the stochastic homogenization of $p$-resistance forms on infinite nested fractals. 

Sections 3 and 4 are preparations for the proof. In Section 3, we show that various moment conditions imposed on the level-$0$ pre-nested fractal are strengthened or preserved under successive boundary traces of forms on higher-level vertex sets after renormalization. In Section 4, we show that the energy of a harmonic function does not concentrate on a single cell, and use this to establish concentration estimates for the boundary trace $D_{n,\alpha}$ of the renormalized energy form $D^{(n)}_\alpha$ on each $0$-cell.

In Section 5, we establish the homogenization result that $D_{n,\alpha}$  converges almost surely to a deterministic form $E_+$. The proof involves a Peirone-type convergence theorem, whose proof is given in Appendix B.

Finally, in Section 6, we prove the almost sure Mosco convergence of $\cD^{(n)}$. When $p=2$, in Section 7, we take a step further to show the almost sure weak convergence of the associated random walks.

We also point out that in Appendix A we prove a law of large numbers on rooted trees, which will be used in Section 4 and Section 5.

Throughout the paper, we use $:=$ as a way of definition.  For two functions $f$ and $g$, we say $f\lesssim g$ if there is a positive constant $C$ such that $f\leq Cg $; we say $f\gtrsim g$ if there is a positive constant $C$ such that $f\geq Cg$; we say $f\asymp g$ if $g\lesssim f\lesssim g$. We also write $f\wedge g$ for $\min{\{f,g\}}$, and $f\vee g$ for $\max {\{f,g\}}$.

\section{\texorpdfstring{Geometric setting and random environments}{Geometric setting and random environments}}
\label{sec:settings}
In this section, we  review the definition of nested fractals and $p$-energies. Then, we introduce our homogenization theorem.

\subsection{Nested fractals}
Let 
\[
\psi_i(x)=r\,U_i x+b_i,\qquad i\in I=\{1,\ldots,N\},
\]
be a finite collection of distinct self-similarities on $\mathbb{R}^d$, $d\geq 1$. Here, we assume that $N\geq 2$, $0<r<1$, and for each $i\in \{1,\ldots,N\}$, $b_i\in\mathbb{R}^d$ and $U_i:\R^d\to\R^d$ is an orthogonal linear transformation. There is a unique non-empty compact set $K\subset\mathbb R^d$ such that 
\[
K=\bigcup_{i=1}^N\psi_i(K).
\]
$K$ is called the self-similar set associated with the iterated function system $\{\psi_i\}_{i=1}^N$.

 In this paper, we also assume that $K$ satisfies the following conditions (A-0)--(A-3). Such a self-similar set $K$ was called a \textit{nested fractal} by Lindstr{\o}m \cite{Lindstrom1990} in 1990.\smallskip 

\noindent(A-0): (\textbf{Open set condition}) There exists a nonempty open set $O$ such that $\bigcup_{i \in I} \psi_i(O) \subset O$, and $\psi_i(O) \cap \psi_j(O)=\emptyset$ for any $i, j \in I, i \neq j$.\smallskip \\
\noindent(A-1): (\textbf{Connectivity}) $K$ is connected.\smallskip 

We write $q_i$ for the fixed
point of $\psi_i$.  A fixed point $q_i$ is called
\emph{essential} if there exist distinct $j,k$ and a fixed point $q_\ell$
such that
$\psi_j(q_i)=\psi_k(q_\ell)$.  The set of essential fixed points is
denoted by $V_0$.  After relabeling the similarities, we assume that
\[
V_0=\{q_1,\ldots,q_M\}\ \hbox{ with }2\le M\le N.
\]
For $n\geq 1$ and $w=w_1w_2\cdots w_n\in I^n$, we write $\psi_w=\psi_{w_1}\circ\psi_{w_2}\circ\cdots\circ\psi_{w_n}$. We also use the convention that $I^0=\{\emptyset\}$ with $\psi_\emptyset$ being the identity map. Set 
\[
V_n:=\bigcup_{w\in I^n}\psi_w(V_0)\quad\hbox{ for }n\geq 1.
\]\smallskip 

\noindent(A-2): (\textbf{Symmetry}) For distinct elements $x, y \in \mathbb{R}^d$, $H_{xy}$ denotes the hyperplane given by $H_{x y}=\left\{z \in \mathbb{R}^d:|z-x|=|z-y|\right\}$, and $U_{x y}$ denotes the reflection with respect to $H_{x y}$. For any $x,y\in V_0$ with $x\neq y$, $n\geq 0$ and $w\in I^n$, there is $v\in I^n$ so that $U_{xy}\big(\psi_w(V_0)\big)=\psi_v(V_0)$ and $v=w$ if $\psi_w(V_0)$ contains elements on both sides of $H_{x y}$.\smallskip 

\noindent(A-3): (\textbf{Nesting}) If $n\geq 1$ and $w,v$ are distinct words in $I^n$, then $\psi_w(K)\cap \psi_v(K)=\psi_w(V_0)\cap \psi_v(V_0)$. \smallskip 

{In addition to the above conditions, we also point out a useful property (see \cite[Section~9, Corollary]{BandtKeller1991} and also \cite[Lemma~3.1.4]{Kohl2022}):}
\begin{equation}\label{eq:outer-boundary-intersection}
	\#\bigl(\psi_i(V_0)\cap V_0\bigr)\le1,
	\qquad i\in I.
\end{equation}

Next, we introduce an infinite fractal $\mathcal{K}$ and its discrete approximations $\mathcal{V}_n$. Fix an index $i_0\in I$ and form the one-sided blow-up
\[
\mathcal K:=\bigcup_{\ell\ge0}\psi_{i_0}^{-\ell}(K).
\]
{Correspondingly, for $n\geq 0$, we define the \textit{level-$n$ infinite pre-nested fractal} as 
\begin{equation}\label{eq:infinite-skeletons}
	\mathcal V_n:
	=\bigcup_{\ell\ge0}\psi_{i_0}^{-\ell}(V_{n+\ell}).
\end{equation}}

Let $I^*=\bigcup_{n=0}^\infty I^n$ be the set of finite words. For $w\in I^*$, define
\begin{equation}\label{eq:direct-cell-coding}
	\Psi_w:=\psi_{i_0}^{-|w|}\circ\psi_w.
\end{equation}
Then, the map $w\longmapsto \Psi_w(K)$ is a bijection from
\[
\cW:=\{\emptyset\}\cup
\{w\in I^*\setminus\{\emptyset\}:w_1\ne i_0\}
\]
onto the family $\{\Psi_v(K):v\in I^*\}$ of cells of $\cK$ at the original scale.
For $w\in\cW$ and $v \in I^n$, we call $\Psi_w\psi_v(K)$ an \textit{$n$-cell}. It is worth pointing out that 
\[
\psi_{i_0}^{-n}(\Psi_w\psi_v(K))=\Psi_{wv}(K), 
\]
since  
\begin{equation}\label{eq:descendant-cell-identity}
	\Psi_{wv}
	=\psi_{i_0}^{-n}\circ \Psi_w\circ\psi_v.
\end{equation}

\subsection{\texorpdfstring{$p$}{p}-resistance forms}
Fix $p\in(1,\infty)$.  We use the notion introduced by Kajino and
Shimizu \cite[Definition~6.1]{KajinoShimizu2026}. {Throughout the paper, we denote by $\R^A$ the set of real-valued functions on $A$, and denote by $C(A)$ the set of continuous functions. Moreover, $C_c(A)$ is the set of compactly supported continuous functions.}

\begin{definition}\label{def:p-resistance}
	Let $X$ be a nonempty set with $\#X\geq 2$, let $\cF$ be a linear subspace of
	$\mathbb R^X$ containing the constants, and let
	$E:\cF\to[0,\infty)$.  The pair $(E,\cF)$ is a
	\emph{$p$-resistance form} on $X$ if the following conditions hold.
	\begin{enumerate}[label=\textup{(RF\arabic*)}]
		\item $E^{1/p}$ is a seminorm on $\cF$, and its kernel is
		$\mathbb R\mathbf1$.
		\item The quotient space $\cF/\mathbb R\mathbf1$, equipped with
		$E^{1/p}$, is a Banach space.
		\item The space $\cF$ separates the points of $X$.
		\item For all $x,y\in X$,
		\[
		R_E(x,y):=
		\sup_{u\in\cF\setminus\mathbb R\mathbf1}
		\frac{|u(x)-u(y)|^p}{E(u)}<\infty.
		\]
		\item The \textit{generalized $p$-contraction property} holds: if
		$n_1,n_2\in\mathbb N$, $q_1\in(0,p]$, $q_2\in[p,\infty]$, 
		and $T=(T_1,T_2,\cdots,T_{n_2}):\mathbb R^{n_1}\to\mathbb R^{n_2}$ satisfies $T(0)=0$ and
		\begin{equation}\label{nonexpansive}
		\|T(x)-T(y)\|_{\ell^{q_2}}
		\le\|x-y\|_{\ell^{q_1}},
		\qquad x,y\in\mathbb R^{n_1},
		\end{equation}
		then for every $(u_1,\ldots,u_{n_1})\in\cF^{n_1}$ the coordinate
		functions $T_k(u_1,\ldots,u_{n_1})$ belong to $\cF$ and
		\[
		\left\|\left(
		E\bigl(T_k(u_1,\ldots,u_{n_1})\bigr)^{1/p}
		\right)_{k=1}^{n_2}\right\|_{\ell^{q_2}}
		\le
		\left\|\left(E(u_j)^{1/p}\right)_{j=1}^{n_1}
		\right\|_{\ell^{q_1}}.
		\]
	\end{enumerate}
\end{definition}

\begin{remark}\label{rem:finite-p-resistance}
On a finite nonempty set $V$ with $\#V\geq 2$, any function $E:\R^V\to [0,\infty)$ satisfying conditions (i), (ii) and (iii) below is automatically a $p$-resistance form. See \cite[Example 2.2]{KajinoShimizuSurvey}.
\begin{enumerate}[label=\textup{(\roman*)}]
	\item $E(u)=0$ exactly for constant $u$;
	\item $E(au+b\mathbf1)=|a|^pE(u)$ and $E(u+v)^{1/p}\le E(u)^{1/p}+E(v)^{1/p}$;
	\item the generalized contraction property \textup{(RF5)} holds.
\end{enumerate}
\end{remark}

\begin{remark}
A function $E:\R^V/\R{\bf 1}\to\R$ is called a \textit{$p$-homogeneous functional} if $E(\lambda u)=\lambda^pE(u)$ for every $\lambda>0$ and $u\in\R^V/\R{\bf 1}$. 
By (RF1), we can identify a $p$-resistance form $E$ with a $p$-homogeneous functional on $\R^V/\R{\bf 1}$.
\end{remark}

For the rest of the section, let $V$ be a nonempty finite set with $\#V\geq 2$ and write
$\cQ_p(V)$ for the class of $p$-resistance forms on it. We also use $\widetilde{\mathcal{Q}}_p(V)$ to denote the collection of functionals $\R^V\to [0,\infty)$ satisfying (ii) and (iii) of Remark \ref{rem:finite-p-resistance}. Note that (ii) already implies that $E({\bf 1})=0$.

We recall the following fact proved in \cite[Theorems~3.7 and~6.13 and Proposition~6.26]{KajinoShimizu2026}.
\begin{proposition}\label{prop:derivative}
	Let $E\in\cQ_p(V)$ and $u,v\in\mathbb R^V$.  The directional derivative
	of $E$ at $u$ in the direction $v$ exists.  Under the normalization
	\[
	E(u;v):=\frac1p\left.\frac{d}{dt}E(u+tv)\right|_{t=0},
	\]
	the map $v\mapsto E(u;v)$ is linear, and
	\begin{equation}\label{eq:euler-holder}
		E(u;u)=E(u),\qquad E(u;\mathbf1)=0,\qquad
		|E(u;v)|\le E(u)^{1/p'}E(v)^{1/p},
		\quad p':=\frac p{p-1}.
	\end{equation}
	If $B\subset V$ is nonempty with $\#B\geq 2$ and $f\in\mathbb R^B$, the minimization of $E(h)$ over
	$h|_B=f$ has a unique solution.  This solution takes its values between
	$\min_Bf$ and $\max_Bf$.
\end{proposition}

For $x\in V$, write $\mathbf1_x$ for the indicator function of $\{x\}$ and define the \textit{$p$-current flowing out of $x$} induced by the potential $u\in \mathbb R^V$ as
\[
J_E(u,x):=E(u;\mathbf1_x).
\]
By Proposition~\ref{prop:derivative},
\begin{equation}\label{eq:current-euler}
	E(u)=\sum_{x\in V}J_E(u,x)u(x),
	\qquad
	\sum_{x\in V}J_E(u,x)=0.
\end{equation}\smallskip

\begin{lemma}\label{lem:state-closure}
Let $V$ be a finite set, let $(\Omega,\mathcal{G},\nu)$ be a measure space, and let $\omega\mapsto E_\omega$ be a measurable map from $\Omega$ to $\widetilde{\mathcal{Q}}_p(V)$ with respect to the pointwise topology, i.e. $E_\omega(u)$ is a measurable function on $\Omega$ for each $u\in\R^V$. 

Assume $\int_\Omega E_\omega(u)\,d\nu(\omega)<\infty$ for any $u\in\R^V$, and define $\int_\Omega E_\omega\,d\nu(\omega):\R^V\to[0,\infty)$ by 
\[
(\int_\Omega E_\omega d\nu(\omega))(u):=\int_\Omega E_\omega(u) d\nu(\omega)\quad\hbox{ for }u\in\R^V.
\]
Then $\int_\Omega E_\omega d\nu(\omega)\in \widetilde{\mathcal{Q}}_p(V)$. 

Moreover, if $\nu(\{\omega\in \Omega:\,E_\omega(u)\neq 0\})>0$ for every $u\in \R^V\setminus \R{\bf 1}$, then  $\int_\Omega E_\omega d\nu(\omega)\in \mathcal{Q}_p(V)$.
\end{lemma}
\begin{proof}
Write $E=\int_\Omega E_\omega d\nu(\omega)$ for short. We need to show that $E$ satisfies (ii) and (iii) of Remark \ref{rem:finite-p-resistance}.

The triangle inequality for $E^{1/p}$ follows from the triangle inequality for $E_\omega^{1/p},\omega\in \Omega$, and Minkowski's inequality:
\begin{align*}	
E(u+v)^{1/p}&=\big(\int_\Omega E_\omega(u+v) d\nu(\omega)\big)^{1/p}=\|E_\omega(u+v)^{1/p}\|_{L^p(\nu)}\\
&\leq \|E_\omega(u)^{1/p}\|_{L^p(\nu)}+\|E_\omega(v)^{1/p}\|_{L^p(\nu)}=
E(u)^{1/p}+E(v)^{1/p}. 
\end{align*}
$E(au+b\mathbf1)=|a|^pE(u)$ is straightforward to show. So (ii) holds. 

Next, let $\boldsymbol u=(u_1,\ldots,u_{n_1})$, and let $T=(T_1,\ldots,T_{n_2}):\R^{n_1}\to\R^{n_2}$ be a contraction as in Definition~\ref{def:p-resistance}. Write 
\[
f_k(\omega):=E_\omega(T_k(\boldsymbol u))^{1/p}
\ \hbox{ and }\ 
g_j(\omega):=E_\omega(u_j)^{1/p}\quad\hbox{ for }\omega\in \Omega.
\]
Then, the Minkowski inequality and the generalized contraction
property of each $E_\omega$ give 
\begin{align*}
\left\|\left(E\big(T_k(\boldsymbol u)\big)^{1/p}\right)_{k=1}^{n_2}\right\|_{\ell^{q_2}}=\left\|\left(\|f_k\|_{L^p(\nu)}\right)_{k=1}^{n_2}\right\|_{\ell^{q_2}}
	&\le
	\left\|\left\|(f_k)_{k=1}^{n_2}\right\|_{\ell^{q_2}}\right\|_{L^p(\nu)}\\
	&\le
	\left\|\left\|(g_j)_{j=1}^{n_1}\right\|_{\ell^{q_1}}\right\|_{L^p(\nu)}\\
	&\le
	\left\|\left(\|g_j\|_{L^p(\nu)}\right)_{j=1}^{n_1}\right\|_{\ell^{q_1}}=\left\|\left(E(u_j)^{1/p}\right)_{j=1}^{n_1}\right\|_{\ell^{q_1}}.
\end{align*}
Here the first inequality uses $q_2 \ge p$, the last uses $q_1 \le p$, and the
middle inequality uses \textup{(RF5)}. This is exactly the generalized contraction property for $E=\int_\Omega E_\omega d\nu(\omega)$.

Finally, under the additional assumption, $E_\omega(u)>0$ on a set of positive $\nu$-measure for each nonconstant $u$, so $E(u)>0$ and condition (i) of Remark \ref{rem:finite-p-resistance} holds.
\end{proof}

\subsection{Renormalization operator $\mathcal{R}$ on $V_0$}\label{sec:2.3}
If $E$ is a $p$-resistance form on $V$ and $\emptyset\ne B\subset V$ with $\#B\geq 2$, define its trace to $B$ by
\begin{equation}\label{eq:finite-trace-definition}
	(\Tr_BE)(u):=
	\inf\{E(h):h\in\mathbb R^V,\ h|_B=u\}
	\quad\hbox{ for }u\in\mathbb R^B.
\end{equation}
By \cite[Theorem~6.13]{KajinoShimizu2026}, $\Tr_BE$ is a
$p$-resistance form on $B$. The minimizer of the right-hand side is called the harmonic extension of $u$ with respect to $E$. Harmonic extension is monotone in the boundary values (see \cite[Proposition~6.26]{KajinoShimizu2026}).

\begin{lemma}\label{harmonicmontone} 
Let $V$ be a finite set, $\emptyset\ne B\subsetneq V$ with $\#B\geq 2$ and $E\in\mathcal{Q}_p(V)$. For $u,v\in\mathbb R^B$, let $h_u,h_v$ be the harmonic extensions of $u,v$ with respect to $E$. If $u\leq v$, then $h_u\leq h_v$.  
\end{lemma}
\begin{proof}
Note that $(s,t)\mapsto(s\wedge t,s\vee t)$ is non-expansive from $\ell^p$ to $\ell^p$ in the sense of \eqref{nonexpansive}, so
\[
E(h_u\wedge h_v)+E(h_u\vee h_v)\leq E(h_u)+E(h_v),
\] 
by (RF5). Noticing that $(h_u\wedge h_v)|_{B}=u$ and $(h_u\vee h_v)|_B=v$, we also have 
\[
E(h_u\wedge h_v)\geq E(h_u)\hbox{ and }E(h_u\vee h_v)\geq E(h_v). 
\]
Hence, $E(h_u\wedge h_v)=E(h_u)$. By the uniqueness of harmonic extensions, we have $h_u\wedge h_v=h_u$, so $h_u\leq h_v$.
\end{proof}

Return to the nested fractal setting. For forms
$E_1,\ldots,E_N\in\cQ_p(V_0)$, consider the sum 
\[
h\mapsto\sum_{i=1}^NE_i(h\circ\psi_i|_{V_0}).
\]
Each summand belongs to $\widetilde{\mathcal Q}_p(V_1)$. If the sum vanishes, $h$ is constant on each $\psi_i(V_0)$. By connectivity (A-1) and nesting (A-3), these constants agree on $V_1$. Lemma~\ref{lem:state-closure} therefore shows that the sum is a $p$-resistance form on $V_1$. Its trace onto $V_0$,
\begin{equation}\label{eq:unnormalized-trace}
	T(E_1,\ldots,E_N)(u)
	:=\inf\left\{\sum_{i=1}^N E_i(h\circ\psi_i|_{V_0}):
	h\in\mathbb R^{V_1},\ h|_{V_0}=u\right\}\quad\text{for }u\in\mathbb R^{V_0},
\end{equation}
is also a $p$-resistance form. 

\begin{remark}
For a finite word $w\in I^*$, we abbreviate $h\circ\psi_w|_{V_0}$ to $h\circ\psi_w$ if no confusion occurs. 
\end{remark}

Let $h$ minimize the energy in \eqref{eq:unnormalized-trace}. For $x\in V_1\setminus V_0$, consider $h+t\mathbf1_x,t\in\mathbb{R}$, which preserves the boundary values. Differentiating at $t=0$ gives Kirchhoff's law,
\begin{equation}\label{eq:kirchhoff}
	\sum_{i:x\in\psi_i(V_0)}
	J_{E_i}\bigl(h\circ\psi_i,\psi_i^{-1}x\bigr)=0. 
\end{equation}

A nested fractal is an affine nested fractal in the terminology of
\cite[Section~8.4]{KajinoShimizu2026}.  Applying
\cite[Theorems~8.38 and~8.50]{KajinoShimizu2026} to equal cell weights gives
a \emph{reference $p$-resistance form} $E_0\in\cQ_p(V_0)$ and an \emph{energy renormalization factor} $\rho>1$ such that
\begin{equation}\label{eq:fixed-normalization}
\cR(E_0):=\,\cR(E_0,\ldots,E_0)=E_0,
\end{equation}
where
\begin{equation}\label{eq:renormalization-operator}
\cR(E_1,\ldots,E_N):=\rho\,T(E_1,\ldots,E_N).
\end{equation}
For $E\in\cQ_p(V_0)$, we also write $\cR(E):=\cR(E,\ldots,E)$ for short.
More generally, for forms $(E_w)_{w\in I^k},k\geq 1$, we write
$\cR^k((E_w)_{w\in I^k})$ for the trace onto $V_0$ of the form 
\[
f\longmapsto\rho^k\sum_{w\in I^k}E_w(f\circ\psi_w)
\]
in $\cQ_p(V_k)$.
In particular, $\cR(E_1,\ldots,E_N)=\cR((E_i)_{i\in I})$.

Let $\{E_i\}_{i=1}^N, \{G_i\}_{i=1}^N$ be $p$-resistance forms on $V_0$, and let $t,s> 0$. We have the following elementary properties of $\cR$:
\begin{enumerate}
	\item(\textbf{Homogeneity}) $\cR(tE_1,\ldots,tE_N)=t\cR(E_1,\ldots,E_N)$;
	\item(\textbf{Order-preserving}) $\cR(E_1,\ldots,E_N)\le \cR(G_1,\ldots,G_N)$ if $E_i\leq G_i, 1\le i\le N$;
	\item(\textbf{Concavity}) $\cR(tE_1+sG_1,\ldots,tE_N+sG_N)\ge t\cR(E_1,\dots,E_N)+s\cR(G_1,\ldots,G_N)$.
\end{enumerate}
If $E_i\ge aE_0$ for every $i$, order preservation and
\eqref{eq:fixed-normalization} give
\begin{equation}\label{eq:lower-trace}
	\cR(E_1,\ldots,E_N)\ge aE_0.
\end{equation}

Consider the quotient space $\mathbb R^{V_0}/\mathbb R\mathbf1$. Define
\begin{equation}\label{eq:unit-section}
	\cS:=\{u\in \mathbb R^{V_0}/\mathbb R\mathbf1:E_0(u)=1\}.
\end{equation} From time to time, we identify each element in $\cS$ with its unique representative $u$ satisfying $\sum_{x\in V_0}u(x)=0$. 
This is a compact sphere of dimension $M-2$. Then 
\begin{equation}\label{eq:form-norm}
	\|G\|:=\sup_{u\in\cS}|G(u)|
\end{equation}
defines a norm on the space of continuous $p$-homogeneous functionals on $\mathbb R^{V_0}/\mathbb R\mathbf1$.
For a $p$-resistance form $E\in\cQ_p(V_0)$,
\[
\|E\|=\sup_{u\notin\mathbb R\mathbf1}\frac{E(u)}{E_0(u)}.
\]

We point out that $\cR$ is continuous with respect to the norm $\|\cdot\|$. 

\begin{lemma}\label{lem:trace-lipschitz} 
There is a constant $C_{\cR}<\infty$, depending only on $E_0$, $p$, $\rho$ and $K$, such that
\begin{equation}\label{eq:trace-lipschitz}
	\|\cR(E_1,\ldots,E_N)-\cR(G_1,\ldots,G_N)\|
	\le C_{\cR}\sum_{i=1}^N\|E_i-G_i\|
\end{equation}
for all $E_i,G_i\in\cQ_p(V_0)$, $1\le i\le N$.
\end{lemma}
\begin{proof}
Fix $u\in\cS$, and let $h_E$ and $h_G$ be the two minimizing extensions in
\eqref{eq:unnormalized-trace}. Then
\begin{align}
	\cR(E_1,\ldots,E_N)(u)-\cR(G_1,\ldots,G_N)(u)
	&\le\rho\sum_i(E_i-G_i)(h_G\circ\psi_i)\nonumber\\
	&\le\rho\sum_i\|E_i-G_i\|E_0(h_G\circ\psi_i)\label{e:2.16}\\
	\cR(G_1,\ldots,G_N)(u)-\cR(E_1,\ldots,E_N)(u)
	&\le\rho\sum_i(G_i-E_i)(h_E\circ\psi_i)\nonumber\\
	&\le\rho\sum_i\|E_i-G_i\|E_0(h_E\circ\psi_i).\label{e:2.17}
\end{align}
Since $\R^{V_0}/\R{\bf 1}$ is finite dimensional and $\operatorname{osc}_{V_0}(f):=\max_{x\in V_0}f(x)-\min_{x\in V_0}f(x),f\in \R^{V_0}$ defines a norm on $\R^{V_0}/\R{\bf 1}$, we know by norm equivalence that 
\begin{equation}\label{e:2.18}
	C_1^{-1}\operatorname{osc}_{V_0}(f)\leq E_0(f)^{1/p}\leq C_1\operatorname{osc}_{V_0}(f)\quad\hbox{ for every }f\in \R^{V_0}. 
\end{equation}
By \eqref{e:2.18}, the fact that $E_0(u)=1$ and Proposition \ref{prop:derivative}, 
\[
\operatorname{osc}_{V_0}(h_E\circ\psi_i)
\le \operatorname{osc}_{V_0}u\leq C_1\ \hbox{ and }\  \operatorname{osc}_{V_0}(h_G\circ\psi_i)
\le \operatorname{osc}_{V_0}u\leq C_1\ \hbox{ for every }i\in I.
\]
Consequently, by \eqref{e:2.18},
\begin{align}
	E_0(h_E\circ\psi_i)
	\leq C_1^p\operatorname{osc}_{V_0}(h_E\circ\psi_i)^p\le C_1^{2p}\quad\hbox{ for every }i\in I,\label{e:2.19}\\
	E_0(h_G\circ\psi_i)
	\leq C_1^p\operatorname{osc}_{V_0}(h_G\circ\psi_i)^p\le C_1^{2p}\quad\hbox{ for every }i\in I.\label{e:2.20}
\end{align}
The statement follows from \eqref{e:2.16}, \eqref{e:2.17}, \eqref{e:2.19} and \eqref{e:2.20}. 
\end{proof}

We also point out that $\cR^k$ is equicontinuous on compact subsets.
\begin{lemma}\label{lem:compact-annulus}
For $0<a\le b<\infty$, define
\[
[aE_0,bE_0]:=\{E\in\cQ_p(V_0):aE_0\le E\le bE_0\}.
\]
\begin{enumerate}
\item[(a)] $\cR([aE_0,bE_0])\subset [aE_0,bE_0]$.
\item[(b)] $\|\cR^k(E)-\cR^k(G)\|\leq (b/a)\|E-G\|$ for any $E,G\in [aE_0,bE_0]$ and $k\geq 1$. 
\item[(c)] $[aE_0,bE_0]$ is compact in the norm \eqref{eq:form-norm}. 
\end{enumerate}
\end{lemma}
\begin{proof}
(a) is an immediate consequence of the order-preserving property. 

(b). For $E,G\in [aE_0,bE_0]$, one has
\begin{align*}
E\le G+\|E-G\|E_0\le (1+\|E-G\|/a)G.
\end{align*}
By the order-preserving and homogeneity properties of $\cR$,  
\begin{align*}
\cR^k(E)&\leq (1+\|E-G\|/a)\cR^k(G)\\
&\leq \cR^k(G)+\|E-G\|(b/a)E_0.
\end{align*}
Similarly, 
\begin{align*}
	\cR^k(G)\leq \cR^k(E)+\|E-G\|(b/a)E_0.
\end{align*}
Hence, $\|\cR^k(E)-\cR^k(G)\|
\le (b/a)\|E-G\|.$

(c). First, we show that $[aE_0,bE_0]$ is closed in the norm \eqref{eq:form-norm}. In fact, if a sequence $(E_n)_{n\geq1}$ in $[aE_0,bE_0]$ converges in the norm \eqref{eq:form-norm} to a $p$-homogeneous functional $E$, then $E_n$ converges pointwise, so $E$ satisfies $aE_0\leq E\leq bE_0$, (ii) and (iii) of Remark~\ref{rem:finite-p-resistance}. Hence $E$ is a $p$-resistance form on $V_0$. 
	
	Next, since $E^{1/p}$ is a norm on $\mathbb R^{V_0}/\mathbb R\mathbf1$ by (RF1),
	\[
	|E(u)^{1/p}-E(v)^{1/p}|\le E(u-v)^{1/p}\le b^{1/p}E_0(u-v)^{1/p}\quad\hbox{ for every }E\in [aE_0,bE_0].
	\]
	Thus $[aE_0,bE_0]$ is uniformly bounded
	and equicontinuous on $\cS$ with respect to the norm $E_0^{1/p}$. The compactness of $[aE_0,bE_0]$ follows from the Arzel\`a--Ascoli Theorem. 
\end{proof}

\subsection{Random environment}\label{sec:random-environment}
Let $(\Omega,\mathcal{A},\mathbb P)$ be a probability space, where $\Omega$ is the sample space, $\mathcal{A}$ is the $\sigma$-algebra and $\mathbb{P}$ is a probability measure. 

For each $\alpha\in \mathcal{W}$, we let $D_{0,\alpha}=D_{0,\alpha}(\omega)\in \mathcal{Q}_p(V_0),\omega\in \Omega$ be a random $p$-resistance form, measurable in the sense of Lemma \ref{lem:state-closure}. Since $\Psi_{i_0w}=\Psi_w$, 
 we use the convention $D_{0,i_0w}=D_{0,w}$ for $w\in I^*$. Throughout the paper, we assume $D_{0,\alpha},\alpha\in\mathcal{W}$, are independent and identically distributed. 
 
{We view $\cD^{(0)}$, which is defined by
\[
\mathcal{D}^{(0)}(\omega)(u)=\sum_{\alpha\in \mathcal{W}}D_{0,\alpha}(\omega)(u\circ \Psi_\alpha)\qquad\hbox{ for }u\in\R^{\cV_0},
\] 
as a random form on $\mathcal{V}_0$ that may take the value $+\infty$. For $n\geq 0$, we renormalize $\cD^{(0)}$ to define a random form on $\cV_n$: 
\begin{equation}\label{eq:infinite-skeleton-energy}
	{\mathcal D}^{(n)}(\omega)(u)
	:=\rho^n
\mathcal{D}^{(0)}(\omega)(u\circ\psi_{i_0}^n)\qquad\hbox{ for }u\in\R^{\cV_n}.
\end{equation}
We will study the asymptotic behavior of $\cD^{(n)}$ in this paper. 

For the proofs, we will restrict ourselves to each $0$-cell. For each $\alpha\in\cW$, define 
\[
D^{(n)}_\alpha(u)=\rho^n\sum_{w\in I^n}
D_{0,\alpha w}\bigl(u\circ\psi_w\bigr)\qquad\hbox{ for }u\in\R^{V_n}, 
\]
and 
\begin{equation}\label{eq:expanded-random-trace}
	D_{n,\alpha}:=\Tr_{V_0}D^{(n)}_\alpha.
\end{equation}
Alternatively, we can define $D_{n,\alpha}$ by induction:
\begin{equation}\label{eq:random-recursion}
D_{n+1,\alpha}
=\cR\bigl(D_{n,\alpha1},\ldots,D_{n,\alpha N}\bigr)\quad\hbox{ for }n\geq 0\hbox{ and }\alpha\in\cW.
\end{equation}
We will always use the convention that $D^{(n)}_{i_0w}=D^{(n)}_w$ and $D_{n,i_0w}=D_{n,w}$ for every $n\ge0$ and $w\in I^*$.  

With the above definitions, for $n\geq 0$,
\[
\cD^{(n)}(u)=\sum_{\alpha\in\cW}D^{(n)}_\alpha(u\circ\Psi_\alpha)\quad\hbox{ for }u\in\R^{\cV_n}, 
\]
}

and
\begin{equation}\label{eq:infinite-trace-decomposition}
	\bigl(\Tr_{\mathcal V_0}{\mathcal D}^{(n)}\bigr)(u)	=
	\sum_{\alpha\in\cW}
	D_{n,\alpha}
	\bigl(u\circ \Psi_\alpha\bigr)\quad\hbox{ for  }u\in C_c(\cV_0),
\end{equation}
where we set $\bigl(\Tr_{\mathcal V_0}{\mathcal D}^{(n)}\bigr)(u):=\inf\left\{
{\mathcal D}^{(n)}(h):h\in C_c(\mathcal V_n),\ h|_{\mathcal V_0}=u
\right\}.$

\subsection{Homogenization theorem}
We now state the homogenization theorem. For $E\in\cQ_p(V_0)$, we define
\[
\llbracket E\rrbracket:=\sup_{u\in\cS}E(u)^{-1}. 
\]
Note that $\llbracket E\rrbracket=\sup_{u\in\R^{V_0}\setminus\R\mathbf1}\frac{E_0(u)}{E(u)}$.

\begin{assumption}\label{ch3:3-2}\label{ass:random}
	{The original environment satisfies
	\begin{equation}\label{eq:random-assumption}
		\mathbb E\llbracket D_{0,\alpha}\rrbracket^{1/(p-1)}<\infty.
	\end{equation}}
\end{assumption}

Whether Assumption~\ref{ass:random} holds does not depend
on the choice of $E_0$ and $\alpha$.

\begin{theorem}\label{ch3:3-3}\label{thm:random-homogenization}
	Let $K$ be a nested fractal, $1<p<\infty$, and let $\rho$ be as in
	\eqref{eq:renormalization-operator}. Let $\{D_{0,\alpha}\}_{\alpha\in\cW}$ be as
	in Section~\ref{sec:random-environment} and satisfy Assumption~\ref{ass:random}. Then there exists a
	$p$-resistance form $E_+\in\cQ_p(V_0)$ with $\cR(E_+)=E_+$ such that for every $\alpha\in \mathcal{W}$,
	\[
	\lim\limits_{n\to\infty}\|D_{n,\alpha}-E_+\|=0\quad\text{ a.s.}
	\]
	If, in addition, there exists some $q>0$ such that
	{$\mathbb E\|D_{0,\alpha}\|^q<\infty,$}
	then 
	\[
	\lim\limits_{n\to\infty}\mathbb E\|D_{n,\alpha}-E_+\|^t=0\quad\hbox{ for every }1\leq t<\infty.
	\]
\end{theorem}

{ As an application of the above theorem, we will show Mosco convergence of $\cD^{(n)}$.
Let $\mu$ be the natural extension of the equally weighted self-similar measure on $K$ to $\cK$, and $\nu_n$ be the discrete approximating measure on $\cV_n$. Let $\cE$ be the self-similar form on $\cK$ induced by $E_+$. Let $\iota_n:L^t(\cV_n,\nu_n)\to L^t(\cK,\mu),1\le t<\infty$ be the natural linear isometry defined in Subsection \ref{sec:self-similar-measures}.  For precise definitions, see Subsections \ref{sec:self-similar-forms}, \ref{sec:self-similar-measures} and \ref{sec:mosco}.

\begin{theorem}\label{thm:mosco-homogenization}
Suppose that Assumption~\ref{ass:random} holds. Then, almost surely,
for every $1\le t<\infty$, $\cD^{(n)}$ is Mosco convergent to $\cE$ in
$L^t(\cK,\mu)$:
\begin{enumerate}[label=\textup{(\arabic*)}]
\item If $u_n\in L^t(\mathcal V_n,\nu_n)$ and
$\iota_nu_n\rightharpoonup u$ weakly in $L^t(\cK,\mu)$, then
\[
\cE(u)\le\liminf_{n\to\infty}\cD^{(n)}(u_n).
\]
\item For every $u\in L^t(\cK,\mu)$, there are
$u_n\in C_c(\mathcal V_n)$ such that
\[
\|\iota_nu_n-u\|_{L^t(\cK,\mu)}\to0,\qquad
\limsup_{n\to\infty}\cD^{(n)}(u_n)\le\cE(u).
\]
\end{enumerate}
\end{theorem}

This is also Corollary~\ref{cor:infinite-mosco} in Subsection~\ref{sec:mosco}.
For $p=2$, we also obtain almost-sure weak convergence of the
associated processes. See Theorem~\ref{thm:process-homogenization}.}

\section{Moment estimates}\label{sec:moments}
In preparation for the proof, we show how moment bounds on $D_{0,\alpha}$ are strengthened or preserved under successive renormalizations. We also point out that all the results apply to $D_{j,\alpha}$ at any fixed starting level $j\geq 0$, with $D_{j+n,\alpha}$ in place of $D_{n,\alpha}$. For simplicity, our results are stated in terms of $D_{n,\alpha},n\geq 1$ and $D_{0,\alpha}$. 

\subsection{Moment estimates of \texorpdfstring{$\|D_{n,\alpha}\|$}{the upper comparison ratio}} First, we
show a moment-improvement property.
\begin{lemma}\label{lem:moment-upgrade}
	Let $q>0$. There is a constant $C_q<\infty$, depending only on $K,E_0,p,q$, such that
	\begin{equation}\label{eq:ch3-3-1}
		\mathbb E\|D_{1,\alpha}\|^{2q}
		\le C_q\bigl(\mathbb E\|D_{0,\alpha}\|^q\bigr)^2.
	\end{equation}
	
\end{lemma}

\begin{proof}
	We first prove that, for $E_1,\ldots,E_N\in\cQ_p(V_0)$,
	\begin{equation}\label{eq:ch3-pair-bound}
		\|\cR(E_1,\ldots,E_N)\|^2
		\lesssim\sum_{1\le i<j\le N}\|E_i\|\|E_j\|.
	\end{equation}
	Choose $1\le j_0\le N$ such that
	$\|E_{j_0}\|=\max_i\|E_i\|$. Let $x_0$ be the unique point in $V_0\cap \psi_{j_0}(V_0)$ if $V_0\cap \psi_{j_0}(V_0)\neq\emptyset$; otherwise let $x_0$ be an arbitrary point in $V_0$. Then, given $u\in\cS$, define $h\in\mathbb R^{V_1}$ as 
	\[
	h(x)=\begin{cases}
		u(x) &\hbox{ if }x\in V_0,\\
		u(x_0)&\hbox{ if }x\in V_1\setminus V_0. 
	\end{cases}
	\]
	Note that $h\circ\psi_{j_0}$ is constant and
	$\operatorname{osc}_{V_0}(h\circ\psi_i)\le\operatorname{osc}_{V_0}(u)$
	for every $i$. This implies that $E_0(h\circ\psi_i)\lesssim E_0(u)=1$ by
	\eqref{e:2.18}. Hence 
	\begin{align*}
		\cR(E_1,\ldots,E_N)(u)
		&\le\rho\sum_{i=1}^N E_i(h\circ\psi_i)
		=\rho\sum_{i\ne j_0}E_i(h\circ\psi_i)\\
		&\lesssim\sum_{i\ne j_0}\|E_i\|.
	\end{align*}
	Thus $\|\cR(E_1,\ldots,E_N)\|\lesssim\sum_{i\ne j_0}\|E_i\|$.
	Since $\|E_{j_0}\|=\max_i\|E_i\|$, we have
	\begin{align*}
		\|\cR(E_1,\ldots,E_N)\|^2
		&\lesssim\left(\sum_{i\ne j_0}\|E_i\|\right)^2\\
		&=\sum_{i\ne j_0}\|E_i\|^2
		+2\sum_{\substack{i<j\\i,j\ne j_0}}\|E_i\|\|E_j\|\\
		&\le\sum_{i\ne j_0}\|E_i\|\|E_{j_0}\|
		+2\sum_{\substack{i<j\\i,j\ne j_0}}\|E_i\|\|E_j\|\\
		&\le2\sum_{i<j}\|E_i\|\|E_j\|.
	\end{align*}
	This proves \eqref{eq:ch3-pair-bound}.
	
	Now put $E_i=D_{0,\alpha i}$ in \eqref{eq:ch3-pair-bound}. We have
	\[
	\|D_{1,\alpha}\|^2\lesssim
	\sum_{i<j}\|D_{0,\alpha i}\|\|D_{0,\alpha j}\|.
	\]
	Then for some constant $C_q'<\infty$,
	\[
	\|D_{1,\alpha}\|^{2q}\le
	C_q'\sum_{i<j}\|D_{0,\alpha i}\|^q\|D_{0,\alpha j}\|^q.
	\]
	Taking expectations on both sides and then using independence yields
	\eqref{eq:ch3-3-1}.	
	
\end{proof}

Lemma \ref{lem:moment-upgrade} says that once $\mathbb E\|D_{0,\alpha}\|^r<\infty$ for some $r>0$, then for each $q\ge1$ we have $\mathbb E\|D_{n,\alpha}\|^q<\infty$ for all sufficiently large $n$.
The following estimate shows that the $q$-th moment bound is uniform in $n$ when $q\ge 1$. 
\begin{lemma}\label{lem:moment2}
	If $\mathbb E\|D_{0,\alpha}\|^q\le C_*$ for some $C_*<\infty$ and $q\geq 1$, then
	\[
	\sup_n\mathbb E\|D_{n,\alpha}\|^q<\infty.
	\]
\end{lemma}
\begin{proof}
	Fix $v\in\cS$ and let $h$ be its harmonic extension to $V_n$ with
	respect to the form $f\mapsto \sum_{w\in I^n}E_0(f\circ\psi_w)$. Then, by the self-similar identity $\cR^n(E_0)=E_0$,
	\[ 
	\rho^n\sum_{|w|=n}E_0(h\circ\psi_w)=E_0(v)=1.
	\]
	As $D_{n,\alpha}=\cR^n\big((D_{0,\alpha w})_{w\in I^n}\big)$, we have 
	\[
	D_{n,\alpha}(v)
	\le\rho^n\sum_{|w|=n}D_{0,\alpha w}(h\circ\psi_w)
	\le\rho^n\sum_{|w|=n}E_0(h\circ\psi_w)\|D_{0,\alpha w}\|.
	\]
	By Jensen's inequality,
	\[
	\mathbb E D_{n,\alpha}(v)^q \le\rho^n\sum_{|w|=n}E_0(h\circ\psi_w)\mathbb E\|D_{0,\alpha w}\|^q\le C_*.
	\]
	Next, we choose a basis $v_1,\ldots,v_{M-1}\in\cS$ of
	$\mathbb R^{V_0}/\mathbb R\mathbf1$. Then, for every $u\in\cS$, there exist
	$c_{u,j}\in\mathbb R,\,1\leq j\leq M-1$ such that $u=\sum_jc_{u,j}v_j$. By norm equivalence on
	$\mathbb R^{V_0}/\mathbb R\mathbf1$, there is a constant $C$ such that
	$\sum_j|c_{u,j}|\le C$ for every $u\in\cS$. Hence
	\begin{equation}
		D_{n,\alpha}(u)^{1/p}
		\le\sum_j|c_{u,j}|D_{n,\alpha}(v_j)^{1/p}
		\le C\max_jD_{n,\alpha}(v_j)^{1/p},
	\end{equation}
	where the first inequality uses the fact that $D_{n,\alpha}^{1/p}$ is a norm on $\R^{V_0}/\R{\bf 1}$. Thus 
	\begin{equation}\label{step2:2}
		\|D_{n,\alpha}\|\le C^p\max_jD_{n,\alpha}(v_j).
	\end{equation}
	Therefore,
	\begin{align*}
		\mathbb E\|D_{n,\alpha}\|^q
		\le C^{pq}\mathbb E\left(\max_jD_{n,\alpha}(v_j)\right)^q
		\le C^{pq}C_*(M-1).
	\end{align*}
\end{proof}

\subsection{Moment estimates of \texorpdfstring{$\llbracket D_{n,\alpha}\rrbracket$}{the inverse comparison ratio}}
Next, we prove uniform integrability of $\llbracket D_{n,\alpha}\rrbracket^{1/(p-1)}$ under the condition $\mathbb E\llbracket D_{0,\alpha}\rrbracket^{1/(p-1)}<\infty$.

Throughout the paper, we use the following notation for the distribution of the energy of harmonic extensions of $u\in\mathcal{S}$. 

\begin{definition}\label{def:energy-distribution}
For every $u\in\mathcal{S}$, define 
\[
a_w(u):=\rho^nE_0(h_u\circ\psi_w)\quad\hbox{ for }n\geq 0\hbox{ and }w\in I^n,
\]
where $h_u$ is the harmonic extension of $u$ with respect to the form $f\mapsto\rho^n\sum_{w\in I^n}E_0(f\circ\psi_w)$. 
\end{definition}

Note that $\sum_{w\in I^n}a_w(u)=1$ for $n\geq 0$ and $\sum_{i\in I} a_{wi}(u)=a_w(u)$ for $w\in I^*$, as $E_0$ is the trace of $f\mapsto\rho^n\sum_{w\in I^n}E_0(f\circ\psi_w)$ and $\cR(E_0)=E_0$.

\begin{lemma}\label{lem:3.3}
For any $n\geq 0$, $\alpha\in\mathcal{W}$, $u\in\mathcal{S}$ and $v\in \R^{V_0}$,
\begin{equation}\label{eq:current-lower-bound}
  |E_0(u;v)|\le 
  \left(\sum_{w\in I^n}a_w(u)\llbracket D_{0,\alpha w}\rrbracket^{1/(p-1)}\right)^{(p-1)/p}
  D_{n,\alpha}(v)^{1/p}.
\end{equation}
\end{lemma}
\begin{proof}
For short, we write 
\begin{align*}
E^{(n)}(g):=\rho^n\sum_{w\in I^n}E_0(g\circ\psi_w)\quad\hbox{ and }\quad 
D^{(n)}(g):=\rho^n\sum_{w\in I^n}D_{0,\alpha w}(g\circ\psi_w)\quad\hbox{ for }g\in \R^{V_n}.  
\end{align*}
Also, recall that $E_0=\Tr_{V_0}E^{(n)}$ and $D_{n,\alpha}=\Tr_{V_0}D^{(n)}$.
Let $h_u\in \R^{V_n}$ be the harmonic extension of $u$ with respect to $E^{(n)}$, and let $f_v\in \R^{V_n}$ be the harmonic extension of $v$ with respect to $D^{(n)}$. Then 
\begin{align*}
 |E_0(u;v)|&=\big|E^{(n)}(h_u;f_v)\big|\\
 &=\big|\rho^n\sum_{w\in I^n}
 E_0(h_u\circ\psi_w;f_v\circ\psi_w)\big|\\ 
 &\leq \rho^n\sum_{w\in I^n}
 E_0(h_u\circ\psi_w)^{(p-1)/p}E_0(f_v\circ\psi_w)^{1/p}\\
 &\le\rho^n\sum_{w\in I^n}
 E_0(h_u\circ\psi_w)^{(p-1)/p}\llbracket D_{0,\alpha w}\rrbracket^{1/p}
 D_{0,\alpha w}(f_v\circ\psi_w)^{1/p}\\
 &\le\left(\rho^n\sum_{w\in I^n}
 E_0(h_u\circ\psi_w)\llbracket D_{0,\alpha w}\rrbracket^{1/(p-1)}\right)^{(p-1)/p}
 \left(\rho^n\sum_{w\in I^n}D_{0,\alpha w}(f_v\circ\psi_w)\right)^{1/p}\\
 &= \left(\sum_{w\in I^n}a_w(u)\llbracket D_{0,\alpha w}\rrbracket^{1/(p-1)}\right)^{(p-1)/p}
 D_{n,\alpha}(v)^{1/p}.
\end{align*} 
The first equality follows from \cite[Theorem~6.13, (6.10)]{KajinoShimizu2026}. The second equality follows from the definition of $E^{(n)}$. The first inequality follows from Proposition~\ref{prop:derivative}, and the last follows from H\"older's inequality.
\end{proof}

\begin{lemma}\label{lem:3.5}
There are $u_1,\cdots,u_{M-1}\in\mathcal S$ and a positive constant $C$, which depends only on $E_0,K,p$, such that
\[
\llbracket D_{n,\alpha}\rrbracket^{1/(p-1)}\le C\max_{1\le i\le M-1}\sum_{w\in I^n}a_w(u_i)\llbracket D_{0,\alpha w}\rrbracket^{1/(p-1)}\quad\hbox{ for }\alpha\in\mathcal{W}\hbox{ and }n\geq 0.
\]
\end{lemma}
\begin{proof}
We can choose $u_1,\ldots,u_{M-1}\in\cS$ such that the linear functionals
$E_0(u_i;\cdot)$ form a basis of
$(\R^{V_0}/\R\mathbf1)^*$, where $(\R^{V_0}/\R\mathbf1)^*$ is the dual space of $\R^{V_0}/\R\mathbf1$. Such a choice is possible: otherwise there
would be a non-zero $v\in \R^{V_0}/\R\mathbf1$ so that $E_0(u;v)=0$ for every $u\in \R^{V_0}/\R{\bf 1}$, and in particular $E_0(v;v)=0$, contradicting (RF1). 

By finite-dimensional norm equivalence and Lemma~\ref{lem:3.3}, for every $v\in \R^{V_0}$, 
\begin{align*}
	E_0(v)^{1/p}
	&\le C_1\max_{1\le i\le M-1}|E_0(u_i;v)|\\
	&\le C_1\left(\max_{1\le i\le M-1}
	\sum_{w\in I^n}a_w(u_i)\llbracket D_{0,\alpha w}\rrbracket^{1/(p-1)}\right)^{(p-1)/p}
	D_{n,\alpha}(v)^{1/p}.
\end{align*}
That implies 
\begin{equation}\label{eq:ch3-3-18}
\llbracket D_{n,\alpha}\rrbracket^{1/(p-1)}
\le C\max_{1\le i\le M-1}
\sum_{w\in I^n}a_w(u_i)\llbracket D_{0,\alpha w}\rrbracket^{1/(p-1)},
\end{equation}
where $C=C_1^{p/(p-1)}$. 
\end{proof}

As a consequence of Lemma \ref{lem:3.5}, we have the following moment estimate of $\llbracket D_{n,\alpha}\rrbracket$. 
\begin{lemma}\label{lem:lower-bound-moment}
Suppose that $\mathbb E\llbracket D_{0,\alpha}\rrbracket^{1/(p-1)}<\infty$.
Then
\begin{equation}\label{eq:strong-uniform-inverse}
 \sup_n\mathbb E\llbracket D_{n,\alpha}\rrbracket^{1/(p-1)}<\infty. 
\end{equation}
Moreover, the family $\{\llbracket D_{n,\alpha}\rrbracket^{1/(p-1)}:n\geq0\}$ is uniformly integrable.
\end{lemma}
\begin{proof}
By taking expectations in Lemma~\ref{lem:3.5} and using
$\sum_{w\in I^n}a_w(u_i)=1$, we have
\[
\mathbb E\llbracket D_{n,\alpha}\rrbracket^{1/(p-1)}
\le C(M-1)\mathbb E\llbracket D_{0,\alpha}\rrbracket^{1/(p-1)}.
\]

Next, we prove uniform integrability. By splitting each term  at $R>0$ in Lemma \ref{lem:3.5},
\begin{align*}
\llbracket D_{n,\alpha}\rrbracket^{1/(p-1)}
\le CR+C\max_{1\le i\le M-1}\sum_{w\in I^n}a_w(u_i)
\llbracket D_{0,\alpha w}\rrbracket^{1/(p-1)}
\mathbf1_{\{\llbracket D_{0,\alpha w}\rrbracket^{1/(p-1)}>R\}}.
\end{align*}
On $\{\llbracket D_{n,\alpha}\rrbracket^{1/(p-1)}>2CR\}$, this implies 
\begin{align*}
\llbracket D_{n,\alpha}\rrbracket^{1/(p-1)}
&\le 2\big(\llbracket D_{n,\alpha}\rrbracket^{1/(p-1)}-CR\big)\\
&\le 2C\sum_{i=1}^{M-1}\sum_{w\in I^n}a_w(u_i)
\llbracket D_{0,\alpha w}\rrbracket^{1/(p-1)}
\mathbf1_{\{\llbracket D_{0,\alpha w}\rrbracket^{1/(p-1)}>R\}},
\end{align*}
and hence 
\begin{align*}
&\sup_n\mathbb E\!\left[
\llbracket D_{n,\alpha}\rrbracket^{1/(p-1)}
\mathbf1_{\{\llbracket D_{n,\alpha}\rrbracket^{1/(p-1)}>2CR\}}\right]\\
&\qquad\le2C(M-1)\mathbb E\!\left[
\llbracket D_{0,\alpha}\rrbracket^{1/(p-1)}
\mathbf1_{\{\llbracket D_{0,\alpha}\rrbracket^{1/(p-1)}>R\}}\right]
\longrightarrow0\qquad \hbox{ as }R\to\infty.
\end{align*}
\end{proof}

\section{Concentration estimates}\label{sec:concentration}
In this section, we establish the concentration estimates needed for stochastic homogenization.

\subsection{Energy distribution estimates} We begin with the following observation concerning the proportion of energy in one cell for harmonic extensions, showing that energy cannot concentrate on a single cell.

\begin{lemma}\label{lem:harmonic-energy-propotion}
	Let $E_1,\ldots,E_N\in\cQ_p(V_0)$, $j\in I$, and
	$\beta_1,\beta_2\in(0,\infty)$. Suppose that $E_j\ge\beta_1E_0$ and
	$E_i\le\beta_2E_0$ for $i\in I\setminus\{j\}$. Then for every
	nonconstant $u\in\mathbb R^{V_0}$,
	\[
	\frac{E_j(h\circ\psi_j)}{\sum_{i\in I}E_i(h\circ\psi_i)}
	\le1-\left(\frac12\wedge
	c\left(\frac{\beta_1}{\beta_2}\right)^{1/(p-1)}\right),
	\]
	where $h$ is the harmonic extension of $u$ with respect to $f\mapsto \sum_{i=1}^N E_i(f\circ \psi_i)$ and $c>0$ depends only on $K,E_0,p$.
\end{lemma}

\begin{proof}
	Set $h_i=h\circ\psi_i$ and $e_i=E_i(h_i)$ for $i\in I$. Set
	$e=\sum_{i=1}^Ne_i$. We have $e>0$ since $u$ is not constant. If $e_j=0$,
	there is nothing to prove. Without loss of generality, we can assume that
	$e_j=1$. Then for any $i\ne j$,
	\[
	\frac{e_j}{e}\le\frac1{1+e_i}
	=1-\frac{e_i}{1+e_i}
	\le1-\left(\frac12\wedge\frac{e_i}2\right).
	\]
	So we just need to bound $e_i$ from below for some $i\in I\setminus\{j\}$.
	Since $e_j=1$, the inequality $E_j\ge\beta_1E_0$ and \eqref{e:2.18} give
	\[
	|h_j(x)-h_j(y)|\le C_1E_0(h_j)^{1/p}
	\le C_1\beta_1^{-1/p},\qquad x,y\in V_0.
	\]
	Recall that $\#(\psi_j(V_0)\cap V_0)\le1$. If this intersection is
	nonempty, let $x_0$ be its unique element; otherwise choose
	$x_0\in\psi_j(V_0)$ arbitrarily. Then, by~\eqref{eq:current-euler},
	\begin{align*}
		1=e_j
		&=\sum_{y\in\psi_j(V_0)}J_{E_j}(h_j,\psi_j^{-1}y)
		\bigl(h(y)-h(x_0)\bigr)\\
		&\le(M-1)C_1\beta_1^{-1/p}
		\max_{y\in\psi_j(V_0)\setminus V_0}
		|J_{E_j}(h_j,\psi_j^{-1}y)|.
	\end{align*}
	Thus there exists $x\in\psi_j(V_0)\setminus V_0$ such that
	\[
	|J_{E_j}(h_j,\psi_j^{-1}x)|
	\ge\beta_1^{1/p}(M-1)^{-1}C_1^{-1}.
	\]
	Let $\Gamma:=\max_{z\in V_1\setminus V_0}\#\{i\in I:z\in\psi_i(V_0)\}<\infty$.
Kirchhoff's law gives an index $i\in I\setminus\{j\}$ with
	$x\in\psi_i(V_0)$ such that
	\[
	|J_{E_i}(h_i,\psi_i^{-1}x)|
	\ge\frac{|J_{E_j}(h_j,\psi_j^{-1}x)|}{\Gamma-1}
	\ge c_1\beta_1^{1/p},
	\]
	where $c_1=(C_1(M-1)(\Gamma-1))^{-1}$. Finally, the energy inequality
	\eqref{eq:euler-holder} yields
	\[
	c_1\beta_1^{1/p}
	\le |J_{E_i}(h_i,\psi_i^{-1}x)|
	\le e_i^{1/p'}E_i(\mathbf1_{\{\psi_i^{-1}x\}})^{1/p}
	\le C\beta_2^{1/p}e_i^{1/p'},
	\qquad\frac1{p'}+\frac1p=1.
	\]
	The last inequality follows from $E_i\le\beta_2E_0$ and \eqref{e:2.18}.
	Thus
	\[
	e_i\ge c_2\left(\frac{\beta_1}{\beta_2}\right)^{1/(p-1)},
	\qquad c_2=\left(\frac{c_1}{C}\right)^{p'}.
	\]
	The lemma holds with $c=c_2/2$.
\end{proof}

As an immediate consequence of Lemma \ref{lem:harmonic-energy-propotion}, we see an exponential upper bound on the proportion of energy on an $n$-cell for the standard energy.
\begin{corollary}\label{coro:4.2}
For $u\in \mathcal{S}$, recall that we define 
\[
a_w(u)=\rho^nE_0(h_u\circ\psi_w)\quad\hbox{ for }w\in I^n,n\geq 0,
\]
where $h_u$ is the harmonic extension of $u$ with respect to $f\mapsto\sum_{w\in I^n}E_0(f\circ \psi_w)$. Then, there is $\theta\in(0,1)$ depending only on $K,E_0,p$ such that 
\[
a_{wi}(u)\leq\theta a_w(u)\ \hbox{ for any }u\in\mathcal S,\, w\in I^*\hbox{ and }i\in I.
\]
\end{corollary}

We next apply Lemma \ref{lem:harmonic-energy-propotion} to the random environment. Fix $\alpha\in\cW$ and $n\ge1$. Recall that $D_{n,\alpha}$ is the trace
onto $V_0$ of the resistance form
\begin{equation}\label{eq:ch3-3-23}
 f\longmapsto\rho^n\sum_{w\in I^n}D_{0,\alpha w}(f\circ\psi_w),
 \qquad f\in\R^{V_n}.
\end{equation}
Thus, for every $u\in\R^{V_0}\setminus\R{\bf 1}$,
\[
 D_{n,\alpha}(u)=\rho^n\sum_{w\in I^n}
                         D_{0,\alpha w}(h\circ\psi_w),
\]
where $h$ is the harmonic extension of $u$ to $V_n$ with respect to
\eqref{eq:ch3-3-23}. Define
\[
 A_w(u)=\frac{\rho^nD_{0,\alpha w}(h\circ\psi_w)}{D_{n,\alpha}(u)},
 \qquad w\in I^n.
\]

\begin{proposition}\label{prop:4.3}	
Assume that there are $a>0$ and $C_*>0$ such that 
\[
D_{0,\alpha}\geq aE_0\hbox{ almost surely  and }\quad\mathbb E\|D_{0,\alpha}\|^2\le C_*.
\]
Let $t>1$. There are constants $C_A\in(0,\infty)$ and
$\lambda\in(0,1)$, depending on $a,K,E_0,p,t$ and $C_*$, such that
\[
 \sup_{u\in\cS}\mathbb E\left[\sum_{|w|=n}A_w(u)^t\right]
 \le C_A\lambda^n\quad\hbox{ for }n\ge1.
\]
\end{proposition}
\begin{proof}
Fix $u\in\cS$ and let $h$ be its harmonic extension to $V_n$ with respect to the form in~\eqref{eq:ch3-3-23}.
Let $e_w=\rho^nD_{0,\alpha w}(h\circ\psi_w)$ for $w\in I^n$, and define
\[
 e_\tau=\sum_{b\in I^{n-m}}e_{\tau b},
 \qquad 0\le m\le n-1,\quad\tau\in I^m.
\]
Then $e_\tau=\sum_{i\in I}e_{\tau i}$. By the definition of traces,
for $|\tau|=m$,
\[
 e_{\tau i}=\rho^{m+1}D_{n-m-1,\alpha\tau i}(h\circ\psi_{\tau i}).
\]
Thus, when $e_\tau>0$, the ratios $e_{\tau i}/e_\tau$ are first-level
cell energy proportions for the effective forms
$D_{n-m-1,\alpha\tau i}$.

Let $L>0$ be chosen later, and let
$\delta=\frac12\wedge c(a/L)^{1/(p-1)}$, with $c$ from Lemma~\ref{lem:harmonic-energy-propotion}.
For a word $w=w_1\cdots w_n\in I^n$, define the event
\[
 G_m(w)=\left\{\|D_{n-m-1,\alpha w_1\cdots w_mi}\|\le L\text{ for each }i\in I\setminus\{w_{m+1}\}\right\}.
\]
On this event, Lemma~\ref{lem:harmonic-energy-propotion} gives
\[
 \frac{e_{w_1\cdots w_mw_{m+1}}}{e_{w_1\cdots w_m}}\le1-\delta
\]
whenever the denominator is positive. Notice that
\[
 A_w(u)=\frac{e_w}{e_\emptyset}
       =\prod_{m=0}^{n-1}
                \frac{e_{w_1\cdots w_mw_{m+1}}}{e_{w_1\cdots w_m}},
\]
with $A_w(u)=0$ when some $e_{w_1\cdots w_m}$ is zero. Thus
\[
 \max_{w\in I^n}A_w(u)\le(1-\delta)^{n/2}
 \quad\text{on }\left\{\sum_{m=0}^{n-1}\mathbf1_{G_m(w)}\ge\frac n2
                         \text{ for every }w\in I^n\right\}.
\]

We next show that this event has high probability.
Markov's inequality and Lemma~\ref{lem:moment2} imply
\[
\mathbb P(\Omega\setminus G_m(w))\le
(N-1)L^{-2}\sup_k\mathbb E\|D_{k,\alpha}\|^2
 =:\varepsilon_L.
\]
For a fixed word $w\in I^n$, the events $G_m(w)$, $0\le m\le n-1$,
are independent.
Choose $L$ so large that $2N\varepsilon_L^{1/2}\le\frac12$. Hence
\[
 \mathbb P\left(\sum_{m=0}^{n-1}\mathbf1_{G_m(w)}<\frac n2\right)
 \le\sum_{j>n/2}\binom nj\varepsilon_L^j
 \le2^n\varepsilon_L^{n/2}.
\]
Consequently,
\[
 \mathbb P\left(\bigcup_{w\in I^n}
 \left\{\sum_{m=0}^{n-1}\mathbf1_{G_m(w)}<\frac n2\right\}\right)
 \le N^n\,2^n\varepsilon_L^{n/2}
 =\bigl(2N\varepsilon_L^{1/2}\bigr)^n\le 2^{-n},
\]
and hence
\[
 \mathbb P\left(\max_{w\in I^n}A_w(u)\le(1-\delta)^{n/2}\right)
 \ge1-2^{-n}.
\]

Finally, as
\[
 \sum_{w\in I^n}A_w(u)^t
 \le\left(\max_{w\in I^n}A_w(u)\right)^{t-1}
          \sum_{w\in I^n}A_w(u)
 =\left(\max_{w\in I^n}A_w(u)\right)^{t-1}\le1,
\]
we have
\begin{align*}
\mathbb E[\sum_{w\in I^n}A_w(u)^t]
	&\leq\mathbb E[\sum_{w\in I^n}A_w(u)^t,\max_{w\in I^n}A_w(u)\leq (1-\delta)^{n/2}] +\mathbb P(\max_{w\in I^n}A_w(u)>(1-\delta)^{n/2})\\
	&\leq (1-\delta)^{(t-1)n/2}+2^{-n}. 
\end{align*}
Thus we can take $C_A=2$ and
$\lambda=\max\{(1-\delta)^{(t-1)/2},\frac12\}$.
\end{proof}

\subsection{Concentration: uniform lower bound case} We first assume that there are $a>0$ and $C_*>0$ such that
\begin{equation}\label{eq:concentration-assumption}
	D_{0,\alpha}\geq aE_0\hbox{ almost surely  and }\quad\mathbb E\|D_{0,\alpha}\|^2\le C_*.
\end{equation}

\begin{theorem}\label{thm:concentration}
Assume \eqref{eq:concentration-assumption}. Then, there are $C_S<\infty$ and $\gamma\in(0,1)$ such that
	\begin{equation}
		\mathbb E\|D_{n,\alpha}-\bar E_n\|\le C_S\gamma^n\quad\hbox{ for }n\ge0,
	\end{equation}
	where $\bar E_n:=\mathbb ED_{n,\alpha}$.
\end{theorem}

Note that by Lemma~\ref{lem:moment2}, \eqref{eq:concentration-assumption} implies $\sup_n\mathbb E\|D_{n,\alpha}\|^2<\infty$. Hence, $\bar E_n:=\mathbb ED_{n,\alpha}$ is well defined and is a $p$-resistance form by Lemma~\ref{lem:state-closure}.\medskip

We will use the following Efron--Stein
inequality (see \cite{EfronStein1981,BoucheronLugosiMassart2013}).

\begin{proposition}[Efron--Stein inequality]\label{Efron--Stein}
Let $\xi_1,\ldots,\xi_n$ be independent random variables, and let
$Y=f(\xi_1,\ldots,\xi_n)\in L^2$. Let $Y^{(i)}$ be obtained by replacing
$\xi_i$ with an independent copy $\xi_i'$, that is,
$Y^{(i)}=f(\xi_1,\ldots,\xi_i',\ldots,\xi_n)$. Then
\[
 \operatorname{Var}(Y)
 \le\frac12\sum_{i=1}^n\mathbb E|Y-Y^{(i)}|^2.
\]
\end{proposition}

\begin{proposition}\label{prop:4.5}
Assume \eqref{eq:concentration-assumption}. 
There exist $C_V<\infty$ and $r_0\in(0,1)$ such that
\begin{equation}\label{eq:ch3-3-24}
 \sup_{u\in\cS}\operatorname{Var}\!\left(\sqrt{D_{n,\alpha}(u)}\right)
 \le C_Vr_0^n.
\end{equation}
\end{proposition}
\begin{proof}
Fix $u\in\cS$, and write $Y=D_{n,\alpha}(u)$ and $A_w=A_w(u)$ for $w\in I^n$. Let $h$ be the harmonic
extension of $u$ to $V_n$ with respect to \eqref{eq:ch3-3-23}.
For $w\in I^n$, replace $D_{0,\alpha w}$ by an independent copy
$D_{0,\alpha w}^{(w)}$, leaving the other forms unchanged. Denote the
resulting versions of $Y$ and $A_w$ by $Y^{(w)}$ and $A_w^{(w)}$,
respectively. Put
\[
 \chi=\frac{\|D_{0,\alpha w}\|}{a}-1\ge0,
 \qquad
 \chi^{(w)}=\frac{\|D_{0,\alpha w}^{(w)}\|}{a}-1\ge0.
\]
Then
\[
 D_{0,\alpha w}^{(w)}-D_{0,\alpha w}
 \le\chi^{(w)}D_{0,\alpha w},
 \qquad
 D_{0,\alpha w}-D_{0,\alpha w}^{(w)}
 \le\chi D_{0,\alpha w}^{(w)}.
\]
By the definition of traces,
\[
 Y^{(w)}-Y\le\rho^n
       (D_{0,\alpha w}^{(w)}-D_{0,\alpha w})(h\circ\psi_w)\le\chi^{(w)}A_wY.
\]
On $\{Y^{(w)}\ge Y\}$, it follows that 
\[
 \left|\sqrt{Y^{(w)}}-\sqrt Y\right|^2
 =\frac{(Y^{(w)}-Y)^2}{(\sqrt{Y^{(w)}}+\sqrt Y)^2}
 \le\frac14(\chi^{(w)})^2A_w^2Y.
\]
For the same reason, on $\{Y>Y^{(w)}\}$,
\[
 \left|\sqrt Y-\sqrt{Y^{(w)}}\right|^2
 \le\frac14\chi^2(A_w^{(w)})^2Y^{(w)}.
\]
Since $\chi^{(w)}$ is independent of $(A_w,Y)$, and $\chi$ is independent
of $(A_w^{(w)},Y^{(w)})$, we get
\[
 \mathbb E\left|\sqrt{Y^{(w)}}-\sqrt Y\right|^2
 \le\frac12(\mathbb E\chi^2)(\mathbb EA_w^2Y)
 \le\frac{C_*}{2a^2}\mathbb E(A_w^2Y).
\]
Applying Propositions~\ref{prop:4.3} and \ref{Efron--Stein}, we see 
\[
\begin{aligned}
 \operatorname{Var}(\sqrt Y)
 &\le \frac12\sum_{w\in I^n}
 \mathbb E\left|\sqrt{Y^{(w)}}-\sqrt Y\right|^2\\
 &\le \frac{C_*}{4a^2}
 \sum_{w\in I^n}\mathbb E(A_w^2Y)
 \mathrel=C\mathbb E\left[Y\sum_{w\in I^n}A_w^2\right]\\
 &\le C(\mathbb EY^2)^{1/2}
          \left(\mathbb E\left(\sum_{w\in I^n}A_w^2\right)^2\right)^{1/2}\\
 &\le CC_*^{1/2}
          \left(\mathbb E\sum_{w\in I^n}A_w^2\right)^{1/2}
 \le CC_*^{1/2}C_A^{1/2}\lambda^{n/2}.
\end{aligned}
\]
\end{proof}

\begin{corollary}\label{coro:4.6}
Assume \eqref{eq:concentration-assumption}. 
There is $C<\infty$ such that 
\begin{equation}\label{eq:ch3-3-25}
 \sup_{u\in\cS}\mathbb E|D_{n,\alpha}(u)-\bar E_n(u)|
 \le Cr_0^{n/2},
\end{equation}
where $r_0$ is the same as in Proposition~\ref{prop:4.5}.
\end{corollary}
\begin{proof}
Write $Y=D_{n,\alpha}(u)$ for a fixed $u\in\cS$, as before. The triangle
inequality and the Cauchy--Schwarz inequality give
\[
\begin{aligned}
 \mathbb E|Y-\mathbb EY|
 &\le2\mathbb E|Y-(\mathbb E\sqrt Y)^2|\\
 &\le2\operatorname{Var}(\sqrt Y)^{1/2}
        \left(\mathbb E(\sqrt Y+\mathbb E\sqrt Y)^2\right)^{1/2}\\
 &\le4\bigl(\operatorname{Var}(\sqrt Y)\mathbb EY\bigr)^{1/2}.
\end{aligned}
\]
The conclusion follows from
$\mathbb EY=\mathbb ED_{n,\alpha}(u)\le(\mathbb EY^2)^{1/2}
\le C_*^{1/2}$ and Proposition~\ref{prop:4.5}.
\end{proof}

Finally, we deduce the norm estimate in Theorem \ref{thm:concentration} from Corollary \ref{coro:4.6}, by using the continuity of the forms. 
For $E\in\cQ_p(V_0)$, define
\[
\operatorname{Lip}_{\cS}(E)
:=\sup_{\substack{u,v\in\cS\\u\ne v}}
\frac{|E(u)-E(v)|}{E_0(u-v)^{1/p}}.
\]

\begin{lemma}\label{lem:Lipschitz}
	For every $E\in\cQ_p(V_0)$,
	$\operatorname{Lip}_{\cS}(E)\le p\|E\|$.
\end{lemma}

\begin{proof}
	Using $|x^p-y^p|\le p\max\{x^{p-1},y^{p-1}\}|x-y|$ for $x,y\ge0$ and
	the triangle inequality for $E^{1/p}$, we have, for $u,v\in\cS$,
	\begin{align*}
		|E(u)-E(v)|
		&\le p\max\{E(u)^{1-1/p},E(v)^{1-1/p}\}
		|E(u)^{1/p}-E(v)^{1/p}|\\
		&\le p\|E\|^{1-1/p}E(u-v)^{1/p}
		\le p\|E\|E_0(u-v)^{1/p}.
	\end{align*}
\end{proof}

\begin{proof}[Proof of Theorem \ref{thm:concentration}]
For $0<\delta<1$, $\cS$ has a $\delta$-net $\mathcal N_\delta$ (i.e., for every $u\in\cS$, there exists $v\in\mathcal N_\delta$ such that $E_0(u-v)^{1/p}<\delta$) with at
most $C\delta^{-(M-2)}$ points, where $C$ depends only on the choice of
$E_0$. Corollary~\ref{coro:4.6} then gives
\[
 \mathbb E\max_{v\in\mathcal N_\delta}
             |D_{n,\alpha}(v)-\bar E_n(v)|
 \lesssim\delta^{-(M-2)}r_0^{n/2}.
\]
For $u\in\cS$, choose $v\in\mathcal N_\delta$ such that
$E_0(u-v)^{1/p}<\delta$. By Lemma~\ref{lem:Lipschitz},
\[
\begin{aligned}
 |(D_{n,\alpha}-\bar E_n)(u)-(D_{n,\alpha}-\bar E_n)(v)|
 &\le\delta\bigl(\operatorname{Lip}_{\cS}(D_{n,\alpha})
                 +\operatorname{Lip}_{\cS}(\bar E_n)\bigr)\\
 &\le\delta p(\|D_{n,\alpha}\|+\|\bar E_n\|).
\end{aligned}
\]
Therefore,
\[
 \mathbb E\|D_{n,\alpha}-\bar E_n\|
 \lesssim\delta^{-(M-2)}r_0^{n/2}+\delta.
\]
Taking $\delta=(r_0^{1/2})^{n/(M-1)}$ proves the assertion.
\end{proof}

\subsection{Concentration: a moment condition for the lower bound}
Next, we will replace the uniform lower bound condition $D_{0,\alpha}\geq aE_0$ by a moment condition on $\llbracket D_{0,\alpha}\rrbracket$.\smallskip 

First, we make a rough estimate of the $L^1$ concentration, without assuming the uniform lower bound. Note that we may replace $D_{0,\alpha}$ with any other random environment. 

\begin{lemma}\label{lem:4.9}
	Assume that $\mathbb E\|D_{0,\alpha}\|^2<C_*$. 
	Then, for every $a>0$, there are positive constants $C_1,C_2$ and $\gamma\in(0,1)$ that depend only on $C_*,E_0,K,p,a$ so that
	\begin{equation}
		\mathbb E\left\|D_{k,\alpha}-\mathbb{E}D_{k,\alpha}\right\|\le C_1\gamma^k+C_2\sqrt{N^k\mathbb{P}(D_{0,\alpha}\not\geq aE_0)}.
	\end{equation}
\end{lemma}
\begin{proof}
	Let	
	\[
	\widetilde D_{0,\beta}
	:=D_{0,\beta}+aE_0\cdot\mathbf1_{\{D_{0,\beta}\not\ge aE_0\}}\quad\hbox{ for }\beta\in \mathcal{W},
	\]
	and 
	\[
	\widetilde D_{k,\alpha}:=\cR^k((\widetilde D_{0,\alpha w})_{w\in I^k})\quad\hbox{ for }k\geq 0.
	\]
	Then
	\begin{align}
		\mathbb E\big\|D_{k,\alpha}-\mathbb{E}D_{k,\alpha}\big\|
		&\leq \mathbb E\big\|D_{k,\alpha}-\widetilde D_{k,\alpha}\big\|
		+\mathbb E\big\|\widetilde D_{k,\alpha}-\mathbb E\widetilde D_{k,\alpha}\big\|
		+\big\|\mathbb E\widetilde D_{k,\alpha}-\mathbb ED_{k,\alpha}\big\|\notag\\
		&\leq 2\mathbb E\big\|D_{k,\alpha}-\widetilde D_{k,\alpha}\big\|
		+\mathbb E\big\|\widetilde D_{k,\alpha}-\mathbb E\widetilde D_{k,\alpha}\big\|.\label{e:4.5}
	\end{align}
	It remains to establish upper bounds for the two terms in \eqref{e:4.5}. 
	
	For the second term of \eqref{e:4.5}, by Theorem~\ref{thm:concentration}, there are $C_1>0$ and $\gamma\in(0,1)$ so that 
	\begin{equation*}
		\mathbb E\left\|\widetilde D_{k,\alpha}
		-\mathbb E\widetilde D_{k,\alpha}\right\|\leq C_1\gamma^k.
	\end{equation*}
	
	For the first term of \eqref{e:4.5}, we note that 
	\[
	\mathbb{P}\left(D_{k,\alpha}	\neq \widetilde D_{k,\alpha}\right)\leq \sum_{w\in I^k}\mathbb{P}(D_{0,\alpha w}\not\ge aE_0)=N^k\mathbb{P}(D_{0,\alpha}\not\geq aE_0),
	\]
	which implies, by using the Cauchy--Schwarz inequality and the second-moment assumption,
	\begin{align*}
		\quad\ \mathbb{E}\big\|D_{k,\alpha}-\widetilde D_{k,\alpha}\big\|
		=\mathbb{E}\left[\big\|D_{k,\alpha}-\widetilde D_{k,\alpha}\big\|;\widetilde D_{k,\alpha}\neq D_{k,\alpha}\right]
		\leq C_2\sqrt{N^k\mathbb{P}(D_{0,\alpha}\not\geq aE_0)}.
	\end{align*}
\end{proof}

We next show that the effective forms are bounded below with probability tending to one. The following lemma uses a law of large numbers from Appendix~\ref{AppendixA}.

\begin{lemma}\label{lem:4.10}
Assume $\mathbb{E}\llbracket D_{0,\alpha}\rrbracket^{1/(p-1)}<\infty$. Then there is $a>0$ so that
\[
\lim_{n\to\infty}\mathbb P(D_{n,\alpha}\not\ge aE_0)=0.
\] 
\end{lemma}
\begin{proof}
Let $Y_w:=\llbracket D_{0,\alpha w}\rrbracket^{1/(p-1)}$ for $w\in I^*$. Then, Lemma \ref{lem:3.5} states that there are $u_1,u_2,\cdots,u_{M-1}\in\mathcal{S}$ such that 
\[
\llbracket D_{n,\alpha}\rrbracket^{1/(p-1)}\leq Z_n:=C_1\max_{1\le i\le M-1}\sum_{w\in I^n}a_w(u_i)Y_w.  
\]

The assumption can be rephrased as $\mathbb{E}Y_w<\infty$. Proposition \ref{lem:tree-averages} and Corollary \ref{coro:4.2} imply that
\[
\lim\limits_{n\to\infty}Z_n=C_1\mathbb{E}Y_w\quad\hbox{almost surely}.
\]
Hence, for every $b>C_1\mathbb{E}Y_w$, 
\[
\lim\limits_{n\to\infty}\mathbb P(D_{n,\alpha}\not\geq b^{-(p-1)} E_0)=\lim\limits_{n\to\infty}\mathbb P(\llbracket D_{n,\alpha}\rrbracket^{1/(p-1)}> b)\leq \lim\limits_{n\to\infty}\mathbb P(Z_n>b)=0.
\]

\end{proof} 

\begin{theorem}\label{thm:concernch4}
Assume that
\[
\mathbb E\llbracket D_{0,\alpha}\rrbracket^{1/(p-1)}<\infty
\quad\hbox{and}\quad
\mathbb E\|D_{0,\alpha}\|^q<\infty
\quad\hbox{for some }q>0.
\]
Then
\begin{equation}\label{eq:concen1}
\lim_{n\to\infty}
\mathbb E\|D_{n,\alpha}-\mathbb ED_{n,\alpha}\|=0.
\end{equation}
\end{theorem}

\begin{proof}
By Lemmas~\ref{lem:moment-upgrade} and~\ref{lem:moment2}, there exists
$j_0\ge0$ such that
\[
\sup_{j\ge j_0}\mathbb E\|D_{j,\alpha}\|^2<\infty.
\]
Choose $a>0$ as in Lemma~\ref{lem:4.10}, and let
\[
P_j:=\mathbb P(D_{j,\alpha}\not\ge aE_0)
\quad\hbox{for }j\ge0.
\]
Applying Lemma~\ref{lem:4.9} with $D_{j,\alpha}$ in place of
$D_{0,\alpha}$ gives
\begin{equation}\label{eq:4.7}
\mathbb E\|D_{j+k,\alpha}-\mathbb ED_{j+k,\alpha}\|
\le C_1\gamma^k+C_2(N^kP_j)^{1/2}
\quad\hbox{for }j\ge j_0,\ k\ge0.
\end{equation}
By Lemma~\ref{lem:4.10}, $P_j\to0$ as $j\to\infty$.
Thus, for every fixed $k\ge0$, letting $j\to\infty$ in
\eqref{eq:4.7} yields
\[
\limsup_{n\to\infty}
\mathbb E\|D_{n,\alpha}-\mathbb ED_{n,\alpha}\|
\le C_1\gamma^k.
\]
Letting $k\to\infty$ proves \eqref{eq:concen1}.
\end{proof}

\section{Stochastic homogenization}\label{sec:homogenization}
{In this section, we prove Theorem~\ref{thm:random-homogenization}.}
We will use the following deterministic result, which will
be proved in Appendix~\ref{sec:peirone-proof}.

\begin{theorem}[Peirone-type convergence]\label{thm:peirone}
	For every $E\in\cQ_p(V_0)$, there exists $E_\infty\in\cQ_p(V_0)$ such that
	\[
	\lim_{n\to\infty}\|\cR^n(E)-E_\infty\|=0,
	\qquad \cR(E_\infty)=E_\infty.
	\]
\end{theorem}

We divide the proof into two steps. Readers may skip the details of the first part on a first reading.

\subsection{\texorpdfstring{$L^2$}{L2}-bounded modification of \texorpdfstring{$D_{n,\alpha}$}{Dn}}\label{sec5.1}
We do not always assume a moment condition on $\|D_{0,\alpha}\|$ in Theorem \ref{thm:random-homogenization}. To apply the concentration theorem,  Theorem~\ref{thm:concentration}, we construct a modification
$\widehat D_{n,\alpha}$ of $D_{n,\alpha}$ such that
$\widehat D_{n,\alpha}=D_{n,\alpha}$ with high probability and
$\sup_n\mathbb E\|\widehat D_{n,\alpha}\|^2<\infty$.

\begin{proposition}\label{prop:L^2modification}
There are events
$\mathcal H_{n,\alpha}\in\sigma\{D_{0,\alpha w}:w\in I^n\}$ such that
\begin{equation}\label{eq:ch3-3-27}
 \delta_n:=\mathbb P(\mathcal H_{n,\alpha})
 \lesssim\exp(-cr^{-n}),
\end{equation}
where we recall that $r\in(0,1)$ is the ratio of the self-similar contractions defining $K$, and 
\begin{equation}\label{eq:ch3-3-29}
	\sup_n\mathbb E\|\widehat D_{n,\alpha}\|^2<\infty,
\end{equation}
where 
\begin{equation}\label{eq:ch3-3-28}
 \widehat D_{n,\alpha}:=
 \begin{cases}
  D_{n,\alpha},&\omega\in\mathcal H_{n,\alpha}^c,\\
  E_0,&\omega\in\mathcal H_{n,\alpha}.
 \end{cases}
\end{equation}
\end{proposition}
\begin{proof}
Set $\mathcal H_{0,\alpha}=\Omega$, and fix $n\geq 1$. For $w,v\in I^n$, we say $w\sim v$ if and
only if $w\ne v$ and $\psi_w(V_0)\cap\psi_v(V_0)\ne\emptyset$. Let $G$ be the graph with vertex set $I^n$ and edge set $\big\{\{w,v\}\subset I^n:\,w\sim v\big\}$.

Let $T>0$ be a large number to be fixed later. We say that $w\in I^n$ is \textit{bad} if $\|D_{0,\alpha w}\|>T$; otherwise it is \textit{good}. For $w\in I^n$, if $w$ is bad, let $S(w)$ be the component of bad words containing $w$; otherwise set $S(w)=\emptyset$. Let 
\[
\mathcal H_{n,\alpha}:=\Big\{\hbox{There are distinct }w,v\in \partial I^n\hbox{ so that }S(w)=S(v)\ne\emptyset\Big\},
\]
where $\partial I^n:=\{w\in I^n:\psi_w(V_0)\cap V_0\neq \emptyset\}$.

First, we show that $\delta_n:=\mathbb P(\mathcal H_{n,\alpha})$ decreases
rapidly if we choose a large enough $T$. The events $\{\|D_{0,\alpha w}\|>T\}$, $w\in I^n$, are independent, with common probability 
\[
\varepsilon:=\mathbb P(\|D_{0,\alpha}\|>T).
\] 
Note that for each $w\in I^n$ and connected subset $F$ of $G$ with $m$ vertices containing $w$, we can find a walk $w=x_0\sim x_1\sim\cdots\sim x_{2m-2}=w$ that goes through each vertex in $F$. Hence, there are at most $\Delta_0^{2m-2}$ connected subsets of $G$ with $m$ vertices containing $w$, where $\Delta_0$ is an upper bound, independent of $n$, for the degrees of vertices in $G$. Therefore,
\[
  \mathbb P(|S(w)|\ge m)\le\Delta_0^{2m-2}\varepsilon^m\le(\Delta_0^2\varepsilon)^m=:\tau^m.
\]
We choose $T$ large enough so that 
\[
\Delta_0\tau^{1/4}<1\,\hbox{ and }\,\tau<1.
\]
Note that any connected subset $F$ containing two distinct vertices of $\partial I^n$ has at least $\lceil c_0r^{-n}\rceil$ elements, where $c_0=(\operatorname{diam}K)^{-1}
\min\{|x-y|:x,y\in V_0,\ x\ne y\}$. Thus,
\begin{align*}
\mathbb P(\mathcal H_{n,\alpha})
&\leq \sum_{w\in \partial I^n}\mathbb P(|S(w)|\ge \lceil c_0r^{-n}\rceil)\\
&\leq M\tau^{\lceil c_0r^{-n}\rceil}.
\end{align*}
This proves \eqref{eq:ch3-3-27}.\medskip 

It remains to estimate $\mathbb E\|\widehat D_{n,\alpha}\|^2$. To achieve this, it suffices to estimate $\mathbb E\left[\|\widehat D_{n,\alpha}\|^2{\bf 1}_{\mathcal{H}_{n,\alpha}^c}\right]$. 

Now, we fix $u\in\R^{V_0}\setminus\R{\bf 1}$ and assume $\mathcal{H}_{n,\alpha}^c$ occurs. We establish an upper bound of $\widehat{D}_{n,\alpha}(u)$. Let $h$ be its harmonic
extension with respect to $f\mapsto\sum_{w\in I^n}E_0(f\circ\psi_w)$.
For each component $S$ of bad cells, if $V_0\cap\bigcup_{w\in S}\psi_w(V_0)\neq\emptyset$, we take $x_S$ to be the unique vertex in $V_0\cap\bigcup_{w\in S}\psi_w(V_0)$; otherwise, we choose an arbitrary vertex $x_S\in \bigcup_{w\in S}\psi_w(V_0)$. Then, we define $h'\in\R^{V_n}$ by
\[
 h'(x)=
\begin{cases}
h(x_S),&x\in\bigcup_{w\in S}\psi_w(V_0)\hbox{ for some component }S\hbox{ of bad cells},\\
h(x),&\text{otherwise}.
\end{cases}
\]
Note that $h'|_{V_0}=u$ as $\mathcal{H}_{n,\alpha}$ does not occur, and $h'$
is constant on each bad cell. For each good word $b$, set
\[
 A_b=\{b\}\cup
 \bigcup_{\substack{S\text{ a bad component}\\S\text{ adjacent to }b}}S.
\]
For $x\in \psi_b(V_0)$, there is $z_x\in \bigcup_{w\in A_b}\psi_w(V_0)$ such that $h'(x)=h(z_x)$. Therefore, 
\[
\operatorname{osc}(h'\circ\psi_b)\le\operatorname{osc}(h|_{\bigcup_{w\in A_b}\psi_w(V_0)})\le\sum_{c\in A_b}\operatorname{osc}(h\circ\psi_c).
\]
By H\"older's inequality,
\[
\operatorname{osc}(h'\circ\psi_b)^p\le |A_b|^{p-1}
 \sum_{c\in A_b}\operatorname{osc}(h\circ\psi_c)^p. 
\]
Then, by \eqref{e:2.18}, 
\[
 E_0(h'\circ\psi_b)\le C|A_b|^{p-1}
 \sum_{c\in A_b}E_0(h\circ\psi_c).
\]
Summing the estimate over all indices, we get 
\begin{align}
 D_{n,\alpha}(u)
 &\le\rho^n\sum_{b\text{ good}}D_{0,\alpha b}(h'\circ\psi_b)\nonumber\\
 &\le T\rho^n\sum_{b\text{ good}}E_0(h'\circ\psi_b)\nonumber\\
 &\le CT\rho^n\sum_{b\in I^n}|A_b|^{p-1}
\sum_{c\in A_b}E_0(h\circ\psi_c),\label{eq:ch3-3-30}
\end{align}
where we set $A_b=\emptyset$ when $b$ is bad. By exchanging the order of summation in \eqref{eq:ch3-3-30}, we conclude 
\begin{equation}\label{eq:ch3-3-31}
D_{n,\alpha}(u){\bf 1}_{\mathcal{H}_{n,\alpha}^c}\leq CT\rho^n\sum_{b,c\in I^n}
E_0(h\circ\psi_c)|A_b|^{p-1}\mathbf1_{\{c\in A_b\}}.
\end{equation}

To establish an upper bound of $\mathbb E\left[\|\widehat D_{n,\alpha}\|^2{\bf 1}_{\mathcal{H}_{n,\alpha}^c}\right]$, we next estimate
\[
 \mathbb E|A_b|^{4(p-1)}
 \quad\text{and}\quad \mathbb P(c\in A_b).
\]
For the first quantity, we have
\[
 |A_b|\le 1+\sum_{v\sim b}|S(v)|.
\]
Therefore,
\[
 |A_b|^{4(p-1)}\lesssim
 1+\sum_{v\sim b}|S(v)|^{4(p-1)}.
\]
Since $\mathbb P(|S(v)|\ge m)\le\tau^m$,
\begin{equation}\label{eqn5.5}
 \mathbb E|A_b|^{4(p-1)}\lesssim
 1+\sum_{v\sim b}\mathbb E|S(v)|^{4(p-1)}\le C<\infty.
\end{equation}
Next, let us estimate $\mathbb P(c\in A_b)$. Write $\ell(b,c)$ for the graph
distance. If $b\ne c$ and $c\in A_b$, then $|S(c)|\ge\ell(b,c)$. Thus,
\begin{equation}\label{eqn5.6}
 \mathbb P(c\in A_b)\le\tau^{\ell(b,c)}.
\end{equation}
Now substitute these estimates into \eqref{eq:ch3-3-31} to obtain
 \begin{align}
\Big(\mathbb E\big[|D_{n,\alpha}(u)|^2{\bf 1}_{\mathcal{H}_{n,\alpha}^c}\big]\Big)^{1/2}&\leq\left(\mathbb E\left[\left(
 CT\rho^n\sum_{b,c\in I^n}E_0(h\circ\psi_c)
 |A_b|^{p-1}\mathbf1_{\{c\in A_b\}}
  \right)^2\right]\right)^{1/2}\nonumber\\
 &\quad\le CT\rho^n\sum_{b,c\in I^n}E_0(h\circ\psi_c)
 \left(\mathbb E\left[|A_b|^{2(p-1)}\mathbf1_{\{c\in A_b\}}\right]\right)^{1/2}\nonumber\\
 &\quad\le C'T\rho^n\sum_{b,c\in I^n}E_0(h\circ\psi_c)\mathbb P(c\in A_b)^{1/4}\nonumber\\
 &\quad\le C''T\rho^n\sum_{c\in I^n}E_0(h\circ\psi_c)\sum_{b\in I^n}\tau^{\ell(b,c)/4},\label{eq:ch3-3-32} 
 \end{align}
where the first inequality uses \eqref{eq:ch3-3-31}, the second uses Minkowski's inequality, the third uses the Cauchy--Schwarz inequality and \eqref{eqn5.5}, and the last uses \eqref{eqn5.6}.
Since
\[
 \begin{aligned}
 \sum_{b\in I^n}\tau^{\ell(b,c)/4}
 &=\sum_{j\ge0}\#\{b\in I^n:\ell(b,c)=j\}\tau^{j/4}\\
 &\le\sum_{j\ge0}(\Delta_0\tau^{1/4})^j
 =\frac1{1-\Delta_0\tau^{1/4}}<\infty,
 \end{aligned}
\]
we can bound \eqref{eq:ch3-3-32} by 
\[
 C''T\rho^n(1-\Delta_0\tau^{1/4})^{-1}
 \sum_{c\in I^n}E_0(h\circ\psi_c)
 =C''T(1-\Delta_0\tau^{1/4})^{-1}E_0(u).
\]
Thus
\[
\mathbb E\widehat D_{n,\alpha}(u)^2
\leq 
E_0(u)^2+\mathbb E\big[|D_{n,\alpha}(u)|^2{\bf1}_{\mathcal{H}_{n,\alpha}^c}\big]
\le
\Big(1+\big(C''T(1-\Delta_0\tau^{1/4})^{-1}\big)^2\Big)E_0(u)^2.
\]
By the same argument as in Lemma~\ref{lem:moment2}, we obtain
\eqref{eq:ch3-3-29}. 
\end{proof}

\begin{lemma}\label{quasi-trace}
	$
	\mathbb E\bigl\|\widehat D_{j+k,\alpha}
	-\cR^k((\widehat D_{j,\alpha w})_{w\in I^k})\bigr\|
	\le C(\delta_{j+k}+N^k\delta_j)^{1/2}
	$
	for each $j,k\geq 0$, where $\delta_n$ is as in Proposition~\ref{prop:L^2modification}. 
\end{lemma}
\begin{proof}
	We note that 
	\begin{align*}
		&\quad\ \mathbb{P}\big(\widehat D_{j+k,\alpha}\neq\cR^k((\widehat D_{j,\alpha w})_{w\in I^k})\big)\\
		&\leq \mathbb{P}\big(\widehat D_{j+k,\alpha}\neq D_{j+k,\alpha}\big)+\mathbb{P}\big(\cR^k((\widehat D_{j,\alpha w})_{w\in I^k})\neq\cR^k((D_{j,\alpha w})_{w\in I^k})\big)\\
		&\leq \delta_{j+k}+N^k\delta_j,
	\end{align*} 
	and $\mathbb E\bigl\|\widehat D_{j+k,\alpha}
	-\cR^k((\widehat D_{j,\alpha w})_{w\in I^k})\bigr\|^2\leq C$ by Proposition \ref{prop:L^2modification} and Lemma \ref{lem:moment2}. So the estimate follows from the Cauchy--Schwarz inequality.
\end{proof}

\subsection{Homogenization}
We prove Theorem \ref{thm:random-homogenization} in this subsection. Throughout this subsection, we assume Assumption~\ref{ch3:3-2}. We use the same notation as in Proposition \ref{prop:L^2modification}. \smallskip

First, we show that $\widehat D_{n,\alpha}-\mathbb E\widehat D_{n,\alpha}$ converges to zero in $L^1$.

\begin{lemma}\label{lem:5.3}
Under Assumption~\ref{ch3:3-2}, there is $a>0$ so that $\lim\limits_{j\to\infty}\mathbb{P}(\widehat D_{j,\alpha}\not\geq aE_0)=0$.
\end{lemma}
\begin{proof}
This is an immediate consequence of Lemma \ref{lem:4.10} and Proposition \ref{prop:L^2modification}. 
\end{proof}

\begin{lemma}\label{lem:5.4}
Under Assumption~\ref{ch3:3-2},
\begin{equation}\label{e:5.6}
 \mathbb E\big\|\widehat D_{n,\alpha}-\mathbb E\widehat D_{n,\alpha}\big\|\longrightarrow0.
\end{equation}
\end{lemma}
\begin{proof}
Choose $a>0$ as in Lemma~\ref{lem:5.3}, and let $P_j:=\mathbb{P}(\widehat D_{j,\alpha}\not\geq aE_0)$. Then, by Proposition~\ref{prop:L^2modification}, Lemma \ref{quasi-trace} and Lemma \ref{lem:4.9},
\begin{align*}
&\quad\ \mathbb E\big\|\widehat D_{j+k,\alpha}-\mathbb E\widehat D_{j+k,\alpha}\big\|\\
&\leq 
 2\mathbb E\big\|\widehat D_{j+k,\alpha}-\cR^k((\widehat D_{j,\alpha w})_{w\in I^k})\big\|+\mathbb E\big\|\cR^k((\widehat D_{j,\alpha w})_{w\in I^k})-\mathbb{E}\cR^k((\widehat D_{j,\alpha w})_{w\in I^k})\big\|\\
&\leq C_1(\delta_{j+k}+N^k\delta_j)^{1/2}+C_2(N^kP_j)^{1/2}+C_3\gamma^k,
\end{align*}
where $C_1,C_2,C_3$ are constants independent of $j,k$.
For fixed $k$, let $j\to\infty$. Lemma~\ref{lem:5.3} and
Proposition~\ref{prop:L^2modification} give $P_j\to0$ and $\delta_j\to0$, so
\[
 \limsup_{n\to\infty}
 \mathbb E\big\|\widehat D_{n,\alpha}-\mathbb E\widehat D_{n,\alpha}\big\|
 \le C_3\gamma^k.
\]
Letting $k\to\infty$ proves \eqref{e:5.6}.
\end{proof}

Now we proceed to show that $\widehat D_{n,\alpha}$ converges in $L^1$ to a fixed
form $E_+\in\cQ_p(V_0)$. In the following, we write $\hat{E}_n:=\mathbb{E}\widehat D_{n,\alpha}$ for $n\geq 0$. We first collect some properties of $\widehat E_n$. 

\begin{lemma}[Uniform boundedness]\label{lem:5.5}
There are $0<a_0<b_0<\infty$ so that $\widehat E_n\in [a_0E_0,b_0E_0]$ for $n\geq 0$.
\end{lemma}
\begin{proof}
First, $\|\hat{E}_n\|=\|\mathbb{E}\widehat D_{n,\alpha}\|\leq\mathbb{E}\|\widehat D_{n,\alpha}\|\leq\sqrt{\mathbb{E}\|\widehat D_{n,\alpha}\|^2}$, so we can take $b_0=\sup_{n\geq 0}\|\hat{E}_n\|<\infty$ by Proposition~\ref{prop:L^2modification}. Next, Lemma~\ref{lem:5.3} gives $\widehat E_n\ge(a/2)E_0$ for all sufficiently large $n$, which proves the lower bound.
\end{proof}

\begin{lemma}[Quasi-Jensen inequality]\label{quasi-Jensen}
There are $\varepsilon_n\ge0$ with $\sum_{n=0}^\infty\varepsilon_n<\infty$
such that
	$\hat{E}_{n+1}\le\cR(\hat{E}_n)+\varepsilon_nE_0$.
\end{lemma}
\begin{proof}
	By applying Lemma \ref{quasi-trace} and then applying the concavity of $\cR$ and Jensen's inequality, we obtain
	\begin{align*}
	\mathbb{E}\hat{D}_{n+1,\alpha}&\leq C(\delta_{n+1}+N\delta_n)^{1/2}E_0+\mathbb{E}\cR((\widehat D_{n,\alpha i})_{i\in I})\\
		&\leq C(\delta_{n+1}+N\delta_n)^{1/2}E_0+\mathcal{R}(\mathbb{E}\hat{D}_{n,\alpha}).  
	\end{align*}
	So it suffices to take $\varepsilon_n=C(\delta_{n+1}+N\delta_n)^{1/2}$.
\end{proof}

\begin{lemma}[Quasi-trace]\label{lem:5.7}
$\big\|\widehat E_{n+1}-\cR(\widehat E_n)\big\|\longrightarrow0 \quad\text{as }n\longrightarrow\infty$.
\end{lemma}
\begin{proof}
By Lemmas \ref{quasi-trace}, \ref{lem:trace-lipschitz} and \ref{lem:5.4}, 
\[
 \begin{aligned}
 \big\|\widehat E_{n+1}-\cR(\widehat E_n)\big\|&=\big\|\mathbb{E}\widehat D_{n+1,\alpha}-\cR(\mathbb E \widehat D_{n,\alpha})\big\|\\
 &\le\bigl\|\mathbb{E}\widehat D_{n+1,\alpha}-\mathbb E\cR((\widehat D_{n,\alpha i})_{i\in I})\bigr\|+\bigl\|\mathbb E\cR((\widehat D_{n,\alpha i})_{i\in I})-\cR(\mathbb E \widehat D_{n,\alpha})\bigr\|\\
 &\le C(\delta_{n+1}+N\delta_n)^{1/2}+C_{\cR}N
       \mathbb E\|\widehat D_{n,\alpha}-\mathbb E \widehat D_{n,\alpha}\|
 \longrightarrow0.
 \end{aligned}
\]
\end{proof}

\begin{theorem}\label{thm:5.8}
There is a fixed form $E_+\in\cQ_p(V_0)$ such that $\cR(E_+)=E_+$ and
\begin{equation}\label{eq:ch3-3-39}
 \mathbb E\|\widehat D_{n,\alpha}-E_+\|\longrightarrow0\quad\hbox{ as }n\to\infty.
\end{equation}
Consequently, $D_{n,\alpha}\to E_+$ in probability.
\end{theorem}
\begin{proof}
We first prove the $L^1$ convergence. By Lemma~\ref{lem:5.4}, we only need to show
that there exists $E_+\in\cQ_p(V_0)$ with $\cR(E_+)=E_+$ such that
$\|\widehat E_n-E_+\|\to0$. Let $\varepsilon_n$ be as in Lemma~\ref{quasi-Jensen}, and set
\[
 \overline\varepsilon_n=\sum_{\ell=n}^\infty\varepsilon_\ell<\infty,
 \qquad\widetilde E_n=\widehat E_n+\overline\varepsilon_nE_0.
\]
In order to prove that $\widehat E_n$ converges to some $E_+$, we only need
to prove that $\widetilde E_n$ converges to $E_+$, by the definition of
$\widetilde E_n$. 

The good properties of $\hat{E}_n$ listed in Lemmas \ref{lem:5.5}--\ref{lem:5.7} are inherited by $\widetilde E_n$: 
\begin{align}
&\widetilde E_n\in[a_0E_0,(b_0+\overline\varepsilon_0)E_0]\hbox{ for }n\geq0,\tag{Uniform boundedness}\label{UB}\\
&\widetilde E_{n+1}\leq \cR(\widetilde E_n)\hbox{ for }n\geq 0,\tag{Jensen}\label{Jen}\\
&\|\cR(\widetilde E_n)-\widetilde E_{n+1}\|\longrightarrow0\hbox{ as }n\to\infty.\tag{Quasi-trace}\label{QT} 
\end{align}
By Theorem~\ref{thm:peirone}, we can define
$\cR^\infty=\lim_{k\to\infty}\cR^k$ pointwise on $\cQ_p(V_0)$. Then, by \eqref{Jen}, 
\[
 \cR^\infty(\widetilde E_{n+1})
 \le\cR^\infty(\cR(\widetilde E_n))
 =\cR^\infty(\widetilde E_n).
\]
Moreover, we also have $\cR^\infty(\widetilde E_n)\in[a_0E_0,(b_0+\overline\varepsilon_0)E_0]$ by \eqref{UB} and Lemma \ref{lem:compact-annulus}(a). By Lemma~\ref{lem:compact-annulus}(c), this decreasing sequence converges in norm, so we can define
\[
E_+:=\lim\limits_{n\to\infty}\cR^\infty(\widetilde E_n).
\]

The continuity of $\cR$ implies $\cR(E_+)=E_+$. It remains to prove $E_+=\lim\limits_{n\to\infty}\widetilde E_n$. Let 
\[
 \mathcal O:=\{G:\widetilde E_{n_j}\to G
                 \text{ for some }n_j\to\infty\}.
\]
By \eqref{UB} and Lemma~\ref{lem:compact-annulus}(c),
$\mathcal O$ is nonempty and compact. It suffices to show that $\mathcal O=\{E_+\}$. To show this, we will use the following two properties of $\mathcal O$: 
\begin{enumerate}
\item[(i)] $\mathcal{O}\subset\cR(\mathcal{O})$. 
\item[(ii)] For each $G\in \mathcal O$, $\cR^\infty(G)=E_+$.  
\end{enumerate}
For any $G\in\mathcal{O}$, there is a subsequence $\widetilde{E}_{n_j}$ such that $\widetilde{E}_{n_j}\to G$ as $j\to\infty$. By passing to a further subsequence, we assume $\widetilde{E}_{n_j-1}\to F$ for some $F\in\mathcal{O}$, so  by the continuity of $\cR$ and \eqref{QT},
\[
G=\lim_{j\to\infty}\widetilde E_{n_j}=\lim_{j\to\infty}\cR(\widetilde E_{n_j-1})=\cR(F)\in \cR(\mathcal{O}).
\]
This proves (i).  By Lemma \ref{lem:compact-annulus}(b), $\cR^\infty=\lim_{k\to\infty}\cR^k$ is continuous on $[a_0E_0,(b_0+\overline{\varepsilon}_0)E_0]$, so
\[
\cR^\infty(G)=\lim_{j\to\infty}\cR^\infty(\widetilde E_{n_j})=E_+.
\] 
This proves (ii).

Now, for each $n$, choose $G_n\in\mathcal O$ such that $\cR^n(G_n)=G$.
Along a subsequence $n_j$, we have $G_{n_j}\to G'\in\mathcal O$. Hence
\[
 \begin{aligned}
 \|G-E_+\|
 &=\|G-\cR^\infty(G')\|\\
 &\le\|G-\cR^{n_j}(G')\|
      +\|\cR^{n_j}(G')-\cR^\infty(G')\|\\
 &=\|\cR^{n_j}(G_{n_j})-\cR^{n_j}(G')\|
      +\|\cR^{n_j}(G')-\cR^\infty(G')\|\\
 &\le\frac{b_0+\overline\varepsilon_0}{a_0}\|G_{n_j}-G'\|
      +\|\cR^{n_j}(G')-\cR^\infty(G')\|
 \longrightarrow0,
 \end{aligned} 
\]
where the last inequality uses Lemma \ref{lem:compact-annulus}(b). Thus $\mathcal O=\{E_+\}$.

Finally, the convergence in probability of $D_{n,\alpha}\to E_+$ follows since $\mathbb P(D_{n,\alpha}\ne\widehat D_{n,\alpha})\le\delta_n\to0$ by Proposition~\ref{prop:L^2modification}, which finishes the proof.
\end{proof}\smallskip

\begin{proposition}\label{prop5.10}
Let $E_+$ be the $p$-resistance form of Theorem \ref{thm:5.8}; then for each $u\in\cS$,
\[
\lim_{n\to\infty}D_{n,\alpha}(u)=E_+(u)\quad\hbox{ a.s.}
\]   
\end{proposition}
\begin{proof}
For simplicity, assume $E_0=E_+$. We break the proof into two parts: 
\begin{align}
	&\limsup_{n\to\infty}D_{n,\alpha}(u)\le E_0(u)\quad\text{a.s.},\label{step1:1}\\
	&\liminf_{n\to\infty}D_{n,\alpha}(u)\ge E_0(u)\quad\ \text{a.s.}\label{step1:2}
\end{align} 

First, we show \eqref{step1:1}. For $n\geq 0$, write $n=k_n+j_n$, where $j_n=\lfloor n/2\rfloor$. We claim that
\begin{equation}\label{e:eas}
	D_{n,\alpha}=
	\cR^{k_n}((\widehat D_{j_n,\alpha w})_{w\in I^{k_n}})
	\quad\text{eventually a.s.}
\end{equation}
Indeed, this follows from
\[
\begin{aligned}
	\sum_{n\ge0}\mathbb P\left( D_{n,\alpha}\neq 
	\cR^{k_n}((\widehat D_{j_n,\alpha w})_{w\in I^{k_n}})\right)
	&\leq\sum_{n\ge0}\mathbb P\left(\bigcup_{w\in I^{k_n}}\mathcal H_{j_n,\alpha w}\right)\\
	& \le\sum_{n\ge0}N^{k_n}\delta_{j_n}\le C\sum_{j\ge0}N^j\exp(-cr^{-j})<\infty
\end{aligned}
\]
and the Borel--Cantelli lemma. Let $h$ be the harmonic extension of $u$ with respect to the form $f\mapsto\sum_{w\in I^n}E_0(f\circ\psi_w)$. Then, \eqref{e:eas} implies
\begin{align*}
	D_{n,\alpha}(u)&\le \rho^{k_n}\sum_{w\in I^{k_n}}\widehat D_{j_n,\alpha w}(h\circ\psi_w)\\
	&\le\left(\sum_{w\in I^{k_n}}
	a_w(u)\|\widehat D_{j_n,\alpha w}\|\right)E_0(u)
	\quad\text{eventually a.s.}
\end{align*}
We also observe that $\mathbb E\|\widehat D_{j_n,\alpha}\|\to1$ as $n\to\infty$ by Theorem~\ref{thm:5.8}. So, to show \eqref{step1:1}, it suffices to prove
\[
\sum_{w\in I^{k_n}}a_w(u)\|\widehat D_{j_n,\alpha w}\|
-\mathbb E\|\widehat D_{j_n,\alpha}\|\longrightarrow 0
\quad\text{a.s.}
\]
This follows from a routine probability argument and the following estimate:
\[
\begin{aligned}
	\sum_{n\ge0}\operatorname{Var}\left(
	\sum_{w\in I^{k_n}}a_w(u)\|\widehat D_{j_n,\alpha w}\|\right)
	\le C\sum_{n\ge0}\sum_{w\in I^{k_n}}a_w(u)^2
	\le C\sum_{n\ge0}\theta^{k_n}<\infty,
\end{aligned}
\]
where we have used $\sup_n\mathbb E\|\widehat D_{n,\alpha}\|^2<\infty$ and Corollary \ref{coro:4.2}. \medskip

Next, we prove \eqref{step1:2}. Let $j,k\ge0$. By Lemma~\ref{lem:3.3},
\[
D_{j+k,\alpha}(u)
\ge E_0(u)\left(\sum_{w\in I^k}
a_w(u)\llbracket D_{j,\alpha w}\rrbracket^{1/(p-1)}\right)^{-(p-1)}.
\]
For each $j\geq 0$, Proposition~\ref{lem:tree-averages} shows 
\[
\sum_{w\in I^k}a_w(u)\llbracket D_{j,\alpha w}\rrbracket^{1/(p-1)}
\longrightarrow\mathbb E\llbracket D_{j,\alpha}\rrbracket^{1/(p-1)}\hbox{ as }k\to\infty\text{ a.s.}
\]
Consequently,
\begin{equation}\label{step1:3}
	\liminf_nD_{n,\alpha}(u)
	\ge E_0(u)\bigl(\mathbb E\llbracket D_{j,\alpha}\rrbracket^{1/(p-1)}\bigr)^{-(p-1)}
	\quad\text{a.s.}
\end{equation}
By Theorem~\ref{thm:5.8}, $D_{j,\alpha}\to E_0$ in probability, so $\llbracket D_{j,\alpha}\rrbracket^{1/(p-1)}\to1$ in probability by the continuity of $E\mapsto\llbracket E\rrbracket$ at $E_0$. Since $\{\llbracket D_{j,\alpha}\rrbracket^{1/(p-1)}:j\ge0\}$ is uniformly integrable by Lemma \ref{lem:lower-bound-moment},
\begin{equation}
	\lim_{j\to\infty}\mathbb E\llbracket D_{j,\alpha}\rrbracket^{1/(p-1)}=1.
\end{equation}
Therefore, letting $j\to\infty$ in \eqref{step1:3} yields \eqref{step1:2}.
\end{proof}

Finally, we complete the proof of Theorem \ref{ch3:3-3} by upgrading the pointwise convergence to uniform convergence.
\begin{proof}[Proof of Theorem~\ref{ch3:3-3}]
First, we prove $\|D_{n,\alpha}-E_+\|\to0$ a.s. under Assumption \ref{ass:random}. 

This is equivalent to proving the uniform convergence of $D_{n,\alpha}$ to $E_+$ on $\cS$ a.s. Let $\cS'$ be a countable dense subset of $\cS$ containing a basis $u_1,\ldots,u_{M-1}$ of $\R^{V_0}/\R\mathbf1$. Proposition \ref{prop5.10} then implies that
\[
D_{n,\alpha}\to E_+\hbox{ pointwise on $\mathcal{S}'$ a.s.} 
\]
At a sample point $\omega$ at which this convergence holds, by Lemma \ref{lem:Lipschitz} and \eqref{step2:2}, 
\begin{equation}
\sup_{n\geq 0}\operatorname{Lip}_{\cS}(D_{n,\alpha}(\omega))\leq p\sup_{n\geq 0}\|D_{n,\alpha}(\omega)\|\le p\cdot C\max_{1\le j\le M-1}\sup_{n\geq 0}D_{n,\alpha}(\omega)(u_j)<\infty,
\end{equation}
so by the Arzel\`a--Ascoli theorem, $D_{n,\alpha}(\omega)$ converges to $E_+$ uniformly on $\cS$. \medskip 

{Next, for any $1\leq t<\infty$, we prove $\mathbb E\|D_{n,\alpha}-E_+\|^t\to0$ under the additional assumption that $\mathbb{E}\|D_{0,\alpha}\|^q<\infty$ for some $q>0$.}

Let $r>t$. By Lemmas~\ref{lem:moment-upgrade} and \ref{lem:moment2}, there is $n_0\ge0$ so that
\[
 \sup_{n\ge n_0}\mathbb E\|D_{n,\alpha}-E_+\|^r<\infty.
\]
Consequently,
\[
 \sup_{n\ge n_0}\mathbb E\left[
 \|D_{n,\alpha}-E_+\|^t;
 \|D_{n,\alpha}-E_+\|^t>A\right]
 \le A^{1-r/t}\sup_{n\ge n_0}\mathbb E\|D_{n,\alpha}-E_+\|^r
 \longrightarrow0\quad\hbox{ as }A\to\infty. 
\]
This shows that $\big\{\|D_{n,\alpha}-E_+\|^t:n\geq n_0\big\}$ is uniformly integrable. Since we have proved $\|D_{n,\alpha}-E_+\|\to 0$ a.s., it follows that $\mathbb E\|D_{n,\alpha}-E_+\|^t\to0$.
\end{proof}

\section{Mosco convergence of random environments}
\label{sec:compact-limit}
\noeqref{eq:isometry-core,eq:averaged-resistance,eq:cell-distance,eq:cell-chain-average,eq:weighted-cell-sums,eq:weighted-cell-reindexing,eq:average-inverse-bound,eq:M1,eq:M2,eq:compactness-approximation,eq:compactness-trace-limit}

In this section, we prove the Mosco convergence of the rescaled random
environments to the self-similar form on a cell
$\Psi_\alpha(K)$, $\alpha\in\cW$. For simplicity, we identify this fixed $\Psi_\alpha(K)$ with $K$. We also prove Mosco convergence for ${\cD^{(n)}}$ as a corollary. Recall \eqref{eq:infinite-skeleton-energy} for the definition of ${\mathcal{D}^{(n)}}$.

Throughout this section, we assume Assumption~\ref{ass:random}, and assume $E_0=E_+$ in Theorem \ref{thm:random-homogenization}.

\subsection{Self-similar forms}\label{sec:self-similar-forms}
For $n\ge0$, set
\[
 \begin{aligned}
 \cE_{0,n}(u)
 &=\rho^n\sum_{w\in I^n}E_0\bigl(u\circ\psi_w|_{V_0}\bigr),\\
 {D_\alpha^{(n)}}(u)
 &=\rho^n\sum_{w\in I^n}
 D_{0,\alpha w}\bigl(u\circ\psi_w|_{V_0}\bigr)
 \quad\hbox{ for }u\in\R^{V_n}.
 \end{aligned}
\]

As $\cR(E_0)=E_0$, it is not hard to see that
\[
 \Tr_{V_k}\cE_{0,n}=\cE_{0,k}\quad\hbox{ for }0\le k\le n.
\]
It follows that $\cE_{0,k}(u|_{V_k})$ is increasing in $k$
for every function $u$ on $K$. Thus we can define
\[
 \cF_0:=\left\{u\in C(K):
 \lim_{k\to\infty}\cE_{0,k}(u|_{V_k})<\infty\right\},
\]
and
\[
 \cE_0(u):=\lim_{k\to\infty}\cE_{0,k}(u|_{V_k})\quad\hbox{ for }u\in\cF_0.
\]
We call $\cE_0$ a \textit{self-similar $p$-energy} on $K$.
See \cite[Section~8]{KajinoShimizu2026} for the fact that
$(\cE_0,\cF_0)$ is a $p$-resistance form on $K$. 

For $n\ge k\geq 0$ and $u\in\R^{V_k}$, by the definition {of} $D_{m,\beta}$ in Section \ref{sec:random-environment}, 
\[
 (\Tr_{V_k}{D_\alpha^{(n)}})(u)
 =\rho^k\sum_{\tau\in I^k}
 D_{n-k,\alpha\tau}\bigl(u\circ\psi_\tau|_{V_0}\bigr).
\]

\subsection{Self-similar measures}\label{sec:self-similar-measures}
Let $\mu$ be the equally weighted self-similar probability measure on $K$ (see \cite{Hutchinson1981} for the construction),
which is characterized by the equality
\begin{equation}\label{e:mu}
 \mu=\frac1N\sum_{i\in I}\mu\circ\psi_i^{-1}.
\end{equation}
$\mu$ is non-atomic and has support $K$. Hence, by \eqref{e:mu}, 
\begin{equation}\label{e:mu2}
 \mu(\psi_w(K))=\mu(\psi_w(K\setminus V_0))=N^{-|w|}\ \hbox{ for every }w\in I^*.
\end{equation}
For each $n\geq 0$, we also define a probability measure $\mu_n$ on $V_n=\bigcup_{w\in I^n}\psi_w(V_0)$ by 
\[
 \mu_n=\frac1{MN^n}\sum_{w\in I^n}\sum_{x\in V_0}
 \delta_{\psi_w(x)},
\]
where $\delta_y$ denotes the delta mass {at} a point $y$. 

Next, for each $n\geq 0$, we introduce an embedding map $\iota_n:L^t(V_n,\mu_n){\to} L^t(K,\mu)$ for $1\le t\le\infty$. For this purpose, we first define $\pi_n:K\to V_n$. 
By the nesting property of $K$, we can find a {Borel partition} $(\mathscr B_w)_{w\in I^n}$ of $K$ such that 
\[
\psi_w(K\setminus V_0)\subset\mathscr B_w\subset\psi_w(K)\ \hbox{ and }\  \bigcup_{w\in I^n}{\mathscr B_w}=K.
\]
Clearly, $\mu(\mathscr{B}_w)=N^{-n}$ for $w\in I^n$. Since $\mu$ is nonatomic, by \cite[Section~1.12(iii)]{Bogachev2007}, for each $w\in I^n$, we can find {a Borel partition $(\mathscr B_{w,x})_{x\in V_0}$ of $\mathscr B_w$} so that  
\[
\bigcup_{x\in V_0}\mathscr B_{w,x}=\mathscr B_w\quad \hbox{ and }\quad \mu(\mathscr B_{w,x})=M^{-1}N^{-n},\, \{\psi_w(x)\}\cap \mathscr B_w\subset \mathscr B_{w,x}\hbox{ for each }x\in V_0. 
\] 
Then, define $\pi_n:K\to V_n$ by
\[
 \pi_n(z)=\psi_w(x)\qquad\text{if }z\in\mathscr B_{w,x}.
\]
It is not {hard} to see that $\pi_n(x)=x$ for each $x\in V_n$, and that $\mu_n$ is the pushforward measure of $\mu$ by $\pi_n$, i.e.
\begin{equation}\label{e:pushforward}
\mu_n=\mu\circ\pi^{-1}_n. 
\end{equation}
For each $u\in\R^{V_n}$, we define $\iota_n u$ as the Borel measurable function on $K$ with the formula
\[
\iota_nu:=u\circ\pi_n. 
\]

\begin{proposition}\label{prop:isometric-embedding}
\begin{enumerate}[label=\textup{(\alph*)}]
\item $\|u\|_{L^t(V_n,\mu_n)}=\|{\iota_nu}\|_{L^t(K,\mu)}$ for any $n\geq 0$, $u\in\R^{V_n}$ and $t\in[1,\infty]$.
\item For $f\in C(K)$,
\begin{equation}\label{eq:isometry-core}
 \|\iota_n(f|_{V_n})-f\|_{L^\infty(K,\mu)}\to0\quad \hbox{ as }n\to\infty. 
\end{equation}
\end{enumerate}
\end{proposition}
\begin{proof}
(a) is an immediate consequence of \eqref{e:pushforward}.

(b) For each $z\in K$, we notice that $z$ and $\pi_n(z)$ belong {to} the same $n$-cell, so 
\[
 |f(\pi_n(z))-f(z)|
 \le\max_{|w|=n}\operatorname{osc}_{\psi_w(K)}f.
\]
It follows that $\|\iota_n(f|_{V_n})-f\|_{L^\infty(K,\mu)}\leq \max_{|w|=n}\operatorname{osc}_{\psi_w(K)}f\to 0$ as $n\to\infty$.
\end{proof}

\subsection{Resistance metrics and related estimates}
\label{sec:resistance-estimates}For a $p$-resistance form
$(\mathcal D,\operatorname{dom}(\mathcal D))$ on $V_n$ or $K$, set
\[
 R_{\mathcal D}(x,y):=\sup_{\substack{u\in\operatorname{dom}(\mathcal D)\\ \mathcal D(u)>0}}
 \frac{|u(x)-u(y)|^p}{\mathcal D(u)}\ \hbox{ and }\ 
 d_{\mathcal D}(x,y):=R_{\mathcal D}(x,y)^{1/(p-1)}. 
\]
By \cite[Corollary~6.32 and Definition~6.33]{KajinoShimizu2026}, $d_{\mathcal D}$ is a metric. 
For short, in the {remainder} of this section, we write
\[
 R_0=R_{\cE_0},\qquad d_0=R_0^{1/(p-1)},\qquad
 R_n^\alpha=R_{{D_\alpha^{(n)}}},\qquad
 d_n^\alpha=(R_n^\alpha)^{1/(p-1)}.
\]
By Propositions~8.40--8.41, Lemma~8.42 and Theorem~8.43 of
\cite{KajinoShimizu2026}, $d_0$ induces the original
Euclidean topology on $K$.

We will use the following basic fact, which follows directly from the
definition.

\begin{proposition}\label{prop:trace-resistance}
Suppose that $(\cE,\cF)$ is a $p$-resistance form on $V$ and that
$B\subset V$ contains at least two points. Then,  $R_{\cE}|_{B\times B}=R_{\Tr_B\cE}$, i.e.
\[
 R_{\cE}(x,y)=R_{\Tr_B\cE}(x,y)\quad\hbox{ for every }x,y\in B.
\]
\end{proposition}

\begin{lemma}\label{lem:cell-distance}
Let $C=\operatorname{diam}_{d_0}V_0<\infty$. Then, (a) and (b) hold.
\begin{enumerate}[label=\textup{(\alph*)}]
\item For $1\le j\le n$ and $w\in I^{j-1}$,
\begin{equation}\label{eq:cell-distance}
 \operatorname{diam}_{d_n^\alpha}\psi_w(V_1)
 \le C\rho^{-j/(p-1)}\sum_{i\in I}
 \llbracket D_{n-j,\alpha wi}\rrbracket^{1/(p-1)}.
\end{equation}
\item For $0\le k<n$,
\begin{equation}\label{eq:cell-chain-average}
 \int_Kd_n^\alpha(\pi_n(z),\pi_k(z))\,d\mu(z)
 \leq CN\sum_{j=k+1}^{n}N^{-j}\rho^{-j/(p-1)}\sum_{w\in I^j}\llbracket D_{n-j,\alpha w}\rrbracket^{1/(p-1)}.
\end{equation}
\end{enumerate}
\end{lemma}
\begin{proof}
(a) Notice that, for $i\in I$ and $f\in\R^{V_n}$,
\[
 {D_\alpha^{(n)}}(f)
 \ge\rho^jD_{n-j,\alpha wi}(f\circ\psi_{wi}|_{V_0})
 \ge\frac{\rho^j}{\llbracket D_{n-j,\alpha wi}\rrbracket}
 E_0(f\circ\psi_{wi}|_{V_0}).
\]
Thus, for $x,y\in V_0$,
\[
 d_n^\alpha(\psi_{wi}x,\psi_{wi}y)
 \le\rho^{-j/(p-1)}\llbracket D_{n-j,\alpha wi}\rrbracket^{1/(p-1)}
 d_0(x,y).
\]
Then
\[
 \begin{aligned}
 \operatorname{diam}_{d_n^\alpha}\psi_w(V_1)
 &\le\sum_{i\in I}\operatorname{diam}_{d_n^\alpha}\psi_{wi}(V_0)\\
 &\le
 \sum_{i\in I}\rho^{-j/(p-1)}\llbracket D_{n-j,\alpha wi}\rrbracket^{1/(p-1)}\bigl(\operatorname{diam}_{d_0}V_0\bigr),
 \end{aligned}
\]
which is \eqref{eq:cell-distance}.

(b) For every $1\le j\le n$, ${w\in I^{j-1}}$ and $z\in\psi_w(K)$, noticing that  $\pi_{j-1}(z),\pi_j(z)\in\psi_w(V_1)$, we have by (a) that 
\[
 d_n^\alpha(\pi_{j-1}(z),\pi_j(z))
 \le C\rho^{-j/(p-1)}\sum_{i\in I}
 \llbracket D_{n-j,\alpha wi}\rrbracket^{1/(p-1)}.
\]
By integrating the above estimate over $K$, we see 
\[
 \begin{aligned}
 &\quad\ \int_Kd_n^\alpha(\pi_{j-1}(z),\pi_j(z))\,d\mu(z)\\
 &=\sum_{w\in I^{j-1}}\int_{\psi_w(K)}
 d_n^\alpha(\pi_{j-1}(z),\pi_j(z))\,d\mu(z)\\
 &\le C\rho^{-j/(p-1)}N^{-(j-1)}
 \sum_{|w|=j-1}\sum_{i\in I}
 \llbracket D_{n-j,\alpha wi}\rrbracket^{1/(p-1)}\\
 &=CN^{1-j}\rho^{-j/(p-1)}
 \sum_{|w|=j}\llbracket D_{n-j,\alpha w}\rrbracket^{1/(p-1)}.
 \end{aligned}
\]
\eqref{eq:cell-chain-average} follows by summing the above inequality over $j=k+1,\ldots,n$.
\end{proof}

The next two lemmas estimate the right-hand side of \eqref{eq:cell-chain-average}. 

\begin{lemma}\label{lem:weighted-cell-sums}
For $k\ge0$ and $1\le m\le M-1$, there are nonnegative weights
$a^{(m,k)}_w$, $w\in I^*$, satisfying the conditions of
Proposition~\ref{lem:tree-averages}, such that for any $n>k$,
\begin{equation}\label{eq:weighted-cell-sums}
 \begin{aligned}
 &\sum_{j=k+1}^{n}N^{-j}\rho^{-j/(p-1)}
 \sum_{w\in I^j}\llbracket D_{n-j,\alpha w}\rrbracket^{1/(p-1)}\\
 &\qquad\qquad\le C\max_{1\leq m\leq M-1}\rho^{-k/(p-1)}\sum_{w\in I^n}
a^{(m,k)}_w\llbracket D_{0,\alpha w}\rrbracket^{1/(p-1)},
 \end{aligned}
\end{equation}
where the constant $C>0$ depends only on {$E_0$, $K$ and $p$}. 
\end{lemma}
\begin{proof}
Take $u_1,\ldots,u_{M-1}\in\mathcal S$ as in Lemma~\ref{lem:3.5}. 
Then, for $|w|=j$ and $j\le n$, 
\[
 \llbracket D_{n-j,\alpha w}\rrbracket^{1/(p-1)}
 \le C_1\max_{1\leq m\leq M-1}\sum_{v\in I^{n-j}}
 a_v(u_m)\llbracket D_{0,\alpha wv}\rrbracket^{1/(p-1)},
\]
where $a_v(u_m)$ is as in Definition~\ref{def:energy-distribution}. 
Thus
\begin{equation}\label{eq:weighted-cell-reindexing}
 \begin{aligned}
 &\quad\ \sum_{j=k+1}^{n}N^{-j}\rho^{-j/(p-1)}
 \sum_{w\in I^j}\llbracket D_{n-j,\alpha w}\rrbracket^{1/(p-1)}\\
 &\le {C_1\sum_{j=k+1}^{n}
 N^{-j}\rho^{-j/(p-1)}\sum_{w\in I^j}\max_{1\leq m\leq M-1}\sum_{v\in I^{n-j}}
 a_v(u_m)\llbracket D_{0,\alpha wv}\rrbracket^{1/(p-1)}}\\
 &\le {C_1(M-1)}\max_{1\leq m\leq M-1}\sum_{j=k+1}^{n}
 N^{-j}\rho^{-j/(p-1)}\sum_{w\in I^j}\sum_{v\in I^{n-j}}
 a_v(u_m)\llbracket D_{0,\alpha wv}\rrbracket^{1/(p-1)}\\
 &={C_1(M-1)}\max_{1\leq m\leq M-1}\sum_{\eta\in I^n}
 \left(\sum_{j=k+1}^{n}N^{-j}\rho^{-j/(p-1)} 
 a_{\eta_{j+1}\cdots\eta_n}(u_m)\right)
 \llbracket D_{0,\alpha\eta}\rrbracket^{1/(p-1)}\\
 &\leq {C_1(M-1)}\max_{1\leq m\leq M-1}\sum_{\eta\in I^n}
 \left(\sum_{j=k+1}^{\infty}\rho^{-j/(p-1)}\widehat a^{(m,j)}_\eta\right)
 \llbracket D_{0,\alpha\eta}\rrbracket^{1/(p-1)}\\
 &={C_1(M-1)}\max_{1\leq m\leq M-1}\rho^{-k/(p-1)}\sum_{\eta\in I^n}
 \left(\sum_{j=k+1}^{\infty}\rho^{-(j-k)/(p-1)}\widehat a^{(m,j)}_\eta\right)
 \llbracket D_{0,\alpha\eta}\rrbracket^{1/(p-1)},
 \end{aligned}
\end{equation}
where, for $1\leq m\leq M-1$, $j\geq 0$ and $\eta\in I^*$,
\[
\widehat a^{(m,j)}_\eta:=\begin{cases}
 N^{-|\eta|},&\hbox{ if }|\eta|\le j,\\
 N^{-j}a_{\eta_{j+1}\cdots \eta_{|\eta|}}(u_m),&\hbox{ if }|\eta|>j.
 \end{cases}
\]
We take 
\[ 
a^{(m,k)}_\eta=(\rho^{1/(p-1)}-1)\sum_{j=k+1}^{\infty}\rho^{-(j-k)/(p-1)}\widehat a^{(m,j)}_\eta.
\]
Then \eqref{eq:weighted-cell-sums} follows from \eqref{eq:weighted-cell-reindexing}. 

Next, we observe that by Corollary \ref{coro:4.2}, for $1\leq m\leq M-1$ and $j\geq 0$,
\[
 \widehat a^{(m,j)}_\emptyset=1,\qquad
 \sum_{i\in I}\widehat a^{(m,j)}_{\eta i}=\widehat a^{(m,j)}_\eta,\qquad
\widehat a^{(m,j)}_{\eta i}\le\max\{N^{-1},\theta\}\widehat a^{(m,j)}_\eta\hbox{ for }i\in I,
\]
where $\theta$ is the constant in Corollary~\ref{coro:4.2}. Since $\{{a^{(m,k)}_\eta}\}_{\eta\in I^*}$ is a weighted average of
$\{\widehat a^{(m,j)}_\eta\}_{\eta\in I^*},j\geq k+1$, 
\[
a^{(m,k)}_\emptyset=1,\qquad
 \sum_{i\in I}a^{(m,k)}_{\eta i}=a^{(m,k)}_\eta,\qquad
a^{(m,k)}_{\eta i}\le\max\{N^{-1},\theta\}a^{(m,k)}_\eta\hbox{ for }i\in I.
\]
\end{proof}

\begin{lemma}\label{lem:average-inverse-bound}
For every $k\ge0$,
\begin{equation}\label{eq:average-inverse-bound}
 \limsup_{n\to\infty}
 \sum_{j=k+1}^{n}\rho^{-j/(p-1)}N^{-j}
 \sum_{|w|=j}\llbracket D_{n-j,\alpha w}\rrbracket^{1/(p-1)}
 \le C\rho^{-k/(p-1)}\quad\hbox{ a.s.}
\end{equation}
\end{lemma}
\begin{proof}
By Assumption~\ref{ass:random} and Proposition~\ref{lem:tree-averages},
for each fixed $k\geq 0$ and $1\leq m\leq M-1$,
\[
\lim_{n\to\infty}\sum_{w\in I^n}{a^{(m,k)}_w}
 \llbracket D_{0,\alpha w}\rrbracket^{1/(p-1)}
=\mathbb E\llbracket D_{0,\alpha}\rrbracket^{1/(p-1)}<\infty\quad\hbox{ a.s.}
\]
This, together with Lemma~\ref{lem:weighted-cell-sums}, implies \eqref{eq:average-inverse-bound}.
\end{proof}

\begin{lemma}\label{lem:averaged-resistance}
There is a constant $C<\infty$ such that for every $k\geq 0$, 
\begin{equation}\label{eq:averaged-resistance}
 \limsup_{n\to\infty}\int_K
 d_n^\alpha(\pi_n(z),\pi_k(z))\,d\mu(z)
 \le C\rho^{-k/(p-1)}\quad  \hbox{ a.s.}
\end{equation}
\end{lemma}
\begin{proof}
Combining \eqref{eq:cell-chain-average} with
Lemma~\ref{lem:average-inverse-bound} gives the result.
\end{proof}

\subsection{Mosco convergence}\label{sec:mosco}
In this subsection, we prove the Mosco convergence of random
environments. See \cite{Mosco1969} and \cite{KuwaeShioya2003} for the definition of Mosco convergence in various settings. 

From now on, we view $\cE_0$ as a functional on
$L^t(K,\mu)$ for ${t\in[1,\infty)}$ by taking $\cE_0(u)=+\infty$ if $u$ does not have a representative in $\cF_0$.

\begin{theorem}[Mosco convergence]\label{thm:mosco}
Suppose that Assumption~\ref{ass:random} holds. Then, almost surely, for every
$\alpha\in\cW$ and $1\le t<\infty$, ${D_\alpha^{(n)}}$ is
Mosco convergent to $\cE_0$ in $L^t(K,\mu)$ in the following sense{:} 
\begin{enumerate}[label=\textup{(\arabic*)}]
\item If $\iota_nu_n\rightharpoonup u$ weakly in $L^t(K,\mu)$, where $u_n\in \mathbb{R}^{V_n}$, then
\begin{equation}\label{eq:M1}
 \cE_0(u)\le\liminf_{n\to\infty}{D_\alpha^{(n)}}(u_n).
\end{equation}
\item For every $u\in L^t(K,\mu)$, there are $u_n\in\R^{V_n}$
such that $\iota_nu_n\to u$ strongly in $L^t(K,\mu)$ and
\begin{equation}\label{eq:M2}
 \limsup_{n\to\infty}{D_\alpha^{(n)}}(u_n)\le\cE_0(u).
\end{equation}
{Such a sequence $u_n$ is called a recovery sequence.}
\end{enumerate}
\end{theorem}
\begin{remark}\label{remark:6.8}
We will see in the proof that the following statement holds.

\begin{enumerate}[label=\textup{(\arabic*)}]
	\item[(2')] {Almost surely, for every $\alpha\in\mathcal{W}$ and every $u\in\cF_0$, there are $u_n\in\R^{V_n}$
	such that $u_n|_{V_0}=u|_{V_0}$, $\iota_nu_n\to u$ uniformly and
	\begin{equation}\label{eq:M2'}
		\limsup_{n\to\infty}{D_\alpha^{(n)}}(u_n)\le\cE_0(u).
	\end{equation}}
\end{enumerate}
\end{remark}

{
Next we state a version of the preceding theorem on the infinite graphs $\mathcal{V}_n$.  Extend $\mu$ to $\cK$ by
\[
 \mu\mathrel:=\sum_{\alpha\in\cW}\mu\circ\Psi_\alpha^{-1},
\]
and define the measure on $\mathcal{V}_n$ by
\[
\nu_n\mathrel:=\sum_{\alpha\in\cW}\mu_n\circ\Psi_\alpha^{-1}.
\]

For $u\in L^t(\mathcal V_n,\nu_n)$, $1\le t<\infty$, extend $\iota_n$ by
\[
 (\iota_nu)\circ\Psi_\alpha=\iota_n(u\circ\Psi_\alpha),
 \qquad \alpha\in\cW.
\]
Then, by Proposition~\ref{prop:isometric-embedding}(a),
\[
 \|\iota_nu\|_{L^t(\cK,\mu)}^t
 =\sum_{\alpha\in\cW}\|\iota_n(u\circ\Psi_\alpha)\|_{L^t(K,\mu)}^t
 =\|u\|_{L^t(\mathcal V_n,\nu_n)}^t.
\]
Set
\[
 \cE(u)=\sum_{\alpha\in\cW}\cE_0(u\circ\Psi_\alpha)
 \quad\hbox{for }u\in C(\cK) \hbox{ and }
 u\circ\Psi_\alpha\in\cF_0, \forall\alpha\in\cW,
\]
with the convention $\cE(u)=+\infty$ otherwise.
Recall \eqref{eq:infinite-skeleton-energy} for the definition of ${\cD^{(n)}}$. We extend it to $L^t(\mathcal V_n,\nu_n),1\le t<\infty$ with the convention ${\cD^{(n)}}(u)=+\infty$ if the sum in \eqref{eq:infinite-skeleton-energy} diverges for some $u\in L^t(\mathcal V_n,\nu_n)$.

\begin{corollary}\label{cor:infinite-mosco}
Suppose that Assumption~\ref{ass:random} holds. Then, almost surely,
for every $1\le t<\infty$, ${\cD^{(n)}}$ is Mosco convergent to $\cE$ in
$L^t(\cK,\mu)$ in the following sense:
\begin{enumerate}[label=\textup{(\arabic*)}]
\item If $\iota_nu_n\rightharpoonup u$ weakly in $L^t(\cK,\mu)$, where $u_n\in L^t(\mathcal{V}_n,\nu_n)$, then
\[
 \cE(u)\le\liminf_{n\to\infty}{\cD^{(n)}}(u_n).
\]
\item For every $u\in L^t(\cK,\mu)$, there are
$u_n\in C_c(\mathcal V_n)$ such that
\[
 \|\iota_nu_n-u\|_{L^t(\cK,\mu)}\to0,\qquad
 \limsup_{n\to\infty}{\cD^{(n)}}(u_n)\le\cE(u).
\]
\end{enumerate}
\end{corollary}
}

To prove Theorem~\ref{thm:mosco} and Corollary~\ref{cor:infinite-mosco}, we first establish the following compactness lemma.

\begin{lemma}\label{lem:weak-moment-compactness}
Almost surely, {the following holds}. For any sequence $u_n\in\R^{V_n},n\geq 0$, if
\[
 \sup_n{D_\alpha^{(n)}}(u_n)<\infty\ \hbox{ and }\ \sup_n\|\iota_nu_n\|_{L^1(K,\mu)}<\infty,  
\]
then for any subsequence, there exist a further subsequence $u_{n_j},j\geq 1$ and a function $f\in\cF_0$ such that
\[
 \begin{aligned}
 &\|\iota_{n_j}u_{n_j}-f\|_{L^1(K,\mu)}\to 0,\\
 &{u_{n_j}|_{V_k}\to f|_{V_k}\text{ pointwise for every }k\ge0},\\
 &\cE_0(f)\le\liminf_{j\to\infty}{D_\alpha^{(n_j)}}(u_{n_j}).
 \end{aligned} 
\]
\end{lemma}
\begin{proof} 
In the following, we fix a sample environment so that  
\begin{equation}\label{e:6.12}
\|D_{n,\beta}-E_0\|\to0\hbox{ as }n\to\infty\quad\hbox{ for every }\beta\in\mathcal{W} 
\end{equation}
and the estimates of Lemma~\ref{lem:averaged-resistance} hold. By Theorem \ref{thm:random-homogenization} and Lemma \ref{lem:averaged-resistance}, it suffices to prove the claim for such a sample environment. \medskip 

\noindent
\textit{Step 1.} We prove that there is a constant $C>0$ so that for every $k\geq 0$, 
\begin{equation}\label{eq:compactness-approximation}
	\limsup_{n\to\infty}\int_K
	|\iota_nu_n(z)-\iota_ku_n(z)|\,d\mu(z)
	\le C\rho^{-k/p},
\end{equation}
where we write $\iota_ku_n:=\iota_k(u_n|_{V_k})$ for short. 

For each $k\ge0$ and $n>k$,
\[
 |\iota_nu_n(z)-{\iota_ku_n(z)}|=|u_n(\pi_n(z))-{u_n(\pi_k(z))}|
 \le{D_\alpha^{(n)}}(u_n)^{1/p}
 d_n^\alpha(\pi_n(z),\pi_k(z))^{(p-1)/p}.
\]
By Jensen's inequality,  
\[
 \begin{aligned}
 &\quad\ \int_K
 |\iota_nu_n(z)-\iota_ku_n(z)|\,d\mu(z)\\
 &\le{D_\alpha^{(n)}}(u_n)^{1/p} 
 \int_Kd_n^\alpha(\pi_n(z),\pi_k(z))^{(p-1)/p}\,d\mu(z)\\
 &\le{D_\alpha^{(n)}}(u_n)^{1/p}
 \left(\int_K
 d_n^\alpha(\pi_n(z),\pi_k(z))\,d\mu(z)\right)^{(p-1)/p}.
 \end{aligned}
\]
\eqref{eq:compactness-approximation} then follows {from} $\sup_n{D_\alpha^{(n)}}(u_n)<\infty$ and
Lemma~\ref{lem:averaged-resistance}. 
\medskip

\noindent\textit{Step 2.}
We construct a subsequence $u_{n_j},j\geq 1$ and $f\in\cF_0$ such that
${u_{n_j}|_{V_k}\to f|_{V_k}}$ pointwise {for every $k\ge0$}.

By \eqref{eq:compactness-approximation} and  $\sup_n\|\iota_nu_n\|_{L^1(K,\mu)}<\infty$, we get  
{
\begin{align*}
 \limsup_{n\to\infty}\|u_n|_{V_k}\|_{L^\infty(V_k,\mu_k)}
 &\le MN^k\limsup_{n\to\infty}\|\iota_ku_n\|_{L^1(K,\mu)}\\
 &\le MN^k\left(\sup_n\|\iota_nu_n\|_{L^1(K,\mu)}
 +\limsup_{n\to\infty}\|\iota_ku_n-\iota_nu_n\|_{L^1(K,\mu)}\right)\\
 &\le MN^k\left(\sup_n\|\iota_nu_n\|_{L^1(K,\mu)}
 +C\rho^{-k/p}\right)<\infty.
\end{align*}
}
Hence $ {\sup_{n\ge k}}\|u_n|_{V_k}\|_{L^\infty}<\infty$.  

Then, by a diagonal argument, for any subsequence, we can find a further subsequence $u_{n_j},j\geq 1$
and a function $f$ on ${\bigcup_{k\ge0}V_k}$ such that
\[
 u_{n_j}|_{V_k}\to f|_{V_k}\hbox{ as }j\to\infty\qquad\text{ for every }k\geq 0.
\]

As a consequence, by \eqref{e:6.12},  
\begin{equation}\label{eq:compactness-trace-limit}
 \cE_{0,k}(f|_{V_k})
 =\lim_{j\to\infty}(\Tr_{V_k}{D_\alpha^{(n_j)}})(u_{n_j}|_{V_k})
 \le\liminf_{j\to\infty}{D_\alpha^{(n_j)}}(u_{n_j}).
\end{equation}
Then, for $x,y\in V_k$, we have
\[
 \begin{aligned}
 |f(x)-f(y)|^p
 \le\cE_{0,k}(f|_{V_k})R_0(x,y)\lesssim R_0(x,y)=d_0(x,y)^{p-1}.
 \end{aligned}
\]
Since $d_0$ induces the original topology of $K$, $f$ extends uniquely
to a continuous function $f\in\cF_0$ with
\[
 \cE_0(f)=\sup_k\cE_{0,k}(f|_{V_k})
 \le\liminf_{j\to\infty}{D_\alpha^{(n_j)}}(u_{n_j}).
\]
\medskip

\noindent\textit{Step 3.}
It remains to show that $\iota_{n_j}u_{n_j}\to f$ strongly in
$L^1(K,\mu)$. Note that
\[
 \begin{aligned}
 \|\iota_{n_j}u_{n_j}-f\|_{L^1(K,\mu)}
 &\le
\|\iota_{n_j}u_{n_j}-\iota_ku_{n_j}\|_{L^1(K,\mu)}+\|\iota_ku_{n_j}-\iota_k(f|_{V_k})\|_{L^1(K,\mu)}
 +\max_{|w|=k}\operatorname{osc}_{\psi_w(K)}f.
 \end{aligned}
\]
By Step~2 and \eqref{eq:compactness-approximation},
\[
 \limsup_{j\to\infty}\|\iota_{n_j}u_{n_j}-f\|_{L^1(K,\mu)}
 \le C\rho^{-k/p}+\max_{|w|=k}\operatorname{osc}_{\psi_w(K)}f.
\]
Letting $k\to\infty$ gives the desired result.
\end{proof}

\begin{proof}[Proof of Theorem~\ref{thm:mosco}]
In the following, we fix a sample environment so that the event
\begin{equation}\label{e:thm6.7}
\|D_{n,\beta}-E_0\|\to0\hbox{ as }n\to\infty\quad\hbox{ for every }\beta\in\mathcal{W}
\end{equation}
and the event of {Lemma} \ref{lem:weak-moment-compactness} occur.\medskip 

We first prove (1). Suppose that $\liminf_{n\to\infty}{D_\alpha^{(n)}}(u_n)<\infty$;
otherwise there is nothing to prove. Pass to a subsequence along which
the energies converge to this lower limit. Since $\iota_nu_n{\rightharpoonup} u$ in $L^t$, 
\[
 \sup_n\|\iota_nu_n\|_{L^1(K,\mu)}
 \le\sup_n\|\iota_nu_n\|_{L^t(K,\mu)}
 <\infty.
\] 
Moreover, $\iota_nu_n\rightharpoonup u$ weakly in $L^1(K,\mu)$. 
By Lemma~\ref{lem:weak-moment-compactness}, there are
a further subsequence and $f\in\cF_0$ such that
\[
\|\iota_{n_j}u_{n_j}-f\|_{L^1}\to0\quad\hbox{ and }\quad 
\cE_0(f)\le\liminf_{j\to\infty}{D_\alpha^{(n_j)}}(u_{n_j})
=\liminf_{n\to\infty}{D_\alpha^{(n)}}(u_n).
\]
{Since $f=u$ by the uniqueness of the weak limit, (1) follows.}\medskip 

Next, we prove (2). First, suppose that $\cE_0(u)=+\infty$. Then, as $C(K)$ is dense in $L^t(K,\mu)$, and by Proposition~\ref{prop:isometric-embedding}(b), we can find $u_n\in \R^{V_n},\,n\geq 0$ so that $\iota_nu_n\to u$ in $L^t(K,\mu)$. \eqref{eq:M2} holds automatically, so (2) holds for this case.

Next suppose that $\cE_0(u)<+\infty$. By
\eqref{e:thm6.7}, we can choose strictly increasing
integers ${n_k\ge k},k\ge1$ such that
\[
 \left|(\Tr_{V_k}{D_\alpha^{(n)}})(u|_{V_k})
       -\cE_{0,k}(u|_{V_k})\right|
 \le k^{-1}\quad\hbox{ for every }n\ge n_k.
\]
For $n\ge n_1$, let $k(n):=\max\{k:n_k\le n\}$ and let $u_n$
be the harmonic extension of $u|_{V_{k(n)}}$ to $V_n$ with respect to
${D_\alpha^{(n)}}$; choose $u_n\in {\R^{V_n}}$ {with $u_n|_{V_0}=u|_{V_0}$} for $n<n_1$.
Then, for $n\ge n_1$,
\[
 {D_\alpha^{(n)}}(u_n)
 =(\Tr_{V_{k(n)}}{D_\alpha^{(n)}})(u|_{V_{k(n)}}).
\]
Therefore,
\[
 \left|{D_\alpha^{(n)}}(u_n)
       -\cE_{0,k(n)}(u|_{V_{k(n)}})\right|
 \le\frac1{k(n)}. 
\]
So
\[
 |{D_\alpha^{(n)}}(u_n)-\cE_0(u)|
 \le\frac1{k(n)}+\cE_0(u)-\cE_{0,k(n)}(u|_{V_{k(n)}})
 \to0\quad\hbox{ as }n\to\infty.
\]
Moreover, the maximum principle gives
\[
 \|\iota_nu_n-u\|_\infty
 \le\max_{|w|=k(n)}\operatorname{osc}_{\psi_w(K)}u
 \to0\quad\hbox{ as }n\to\infty. 
\]
Thus $\|\iota_nu_n-u\|_{L^t}\to0$. This proves (2).
\end{proof}

{
\begin{proof}[Proof of Corollary~\ref{cor:infinite-mosco}]
Fix an environment for which Theorem~\ref{thm:mosco},
Remark~\ref{remark:6.8} and Lemma~\ref{lem:weak-moment-compactness}
hold for every $\alpha\in\cW$.

\noindent\textit{Step 1.}
Suppose that $\iota_nu_n\rightharpoonup u$ in $L^t(\cK,\mu)$.
Passing to a subsequence, assume that
\[
 \lim_{n\to\infty}{\cD^{(n)}}(u_n)
 =\liminf_{n\to\infty}{\cD^{(n)}}(u_n)<\infty.
\]
On each $0$-cell $\Psi_\alpha(K)$, $\alpha\in\cW$,
\[
 \sup_n\|\iota_n(u_n\circ\Psi_\alpha)\|_{L^1(K,\mu)}
 \le\sup_n\|\iota_nu_n\|_{L^t(\cK,\mu)}<\infty.
\]
Lemma~\ref{lem:weak-moment-compactness} and a diagonal argument
give continuous representatives of $u$ on all $0$-cells and convergence
of a subsequence $u_{n_j}$ at every vertex of $\mathcal V_0$, so $u\in C(\cK)$.
By Theorem~\ref{thm:mosco}(1) and Fatou's lemma,
\[
 \begin{aligned}
 \cE(u)
 &=\sum_{\alpha\in\cW}\cE_0(u\circ\Psi_\alpha)\\
 &\le\sum_{\alpha\in\cW}
       \liminf_{n\to\infty}{D_\alpha^{(n)}}(u_n\circ\Psi_\alpha)\\
 &\le\liminf_{n\to\infty}
       \sum_{\alpha\in\cW}{D_\alpha^{(n)}}(u_n\circ\Psi_\alpha)\\
 &=\liminf_{n\to\infty}{\cD^{(n)}}(u_n).
 \end{aligned}
\]
This proves (1).\medskip

\noindent\textit{Step 2.}
Suppose first that $\cE(u)<\infty$. We approximate $u$ by functions supported on finitely many
$0$-cells.  Since $\mu(\Psi_\alpha(K))=1$ and $\operatorname{diam}_{d_0}(K)<\infty$,
\[
 \begin{aligned}
 \|u\|_{L^\infty(\Psi_\alpha(K))}
 &\le\int_{\Psi_\alpha(K)}|u|\,d\mu
      +\operatorname{osc}_{\Psi_\alpha(K)}u\\
 &\le\|u\|_{L^t(\Psi_\alpha(K),\mu)}
      +C\cE_0(u\circ\Psi_\alpha)^{1/p},
 \end{aligned}
\]
where $C=(\operatorname{diam}_{d_0}(K))^{(p-1)/p}$.
Note that
\[
 \sum_{\alpha\in\cW}\|u\|_{L^t(\Psi_\alpha(K),\mu)}^t
 =\|u\|_{L^t(\cK,\mu)}^t<\infty,\quad
 \sum_{\alpha\in\cW}\cE_0(u\circ\Psi_\alpha)=\cE(u)<\infty.
\]
Thus, for every $\varepsilon>0$,
\[
 \begin{aligned}
 \#\{\alpha:\|u\|_{L^\infty(\Psi_\alpha(K))}>\varepsilon\}
 &\le (2/\varepsilon)^t\|u\|_{L^t(\cK,\mu)}^t
      +(2C/\varepsilon)^p\cE(u)<\infty.
 \end{aligned}
\]
Put
\[
 u_\varepsilon=(u-\varepsilon)_+-(-u-\varepsilon)_+.
\]
Then $u_\varepsilon$ is nonzero on only finitely many $0$-cells.
The dominated convergence theorem and \textup{(RF5)} yield
\[
 \|u_\varepsilon-u\|_{L^t(\cK,\mu)}\to0,\qquad
 \cE(u_\varepsilon)\le\cE(u).
\]

For fixed \(\varepsilon\), apply the recovery sequence construction in the proof of Theorem 6.7 to \(u_\varepsilon\circ\Psi_\alpha\). These recovery sequences define
$u_{n,\varepsilon}\in C_c(\mathcal V_n)$, with
\[
 \begin{aligned}
 \|\iota_nu_{n,\varepsilon}-u_\varepsilon\|_{L^t(\cK,\mu)}^t
 &=\sum_{\alpha\in\cW}
   \|\iota_n(u_{n,\varepsilon}\circ\Psi_\alpha)
          -u_\varepsilon\circ\Psi_\alpha\|_{L^t(K,\mu)}^t
 \longrightarrow0,\\
 \limsup_{n\to\infty}{\cD^{(n)}}(u_{n,\varepsilon})
 &\le\sum_{\alpha\in\cW}\cE_0(u_\varepsilon\circ\Psi_\alpha)
 =\cE(u_\varepsilon)\le\cE(u).
 \end{aligned}
\]
Both sums have only finitely many nonzero terms. It is then routine to prove (2).

If $\cE(u)=+\infty$, by the density of $C_c(\cK)$ in $L^t(\cK,\mu)$
and Proposition~\ref{prop:isometric-embedding}(b), we can find $u_n\in C_c(\mathcal{V}_n)$ such that $\iota_nu_n\to u$ in $L^t(\cK,\mu)$.
\end{proof}
}

\section{Convergence of processes in random environments}
\label{sec:process-convergence}
From now on, we focus on the $p=2$ case. We assume that Theorem~\ref{thm:random-homogenization} holds and take $E_0=E_+$ for simplicity. Note that Assumption~\ref{ass:random} becomes $\mathbb E\llbracket D_{0,\alpha}\rrbracket<\infty$. Recall that 
\[
\cE(u):=\sum_{\alpha\in\cW}\cE_0(u\circ\Psi_\alpha)\quad\hbox{ for }u\in C(\cK).
\]
Let 
\[
\cF:=L^2(\cK,\mu)\cap\{u\in C(\cK):\,\cE(u)<\infty\},
\]
and set $\cE(u,v):=\frac14\big(\cE(u+v)-\cE(u-v)\big)$ when $u,v\in\cF${;} then $(\cE,\cF)$ is a strongly local, regular Dirichlet form on $L^2(\cK,\mu)$ (see \cite[Theorem 2.7]{FHK}). Let $\{X_t, t\geq 0; \mathbb P_x, x\in \cK\}$ be the Hunt process associated with $(\cE,\cF)$ on $L^2(\cK,\mu)$. The two-sided sub-Gaussian heat kernel estimates for $(\cE,\cF)$ follow from the estimates of Fitzsimmons, Hambly and Kumagai \cite[Theorem 1.1]{FHK} and the stability theorem \cite[Theorem 2.19]{BBK}. Further, by \cite[Theorem 3.1]{BGK}, the heat kernel of $X$ admits a jointly continuous version. The process $X$ can thus be chosen so that $\mathbb P_x$ is defined for every $x\in\cK$.
We denote the law of ${X}$ with starting point $x\in\cK$ by $\mathcal{L}_x(X)$, i.e. $\mathcal{L}_x(X):=\mathbb{P}_x\circ X^{-1}$.\medskip

We also consider the random environment $\cD^{(n)}=\cD^{(n)}(\omega)$ as introduced in Subsection \ref{sec:random-environment}. 
$(\cD^{(n)},\cF^{(n)})$ is a regular Dirichlet form on $L^2(\mathcal{V}_n,\nu_n)$, where $\cD^{(n)}(u,v):=\frac14\big(\cD^{(n)}(u+v)-\cD^{(n)}(u-v)\big)$ and 
\begin{align*}
\cF^{(n)}:=L^2(\mathcal{V}_n,\nu_n)\cap\{u\in \R^{\mathcal{V}_n}:\,\cD^{(n)}(u)<\infty\}.  
\end{align*}
 Let $\{X^{(n)}_t, t\geq 0;\mathbb P^{(n)}_x, x\in {\mathcal V}_n\}$ be the Hunt process associated with $(\cD^{(n)},\cF^{(n)})$ on $L^2({\mathcal V}_n,\nu_n)$. We denote the law of ${X^{(n)}}$ {with} starting point $x$ by ${\mathcal{L}_x(X^{(n)})}$, i.e. ${\mathcal{L}_x(X^{(n)})}:=\mathbb{P}_x^{(n)}\circ(X^{(n)})^{-1}$.\medskip 

By the scaling relations $\cD^{(n)}(\bullet)=\rho^n\cD^{(0)}(\bullet\circ{\psi_{i_0}^{n}})$ and $\nu_n(\bullet)=N^{-n}\nu_0(\psi_{i_0}^{-n}(\bullet))$,  we can couple $X^{(0)}$ and $X^{(n)}$ by 
\[
X^{(n)}_t=\psi_{i_0}^n(X^{(0)}_{(N\rho)^n t}). 
\]
The following lemma shows that $X^{(n)}$ does not blow up in finite time.

\begin{lemma}
$\mathbb{P}$-a.s., $(\cD^{(n)},\cF^{(n)})$ is a recurrent Dirichlet form on $L^2(\mathcal{V}_n,\nu_n)$, and hence $X^{(n)}_t$ is a recurrent process on $\mathcal{V}_n$ by \cite[Theorem 4.7.1]{FukushimaOshimaTakeda2011}. 
\end{lemma}
\begin{proof}
It suffices to prove the statement for $n=0$. For $k\geq 1$, we define $f_k\in \R^{V_k}$ such that $f_k|_{\psi_{i_0}(V_0)}=1$, $f_k|_{V_1\setminus\psi_{i_0}(V_0)}=0$ and $f_k$ is harmonic on $V_k\setminus V_1$ with respect to ${D_\emptyset^{(k)}}$. Then, we rescale $f_k$ to get  $u_k\in\R^{{\mathcal V}_0}$ by
\[
u_k(x)=\begin{cases}
f_k(\psi_{i_0}^k(x))\quad&\hbox{ if }x\in \psi_{i_0}^{-k}(V_k),\\
0\quad&\hbox{ if }x\in{\mathcal V}_0\setminus \psi_{i_0}^{-k}(V_k).
\end{cases}
\]
Then, $u_k\to{\bf 1}$ pointwise on ${\mathcal V}_0$. Moreover, as $\|D_{k-1,i}-E_0\|\to 0$ $\mathbb{P}$-a.s. for $1\leq i\leq N$, 
\begin{align*}
\cD^{(0)}(u_k)=\rho^{-k}\cD^{(k)}(u_k\circ\psi_{i_0}^{-k})=\rho^{-k}{D_\emptyset^{(k)}}(f_k)=\rho^{-k+1}\sum_{i\in I}D_{k-1,i}(f_k\circ\psi_i|_{V_0})\to 0\qquad\mathbb{P}\hbox{-}a.s.
\end{align*}
This proves the recurrent property of $(\cD^{(n)},\cF^{(n)})$ by \cite[Theorem 1.6.3]{FukushimaOshimaTakeda2011}.
\end{proof}

The main goal of this section is to prove the $\mathbb{P}$-a.s. weak convergence of $X^{(n)}$ to $X$.
 
\begin{theorem}\label{thm:process-homogenization}
Suppose that $\mathbb E\llbracket D_{0,\alpha}\rrbracket<\infty$. For every sequence $x_n\in\mathcal{V}_n,n\geq 0$ and $x\in \cK$ such that $x_n\to x\in\cK$,  
\begin{equation}\label{eq:quenched-process-limit}
\mathbb{P}\big(\mathcal{L}_{x_n}(X^{(n)})\Longrightarrow \mathcal{L}_x(X)\big)=1, 
\end{equation}
where $\Longrightarrow$ means the weak convergence of probability measures on $D([0,\infty),\cK)$ equipped with the Skorokhod $J_1$ topology for the Euclidean metric on $\cK$.
\end{theorem}
\begin{proof}
This will follow from Propositions \ref{prop:process-fdd} and \ref{lem:process-tightness}.
\end{proof}

\subsection{Convergence of finite-dimensional distributions}
\label{sec:process-fdd}
In this subsection, we prove the convergence of finite-dimensional distributions of $X^{(n)}$ {to those of} $X$. 

From now on, we let $(T_t)_{t\geq0}$ be the semigroup associated with $(\cE,\cF)$, and for $n\geq 0$, let $(T^{(n)}_t)_{t\geq0}$ be the semigroup associated with $(\cD^{(n)},\cF^{(n)})$. By {\cite[Theorem 4.2.3]{FukushimaOshimaTakeda2011}}, we have the following pointwise realization of $T_t$. More precisely, $T_tf=P_tf$ in $L^2(\cK,\mu)$, where
\begin{align*}
P_tf(x):=\mathbb{E}_xf(X_t)\quad\hbox{ for any Borel function }f\in L^2(\cK,\mu)\hbox{ and }x\in\cK. 
\end{align*}
Note that the joint continuity and semigroup property of the heat kernel give $P_tf\in C(\cK),\ t>0$.
Moreover, $T^{(n)}_tf=P^{(n)}_tf$ on ${\mathcal V}_n$, where we {set} 
\begin{align*}
P_t^{(n)}f(x):=\mathbb{E}^{(n)}_x f(X^{(n)}_t)\quad \hbox{ for }n\geq 0,\,f\in {L^2(\mathcal V_n,\nu_n)},\,x\in{\mathcal V}_n. 
\end{align*}
We note that $T^{(n)}_t$ converges to $T_t$ strongly in $L^2$ as a consequence of the Mosco convergence.  

\begin{lemma}\label{lem: semigroup}
\(\mathbb P\)-a.s., $\iota_nP^{(n)}_tf_n\to P_tf\text{ in }L^2(\cK,\mu)$ for every $t>0$, every sequence $f_n\in\R^{{\mathcal V}_n}$ and $f\in L^2(\cK,\mu)$ such that $\iota_nf_n\to f$ in $L^2(\cK,\mu)$.
\end{lemma}
\begin{proof}
This is an immediate consequence of	Proposition~\ref{prop:isometric-embedding}, {Corollary~\ref{cor:infinite-mosco}} and \cite[Theorem~2.4]{KuwaeShioya2003}.
\end{proof}

We need to strengthen the conclusion of Lemma \ref{lem: semigroup} to obtain {pointwise convergence}. To achieve this, we first establish some resistance estimates. 

\begin{lemma}\label{lem:process-upper-tail}
	Under the hypotheses of Proposition~\ref{lem:tree-averages},
	\begin{equation}\label{eq:tree-upper-tail}
		\sum_{m\ge0}\mathbb P\left(
		\sum_{v\in I^m}a_vY_v>\mathbb EY+t\right)
		\le\frac{2\mathbb EY}{(1-\theta)t}\quad\hbox{ for }t>0.
	\end{equation}
\end{lemma}
\begin{proof}
Let $T_m=\sum_{v\in I^m}a_vY_v\mathbf1_{\{a_vY_v\le t\}}$. By Lemma~\ref{lem:A2},
\[
\mathbb ET_m\le\mathbb EY,\qquad
\sum_{v\in I^*}\mathbb P(a_vY_v>t)\le\frac{\mathbb EY}{(1-\theta)t}, 
\qquad
\sum_m\operatorname{Var}(T_m)\le\frac{t\mathbb EY}{1-\theta}. 
\]
If $a_vY_v\le t$ for every $v\in I^m$, then $T_m=\sum_{v\in I^m}a_vY_v$. Since $\mathbb ET_m\le\mathbb EY$, it follows that
\[
\left\{\sum_{v\in I^m}a_vY_v>\mathbb EY+t\right\}
\subset
\bigcup_{v\in I^m}\{a_vY_v>t\}
\cup\{T_m-\mathbb ET_m>t\}.
\]
Combining this with the above bounds and Chebyshev's inequality, we obtain
\[
\begin{aligned}
\sum_{m\ge0}\mathbb P\left(\sum_{v\in I^m}a_vY_v>\mathbb EY+t\right)
	&\le\sum_{v\in I^*}\mathbb P(a_vY_v>t)
	+\sum_{m\ge0}\mathbb P(T_m-\mathbb ET_m>t)\\
	&\le\frac{\mathbb EY}{(1-\theta)t}
	+\frac1{t^2}\sum_{m\ge0}\operatorname{Var}(T_m)\\
	&\le\frac{2\mathbb EY}{(1-\theta)t}.
\end{aligned}
\]
\end{proof}

\begin{lemma}\label{lem:process-chain-resistance}
There is a constant $C>0$ depending only on $E_0$, $\cK$ and $\mathbb{E}\llbracket D_{0,\alpha}\rrbracket$ such that 
\begin{equation}\label{eq:prescribed-resistance}
\mathbb{P}\big(\limsup_{n\to\infty}{R_n^\alpha}\big(x_n,\pi_k(x_n)\big)
 \le C\rho^{-k/2}\big)=1
\end{equation}
for  any sequence $x_n\in V_n,n\geq 0$ and $k\geq 0$.
\end{lemma}
\begin{proof}
Since $R_n^\alpha=d_n^\alpha$ for $p=2$ and $\pi_n(x_n)=x_n$, the triangle inequality gives
\begin{equation}\label{e:lemma7.5-1}
R_n^\alpha\big(x_n,\pi_k(x_n)\big)
\le\sum_{j=k+1}^{n}R_n^\alpha\big(\pi_{j-1}(x_n),\pi_j(x_n)\big).
\end{equation}
For each $j$ with $k+1\le j\le n$, choose $w\in I^{j-1}$ so that $\{\pi_{j-1}(x_n),\pi_j(x_n)\}\subset\psi_w(V_1)$. Writing $m=n-j$, denote this word by $w^{(m,j)}$ and set
\begin{equation}\label{e:lemma7.5-2}
Z^{(m,j)}:=\sum_{i\in I}\llbracket D_{m,\alpha w^{(m,j)}i}\rrbracket.
\end{equation}
By Lemma~\ref{lem:cell-distance}(a),
\[
R_n^\alpha\big(\pi_{j-1}(x_n),\pi_j(x_n)\big)
\le C_1\rho^{-j}Z^{(n-j,j)}.
\]
Summing over $j$ gives
\[
R_n^\alpha\big(x_n,\pi_k(x_n)\big)
\le C_1\sum_{j=k+1}^n\rho^{-j}Z^{(n-j,j)}.
\]
In the following, we estimate each $Z^{(m,j)},m\geq0,j\geq k+1$.  \smallskip 

\noindent\textit{Claim 1}. There are $a_v\geq 0,v\in I^*$ satisfying {the weight conditions of Proposition~\ref{lem:tree-averages}}, and identically distributed random variables $Y^{(j)}_v,j\geq k+1,v\in I^*$  such that 
\begin{align}
&{\mathbb{E}Y^{(j)}_v=C_1N(M-1)\mathbb{E}}\llbracket D_{0,\alpha}\rrbracket;\label{e:lemma7.1-2}\\	
&\hbox{for each $m\geq 0,j\geq k+1$, it holds that }Y^{(j)}_v,\,v\in I^m\hbox{ are independent};\label{e:lemma7.1-3}\\
&Z^{(m,j)}\leq \sum_{v\in I^m}a_v Y^{(j)}_v\quad\hbox{ for }j\geq k+1,m\geq 0,\label{e:lemma7.1-4} 
\end{align}
{where $C_1$ is the constant in Lemma \ref{lem:3.5}.}
\begin{proof}[Proof of Claim 1]
For each $j\geq k+1$ and $v\in I^m,m\geq 0$, let
\[
Y^{(j)}_v={C_1(M-1)}\sum_{i\in I}\llbracket {D_{0,\alpha w^{(m,j)}iv}}\rrbracket,
\] 
where ${w^{(m,j)}}\in I^{j-1}$ is the word defining ${Z^{(m,j)}}$ in  \eqref{e:lemma7.5-2}. Then, \eqref{e:lemma7.1-2} and \eqref{e:lemma7.1-3} follow immediately. Next, we let $u_1,u_2,\cdots,u_{M-1}$ be the functions of Lemma~\ref{lem:3.5}, and set 
\[
a_v=\frac1{M-1}\sum_{s=1}^{M-1}a_v(u_s).
\]
{The conditions of Proposition~\ref{lem:tree-averages} follow} from {Corollary} \ref{coro:4.2}. By Lemma~\ref{lem:3.5}, 
\[
Z^{(m,j)}=\sum_{i\in I}\llbracket {D_{m,\alpha w^{(m,j)}i}}\rrbracket\leq {C_1}\sum_{i\in I}\sum_{1\leq s\leq M-1}\sum_{v\in I^m}a_v(u_s)\llbracket {D_{0,\alpha w^{(m,j)}iv}}\rrbracket=\sum_{v\in I^m}a_v Y^{(j)}_v,
\]
which is \eqref{e:lemma7.1-4}.
\end{proof}

By Lemma~\ref{lem:process-upper-tail} with $t=\rho^{j/2}$ and Claim 1, 
\[
 \begin{aligned}
 &\quad\ \sum_{j\ge k+1}\sum_{m\ge0}\mathbb P\left({Z^{(m,j)}}>{C_1N(M-1)}\mathbb E\llbracket D_{0,\alpha}\rrbracket+\rho^{j/2}\right)\\
 &\leq \sum_{j\ge k+1}\sum_{m\ge0}\mathbb P\left(\sum_{v\in I^m}a_vY^{(j)}_v>{C_1N(M-1)}\mathbb E\llbracket D_{0,\alpha}\rrbracket+\rho^{j/2}\right)\le C_2\sum_{j\ge1}\rho^{-j/2}\mathbb E\llbracket D_{0,\alpha}\rrbracket<\infty.
 \end{aligned}
\]
By the Borel--Cantelli lemma, $\mathbb{P}$-a.s.,  $Z^{(m,j)}>{C_1N(M-1)}\mathbb E\llbracket D_{0,\alpha}\rrbracket+\rho^{j/2}$ for only finitely many $m\geq 0,j\geq k+1$. Hence, 
\[
{R_n^\alpha}(x_n,\pi_k(x_n))
 \le C_1\sum_{j=k+1}^n\rho^{-j}({C_1N(M-1)}\mathbb E\llbracket D_{0,\alpha}\rrbracket+\rho^{j/2})
 \le C\rho^{-k/2}
\] 
for all sufficiently large $n$. 
\end{proof}

\begin{lemma}\label{lem:prescribed-evaluation}
Let $x_n\in V_n,n\geq 0$ and $x\in K$ such that $x_n\to x$. 
\(\mathbb P\)-a.s., for every sequence $u_n{\in}\R^{V_n},n\geq 0$,
\[
 \left\{
 \begin{aligned}
 &\sup_n {D_\alpha^{(n)}}(u_n)<\infty,\\
 &\iota_nu_n\to u\quad\text{in }L^2(K,\mu)
 \end{aligned}\right.
 \quad\Longrightarrow\quad
 \left\{
 \begin{aligned}
 &u\in C(K),\\
 &u_n(x_n)\to u(x).
 \end{aligned}\right.
\]
\end{lemma}
\begin{proof}
{We work on a sample environment for which Lemma~\ref{lem:process-chain-resistance} and  Lemma~\ref{lem:weak-moment-compactness} hold.} By Lemma~\ref{lem:weak-moment-compactness}, $u\in\cF_0\subset C(K)$, and every subsequence admits a further subsequence $n_j,j\geq 1$ such that $u_{n_j}|_{V_k}\to u|_{V_k}$ for every fixed $k\geq 0$. 
On this subsequence, 
\[
 \begin{aligned}
&\quad\  |u_{n_j}(x_{n_j})-u(x)|\\
 &\le |u_{n_j}(x_{n_j})-u_{n_j}(\pi_k(x_{n_j}))|
       +|u_{n_j}(\pi_k(x_{n_j}))-u(\pi_k(x_{n_j}))|+|u(\pi_k(x_{n_j}))-u(x)|\\
 &\le\sqrt{{D_\alpha^{(n_j)}}(u_{n_j})R_{n_j}^\alpha(x_{n_j},{\pi_k(x_{n_j})})}
       +\max_{z\in V_k}|u_{n_j}(z)-u(z)|+|u({\pi_k(x_{n_j})})-u(x)|.
 \end{aligned}
\]
Since $|{\pi_k(x_{n_j})}-x_{n_j}|\le r^k\operatorname{diam}K$,
Lemma~\ref{lem:process-chain-resistance} yields
\[
 \limsup_{j\to\infty}|u_{n_j}(x_{n_j})-u(x)|
 \le C\left(\sup_{n_j}{D_\alpha^{(n_j)}}(u_{n_j})\right)^{1/2}\rho^{-k/4}
 +\sup_{w\in I^k}\operatorname{osc}_{\psi_wK}u.
\]
Letting $k\to\infty$, we see $u_{n_j}(x_{n_j})\to u(x)$. Since the argument works for any subsequence, $u_n(x_n)\to u(x)$. 
\end{proof}

Applying Lemma {\ref{lem:prescribed-evaluation}} on each $0$-cell gives the same conclusion on the infinite fractal $\cK$. 
\begin{corollary}\label{coro:prescribed-evaluation}
Let $x_n\in {\mathcal V}_n,n\geq 0$ and $x\in \cK$ such that $x_n\to x$ as $n\to\infty$. 
\(\mathbb P\)-a.s., for every sequence $u_n\in{\R^{\mathcal V_n}},n\geq 0$,
\[
\left\{
\begin{aligned}
	&\sup_n\mathcal D^{(n)}(u_n)<\infty,\\
	&\iota_nu_n\to u\quad\text{in }L^2(\cK,\mu)
\end{aligned}\right.
\quad\Longrightarrow\quad
\left\{
\begin{aligned}
	&u\in C(\cK),\\
	&u_n(x_n)\to u(x).
\end{aligned}\right.
\]
\end{corollary}

\begin{lemma}\label{lem:process-semigroup}
Let $x_n\in {\mathcal V}_n,n\geq 0$ and $x\in\cK$ be such that $x_n\to x\in {\cK}$.
Then, \(\mathbb P\)-a.s., for every sequence $f_n\in {\mathbb{R}^{\mathcal V_n}}$ and $t>0$,
\[
\iota_nf_n\to f\text{ in }{L^2(\cK,\mu)}\quad\Longrightarrow\quad
P_t^{(n)}f_n(x_n)\to P_tf(x). 
\]
\end{lemma}
\begin{proof}
{We work on a sample environment for which Lemma~\ref{lem: semigroup} and Corollary~\ref{coro:prescribed-evaluation} hold.} By Lemma \ref{lem: semigroup}, $\iota_nP_t^{(n)}f_n\to P_tf$ in $L^2(\cK,\mu)$. Moreover, a routine spectral theorem argument {gives}
\begin{equation}\label{eq:semigroup-energy-bound}
	\sup_{n\geq 1}{\mathcal D^{(n)}}(P_t^{(n)}f_n)
	\le \sup_{n\geq 1}(2et)^{-1}\|{\iota_n} f_n\|_{L^2(\cK,\mu)}^2<\infty.
\end{equation}
The lemma follows from the above two observations by using Corollary~\ref{coro:prescribed-evaluation}.
\end{proof}

Finally, we can prove the convergence of {the finite-dimensional laws}. We write $	{\mathcal L_{x_n}}
\bigl(X^{(n)}_{t_1},\ldots,X^{(n)}_{t_\ell}\bigr)$ for the finite-dimensional law of $X^{(n)}$ and write $\mathcal L_x\bigl(X_{t_1},\ldots,X_{t_\ell}\bigr)$ for the finite-dimensional law of $X$. 
\begin{proposition}\label{prop:process-fdd}
	Let $x_n\in {\mathcal V_n}$ with $x_n\to x\in {\cK}$.
	\(\mathbb P\)-a.s., we have for every $\ell\ge1$ and
	$0\le t_1<\cdots<t_\ell<\infty$,
	\[
	\mathcal L_{x_n}
	\bigl(X^{(n)}_{t_1},\ldots,X^{(n)}_{t_\ell}\bigr)
	\Longrightarrow
	\mathcal L_x\bigl(X_{t_1},\ldots,X_{t_\ell}\bigr).
	\]
\end{proposition}
\begin{proof}
We fix {a} sample environment $\omega$ so that the event of Lemma~\ref{lem:process-semigroup} occurs, and prove the statement for {the} $\ell=2$ case. The case $\ell=1$ follows by taking a marginal, and the cases $\ell>2$ follow by induction. 

 For $u\in C_c(\cK)$, we write $\iota_n u=\iota_n(u|_{\mathcal {V}_n})$ and $P^{(n)}_tu=P^{(n)}_t(u|_{\mathcal{V}_n})$. Fix $f,g\in C_c(\cK)$ and $s,t>0$. By Lemma \ref{lem: semigroup}, 
\[
\iota_n(P^{(n)}_t g)\to P_tg\qquad \hbox{ in }L^2(\cK,\mu).
\]
Also noticing that $\|\iota_nf-f\|_\infty\to 0$ as $n\to\infty$, we have 
\begin{align*}
&\quad\ \|\iota_n((f|_{\mathcal{V}_n})P^{(n)}_t g){-} fP_tg\|_{L^2(\cK,\mu)}\\
&=\|(\iota_nf)(\iota_n P^{(n)}_tg)-fP_tg\|_{L^2(\cK,\mu)}\\
&\leq \|(\iota_nf)(\iota_n P^{(n)}_t g)-(\iota_nf)P_tg\|_{L^2(\cK,\mu)}+\|(\iota_nf)P_tg-fP_tg\|_{L^2(\cK,\mu)}\\
&\leq \|f\|_{L^\infty}\|\iota_n P^{(n)}_t g-P_tg\|_{L^2(\cK,\mu)}+\|\iota_nf-f\|_{L^\infty}\|P_tg\|_{L^2(\cK,\mu)}\\
&\to 0\quad\hbox{ as }n\to\infty. 
\end{align*}
That is{,} $\iota_n((f|_{\mathcal{V}_n})P^{(n)}_tg)\to fP_tg$ in $ L^2(\cK,\mu)$ as $n\to\infty$. Then, by {Lemma~\ref{lem:process-semigroup}},
\[
\begin{aligned}
\mathbb E^{(n)}_{x_n}\bigl[f(X^{(n)}_s)g(X^{(n)}_{s+t})\bigr]=P_s^{(n)}((f|_{\mathcal{V}_n})P_t^{(n)}g)(x_n)
\to P_s(fP_tg)(x)=\mathbb E_x\bigl[f(X_s)g(X_{s+t})\bigr]. 
\end{aligned}
\]
This proves $
\mathcal L_{x_n}
\bigl(X^{(n)}_s,X^{(n)}_{s+t}\bigr)
\Longrightarrow
\mathcal L_x\bigl(X_s,X_{s+t}\bigr).
$
\end{proof}

\subsection{Tightness}
\label{sec:process-tightness}
Next, we prove tightness of the laws of $X^{(n)}$. Recall the following Aldous tightness criterion. 

\begin{lemma}\label{lem:process-stopping-times}
Fix a sample point $\omega$ and let $x_n\in {\mathcal V}_n,n\geq 0$. 
Suppose that $\{x_n,n\geq 0\}$ is bounded, and for every $T,\varepsilon>0$,
{there is a compact set $F\subset\cK$ such that
\begin{equation}\label{compact containing condition}
 \limsup_{n\to\infty}\mathbb P_{x_n}^{(n)}
 \bigl(X_s^{(n)}\notin F\text{ for some }s\in[0,T]\bigr)
 <\varepsilon,
 \end{equation}
and}
\begin{equation}\label{e:Aldous}
\lim_{\delta\to0}\limsup_{n\to\infty}\sup_{\tau\le T}\sup_{0\le h\le\delta}\mathbb P^{(n)}_{x_n}\left(|X^{(n)}_{\tau+h}-X^{(n)}_\tau|>\varepsilon\right)=0,
\end{equation}
where $\tau$ ranges over stopping times for the natural filtration $\sigma(X_s^{(n)}:s\in[0,t])$. Then ${\mathcal L_{x_n}}(X^{(n)}),n\geq 0$
is tight in $D([0,\infty),{\cK})$ with the $J_1$ topology.

\end{lemma}
\begin{proof}
This follows from \cite[Chapter~3, Theorem~8.6,
Remark~8.7(b) and Theorem~10.2(a)]{EthierKurtz1986}.
\end{proof}

For $B\subset \cK$ and a stochastic process $Z_t$, we write
\[
\sigma_B=\inf\{s>0:Z(s)\in B\}\quad\hbox{ and }\quad\tau_B=\inf\{s\geq 0:Z(s)\notin B\}
\]
for the first hitting time and exit time.

\begin{lemma}\label{lem:process-crossing}
\(\mathbb P\)-a.s., the following events occur.
\begin{enumerate}
\item[(a)]   For every integer $k\geq 0$ and $z\in{\mathcal V}_k$,
\[
\lim_{\delta\to 0}\limsup_{n\to\infty}\mathbb P^{(n)}_z(\sigma_{{\mathcal V}_k\setminus\{z\}}\le\delta)=0.
\]
\item[(b)] For each $j\geq 0$, let $B_j=\psi^{-j}_{i_0}(K)\setminus\psi_{i_0}^{-j}(V_0\setminus\{q_{i_0}\})$. Then, for every $z\in {\bigcup_{\ell\ge0}\mathcal V_\ell}$ and $T>0$, 
\[
\lim_{j\to\infty}\limsup_{n\to\infty}\mathbb P^{(n)}_z(\tau_{B_j}\le T)=0.
\]
\end{enumerate}
\end{lemma}
\begin{proof}
Fix a sample environment so that Theorem~\ref{thm:mosco}, Remark~\ref{remark:6.8}, Corollary~\ref{cor:infinite-mosco} and Lemma~\ref{lem:weak-moment-compactness} hold. Fix $k\geq 0$ and $z\in {\mathcal V}_k$. Choose $V\subset{\mathcal V}_k\setminus\{z\}$ so that ${\mathcal V}_k\setminus V$ is a finite set. Let $f\in\R^{V}$ be a finitely supported function. 

For $n\geq k$ and $\lambda>0$, we define $h^{(n)}_\lambda\in \R^{{\mathcal V}_n}$ by 
\[
h^{(n)}_\lambda(y)=\mathbb E^{(n)}_y[f(X^{(n)}_{\sigma_{V}})e^{-\lambda\sigma_{V}}]\quad\hbox{ for }y\in {\mathcal V}_n.
\]
By \cite[Theorem 4.3.1]{FukushimaOshimaTakeda2011}, $h^{(n)}_\lambda$ is the harmonic extension of $f$ with respect to the Dirichlet form $(\cD^{(n)})_\lambda$, which is defined by 
\[
(\cD^{(n)})_\lambda(u,v):=\cD^{(n)}(u,v)+\lambda\int_{{\mathcal V}_n}uvd\nu_n\quad\hbox{ for }u,v\in \cF^{(n)}.
\]

For $\lambda>0$, we define $h_\lambda\in\cF$ by 
\[
h_\lambda(y)=\mathbb E_y[f(X_{\sigma_V})e^{-\lambda\sigma_V}]\quad\hbox{ for }y\in \cK. 
\]
$h_\lambda$ is the harmonic extension of $f$ with respect to the Dirichlet form $(\cE)_\lambda$, which is defined by
\[
(\cE)_\lambda(u,v):=\cE(u,v)+\lambda\int_{\cK}uvd\mu\quad\hbox{ for }u,v\in \cF.
\]

\noindent\textit{Claim 1.} For any $\lambda>0$, $\lim\limits_{n\to\infty}h^{(n)}_\lambda(z)=h_\lambda(z)$.
\begin{proof}[Proof of Claim 1] Note that $h^{(n)}_\lambda,n\geq k$ and $h_\lambda$ are all compactly supported in a common neighborhood of $z$ in $\cK$.

{On each of the finitely many $0$-cells meeting the common support, take the recovery sequence from the proof of Theorem~\ref{thm:mosco}.} This gives
a sequence $v^{(n)}\in\cF^{(n)}$ such that
$v^{(n)}|_{\mathcal V_k}=h_\lambda|_{\mathcal V_k}$ {for all sufficiently large $n$},
$\iota_nv^{(n)}\to h_\lambda$ in $L^2(\cK,\mu)$, and
$\limsup\limits_{n\to\infty}\cD^{(n)}(v^{(n)})\leq\cE(h_\lambda)$. Hence,
\begin{equation}\label{e:lemma7.11b}
\begin{split}
\limsup\limits_{n\to\infty}(\cD^{(n)})_\lambda(h^{(n)}_\lambda)&\leq \limsup\limits_{n\to\infty}(\cD^{(n)})_\lambda(v^{(n)})\\
&= \limsup\limits_{n\to\infty}\Big(\cD^{(n)}(v^{(n)})+\lambda\|v^{(n)}\|_{L^2(\mathcal V_n,\nu_n)}^2\Big)\\
&\leq \cE(h_\lambda)+\lambda\|h_\lambda\|_{L^2(\cK,\mu)}^2=(\cE)_\lambda(h_\lambda),
\end{split}
\end{equation}
where the first inequality uses the harmonicity of $h^{(n)}_\lambda$.
{In particular,
\[
 \sup_n\cD^{(n)}(h_\lambda^{(n)})<\infty.
\]
Also note that $|h_\lambda^{(n)}|\le\|f\|_\infty$.} 

Therefore, by Lemma~\ref{lem:weak-moment-compactness}, from any subsequence we obtain a further subsequence $n_j,j\geq 1$ and $h\in\cF$ such that
\begin{equation}\label{e:lemma7.11c}
\begin{aligned}
 &\iota_{n_j}h^{(n_j)}_\lambda\longrightarrow h\quad\text{in }L^2(\cK,\mu),\\
 &h^{(n_j)}_\lambda|_{\mathcal V_\ell}\to h|_{\mathcal V_\ell}
 \quad\text{pointwise for every fixed }\ell\ge0,\qquad h|_V=f.
\end{aligned}
\end{equation}
{Here the $L^1$ convergence given by Lemma~\ref{lem:weak-moment-compactness} implies $L^2$ convergence, since
\[
 \|\iota_{n_j}h^{(n_j)}_\lambda-h\|_{L^2(\cK,\mu)}^2
 \le 2\|f\|_\infty\,\|\iota_{n_j}h^{(n_j)}_\lambda-h\|_{L^1(\cK,\mu)} .
\]}
Hence, by Mosco convergence
\begin{equation}\label{e:lemma7.11a}
\begin{split}
\liminf\limits_{j\to\infty}(\cD^{(n_j)})_\lambda(h^{(n_j)}_\lambda)&= \liminf\limits_{j\to\infty}\Big(\cD^{(n_j)}(h^{(n_j)}_\lambda)+\lambda\|h^{(n_j)}_\lambda\|_{L^2(\mathcal V_{n_j},\nu_{n_j})}^2\Big)\\
&\geq \cE(h)+\lambda\|h\|_{L^2(\cK,\mu)}^2=(\cE)_\lambda(h).
\end{split}
\end{equation}

Now, we have $(\cE)_\lambda(h)\leq(\cE)_\lambda(h_\lambda)$ by \eqref{e:lemma7.11a} and \eqref{e:lemma7.11b}, and $h|_V=f$ by \eqref{e:lemma7.11c}. Hence, $h=h_\lambda$ by the uniqueness of harmonic extension. It follows that $h_\lambda^{(n_j)}(z)\to h_\lambda(z)$ by \eqref{e:lemma7.11c}.

The argument works for any subsequence, so the claim holds.
\end{proof}

(a).  We choose $V={\mathcal V}_k\setminus\{z\}$ and finitely supported $f\in\R^V$ satisfying $f|_{{\psi_{i_0}^k\Psi_w(V_0)\setminus\{z\}}}=1$ for each $w\in\cW$ such that $z\in\psi_{i_0}^k\Psi_w(V_0)$. Then, since 
\[
h_\lambda(z)=\mathbb E_z[f(X_{\sigma_{{\mathcal V}_k\setminus\{z\}}})e^{-\lambda\sigma_{{\mathcal V}_k\setminus\{z\}}}]\ \hbox{ and }\ h^{(n)}_\lambda(z)=\mathbb E^{(n)}_z[f({X^{(n)}_{\sigma_{{\mathcal V}_k\setminus\{z\}}}})e^{-\lambda\sigma_{{\mathcal V}_k\setminus\{z\}}}]\hbox{ for }n\geq k,
\] 
we have
\[
\mathbb P^{(n)}_z(\sigma_{{\mathcal V}_k\setminus \{z\}}\le\delta)\le e^{\lambda\delta}h^{(n)}_\lambda(z)\quad\hbox{ for }\lambda>0,\delta>0,n\geq k,
\]
and $h_\lambda(z)=\mathbb{E}_z[e^{-\lambda\sigma_{{\mathcal V}_k\setminus\{z\}}}]\to 0$ as $\lambda\to\infty$. By Claim 1,   
\[
\lim_{\delta\to0}\limsup_{n\to\infty}\mathbb{P}_{{z}}^{(n)}(\sigma_{{\mathcal V}_k\setminus\{z\}}\le\delta)\leq \lim_{\delta\to0}\limsup_{n\to\infty}e\cdot h^{(n)}_{1/\delta}(z)=\lim_{\delta\to0}e\cdot h_{1/\delta}(z)=0. 
\]

(b). For {sufficiently large} $j$, we choose $V={\mathcal V}_k\setminus B_j$ and finitely supported $f\in\R^V$ so that $f|_{\psi_{i_0}^{-j}(V_0\setminus\{q_{i_0}\})}=1$. Then, by {Claim} 1, 
\[
\limsup_{n\to\infty}\mathbb P^{(n)}_z(\tau_{B_j}\le T){\le}\limsup_{n\to\infty}e^{T}h^{(n)}_1(z)=e^Th_1(z)=e^T\mathbb{E}_z[e^{-\tau_{B_j}}]. 
\]
(b) follows immediately. 
\end{proof}

\begin{proposition}\label{lem:process-tightness}
Let $x_n\in {\mathcal V_n},n\geq 0$ be a bounded sequence of points. Then
${\mathcal L_{x_n}}(X^{(n)}),\,n\ge0$ is tight in
$D([0,\infty),{\cK})$, \(\mathbb P\)-a.s.
\end{proposition}
\begin{proof}
Fix a sample environment so that Lemma~\ref{lem:process-crossing} holds. Let $B_j$ be as in Lemma~\ref{lem:process-crossing}(b). Write $B(x,\varepsilon):=\{y\in\cK:|x-y|\leq \varepsilon\}$ for the closed ball centered at $x$ with radius $\varepsilon$. 

For each $\varepsilon>0$ and $j\geq 0$, by choosing {$\ell$ large enough} so that $2r^\ell\operatorname{diam}(K)\leq\varepsilon$, we see 
\begin{equation}\label{e:tightness1}
\begin{split}
\lim_{\delta\to 0}\limsup_{n\to\infty}\sup_{{x\in B_j\cap\mathcal V_n}}\sup_{0\leq h\leq \delta}\mathbb{P}^{(n)}_x(|X^{(n)}_h-x|>\varepsilon)&\leq \lim_{\delta\to 0}\limsup_{n\to\infty}\sup_{{x\in B_j\cap\mathcal V_n}}\mathbb{P}^{(n)}_x(\tau_{B(x,\varepsilon)}\leq\delta)\\
&\leq \lim_{\delta\to 0}\limsup_{n\to\infty}\sup_{{z\in\mathcal V_\ell\cap(B_j\cup\psi_{i_0}^{-j}(V_0))}}\mathbb{P}^{(n)}_z(\sigma_{{\mathcal V}_\ell\setminus\{z\}}\leq\delta)\\
&=0 
\end{split}
\end{equation}
by the strong Markov property and Lemma \ref{lem:process-crossing}(a).  
Moreover, by the strong Markov property and Lemma \ref{lem:process-crossing}(b), for every $i\ge 0$,
\[
\limsup_{n\to\infty}\sup_{{x\in B_i\cap\mathcal V_n}}\mathbb{P}^{(n)}_x(\tau_{B_j}<T)\to 0\hbox{ as }j\to\infty.  
\]
 \eqref{compact containing condition} and \eqref{e:Aldous} follow from the above observations. {We finish the proof by Lemma~\ref{lem:process-stopping-times}.}
\end{proof}

\appendix
\renewcommand{\appendixname}{Appendix}

\section{A law of large numbers}\label{AppendixA}
In this appendix, we prove a strong law of large numbers for random variables with indices in $\bigcup_{n\ge0}I^n=:I^*$.\medskip

\begin{proposition}\label{lem:tree-averages}
Let $a_w\ge0, w\in I^*$. Assume   
\[
a_\emptyset=1,\qquad \sum_{i\in I}a_{wi}=a_w,\qquad
a_{wi}\le\theta a_w\quad\hbox{ for each }w\in I^*\hbox{ and }i\in I,
\]
for some $\theta\in(0,1)$.  Let $Y_w$, $n\ge0$, $w\in I^n$, be random
	variables with the same law as a nonnegative integrable random variable
	$Y$. Suppose further that, for each $n\geq 0$, the variables $Y_w$, $w\in I^n$, are independent. Then for every $\varepsilon>0$,
	\begin{equation}\label{eq:A.1}
		\sum_{n=0}^\infty\mathbb P\left(
		\left|\sum_{w\in I^n}a_wY_w-\mathbb EY\right|>\varepsilon\right)
		<\infty,
	\end{equation}
	and so $\sum_{w\in I^n}a_wY_w\to\mathbb EY$ almost surely and in
	$L^1$.
\end{proposition}

\begin{lemma}\label{lem:A2}
Let $a_w\ge0, w\in I^*$. Assume   
\[
a_\emptyset=1,\qquad \sum_{i\in I}a_{wi}=a_w,\qquad
a_{wi}\le\theta a_w\quad\hbox{ for each }w\in I^*\hbox{ and }i\in I,
\]
for some $\theta\in(0,1)$. Then, for every $x>0$,
\begin{enumerate}[label=(\alph*)]
\item $\sum_{w\in I^*}\mathbf1_{\{a_wx>1\}}\le\frac{x}{1-\theta}$.
\item $\sum_{w\in I^*}a_w^2\mathbf1_{\{a_wx\le1\}}
\le\frac{\min\{1,x^{-1}\}}{1-\theta}$. 
\end{enumerate}	
\end{lemma}
\begin{proof}
Let $(W_n)_{n\ge0}$ be a Markov chain on the
state space $\{w\in I^*:a_w>0\}$, with initial state $W_0=\emptyset$ and transition probabilities $P_{w,wi}=a_{wi}/a_w$.
Then $\mathbb P(W_n=w)=a_w$ for $w\in I^n$.

(a). If $0<x\le1$, the left-hand side is zero. For $x>1$, let
\[
\tau=\max\{n\geq 0:\,a_{W_n}x>1\},
\]
which is finite because $a_{W_n}\le\theta^n$. Then,
\begin{align*}
\sum_{w\in I^*}\mathbf1_{\{a_wx>1\}}&=\sum_{\substack{w\in I^*\\a_w>0}}a_w\cdot a_w^{-1}\mathbf1_{\{a_wx>1\}}\\
&=\mathbb E\sum_{n\ge0}a_{W_n}^{-1}
	\mathbf1_{\{a_{W_n}x>1\}}
	=\mathbb E\sum_{n=0}^{\tau}a_{W_n}^{-1}\\
&=\mathbb E\sum_{n=0}^{\tau}a_{W_{\tau-n}}^{-1}\le\mathbb E\sum_{n=0}^{\infty}x\theta^n
	=\frac{x}{1-\theta}.
\end{align*}

(b). Let 
\[
\sigma=\min\{n\geq 0:\,a_{W_n}x\leq 1\}.
\]
Then,
\begin{align*}
\sum_{w\in I^*}a_w^2\mathbf1_{\{a_wx\le1\}}&=\sum_{w\in I^*}a_w\cdot a_w\mathbf1_{\{a_wx\le1\}}\\
&=\mathbb E\sum_{n\ge0}a_{W_n}
	\mathbf1_{\{a_{W_n}x\le1\}}
	=\mathbb E\sum_{n\ge \sigma}a_{W_n}\\
&\le\mathbb E\sum_{n\ge \sigma}(1\wedge x^{-1})\theta^{n-\sigma}
	=\frac{\min\{1,x^{-1}\}}{1-\theta}.
\end{align*}
\end{proof}

\begin{proof}[Proof of Proposition \ref{lem:tree-averages}]
We use the truncation method. Let 
\[
T_n:=\sum_{w\in I^n}a_wY_w{\bf 1}_{\{a_wY_w\leq 1\}}\ \hbox{ for }n\ge 0.
\]
Then,
\begin{align} 
	&\qquad\ \sum_{n=0}^\infty\mathbb P\left(
	\left|\sum_{w\in I^n}a_wY_w-\mathbb EY\right|>\varepsilon\right)\nonumber\\
	&\quad\le\sum_{n=0}^\infty\left[
	\mathbb P\left(\left|\sum_{w\in I^n}a_wY_w-\mathbb ET_n\right|
	>\frac\varepsilon2\right)
	+\mathbb P\left(|\mathbb EY-\mathbb ET_n|>\frac\varepsilon2\right)
	\right]\nonumber\\
	&\quad\le\sum_{n=0}^\infty\left[\mathbb P\left(T_n\ne\sum_{w\in I^n}a_wY_w\right)+
	\mathbb P\left(|T_n-\mathbb ET_n|>\frac\varepsilon2\right)
	+\mathbb P\left(|\mathbb EY-\mathbb ET_n|>\frac\varepsilon2\right)
	\right].\label{eqnA2}
\end{align}

The series of first terms in \eqref{eqnA2} converges, as shown in the following calculation:
	\begin{align}
		\sum_{n\ge0}\mathbb P\left(T_n\ne\sum_{w\in I^n}a_wY_w\right)
		&\le\sum_{n\ge0}\sum_{w\in I^n}\mathbb P(a_wY_w>1)\notag\\
		&=\sum_{w\in I^*}\mathbb P(a_wY_w>1)\notag\\
		&=\sum_{w\in I^*}\mathbb P(a_wY>1)\notag\\
		&=\mathbb E\left(\sum_{w\in I^*}\mathbf1_{\{a_wY>1\}}\right)\\
		&\leq (1-\theta)^{-1}\mathbb{E}Y,
		\label{eq:A.3}
	\end{align}
where the last inequality is due to Lemma \ref{lem:A2}(a). 

The second term of \eqref{eqnA2} is dealt with by the classical second-moment estimate.
	\begin{align}
		\sum_{n\ge0}\operatorname{Var}(T_n)
		&=\sum_{n\ge0}\operatorname{Var}\left(
		\sum_{w\in I^n}a_wY_w\mathbf1_{\{a_wY_w\le1\}}\right)\notag\\
		&=\sum_{n\ge0}\sum_{\substack{w\in I^n}}a_w^2
		\operatorname{Var}\bigl(Y_w\mathbf1_{\{a_wY_w\le1\}}\bigr)\notag\\
		&=\sum_{n\ge0}\sum_{\substack{w\in I^n}}a_w^2
		\operatorname{Var}\bigl(Y\mathbf1_{\{a_wY\le1\}}\bigr)\notag\\
		&\le\sum_{n\ge0}\sum_{\substack{w\in I^n}}a_w^2
		\mathbb E\bigl(Y^2\mathbf1_{\{a_wY\le1\}}\bigr)\notag\\
		&=\mathbb E\left(Y^2\sum_{w\in I^*}a_w^2
		\mathbf1_{\{a_wY\le1\}}\right)\\
		&\leq (1-\theta)^{-1}\mathbb{E}Y,
	\end{align}
	where the last inequality uses Lemma \ref{lem:A2}(b). 
	Chebyshev's inequality gives
	\begin{equation}\label{eq:A.4}
		\sum_{n\ge0}\mathbb P(|T_n-\mathbb ET_n|>\frac{\varepsilon}{2})
		\le\sum_{n\ge0}4\varepsilon^{-2}\operatorname{Var}(T_n)<\infty.
	\end{equation}
	
	Finally, for the third term of \eqref{eqnA2},
	\begin{align}
		0\le\mathbb EY-\mathbb ET_n
		&=\sum_{w\in I^n}a_w\mathbb EY-\sum_{w\in I^n}a_w
		\mathbb E(Y\mathbf1_{\{a_wY\le1\}})\notag\\
		&=\sum_{w\in I^n}a_w\mathbb E(Y\mathbf1_{\{a_wY>1\}})\notag\\
		&\le\sum_{w\in I^n}a_w\mathbb E(Y\mathbf1_{\{Y>\theta^{-n}\}})
		=\mathbb E(Y\mathbf1_{\{Y>\theta^{-n}\}})\longrightarrow0.
	\end{align}
 Hence, 
 \[
 \mathbb P\left(|\mathbb EY-\mathbb ET_n|>\frac\varepsilon2\right)=0\hbox{ for all large enough }n. 
 \]
 
Combining \eqref{eqnA2} with the above estimates yields \eqref{eq:A.1}.
By the Borel--Cantelli lemma, \eqref{eq:A.1} implies that 
	\[
	\sum_{w\in I^n}a_wY_w\longrightarrow\mathbb EY
	\qquad\text{almost surely}.
	\]
	Using the basic relation $|x-y|=x+y-2\min\{x,y\}$ for $x,y\ge0$, we have
	\[
	\mathbb E\left|\sum_{w\in I^n}a_wY_w-\mathbb EY\right|
	=2\mathbb EY-2\mathbb E\min\left\{
	\sum_{w\in I^n}a_wY_w,\mathbb EY\right\}
	\longrightarrow0,
	\]
	where we use the dominated convergence theorem.
\end{proof}

\section{Proof of the Peirone-type theorem}
\label{sec:peirone-proof}
In this appendix, we prove the Peirone-type convergence theorem (Theorem \ref{thm:peirone}). Throughout this appendix, we will use the following notation. 

\begin{definition}
\begin{enumerate}
\item For $1\le i\le M$, define 
\begin{equation}
	\R_{i,+}^{V_0}:=\{u\in\mathbb R^{V_0}:u(q_i)=0\hbox{ and } u(q_j)\ge0\hbox{ for }j\ne i\},
\end{equation}
and 
\begin{equation}
\Delta_i=\Big\{u\in\mathbb R_{i,+}^{V_0}:\sum_{j\ne i}u(q_j)=1\Big\}.
\end{equation}

\item For $u\in\R^{V_0}$ and $E\in\mathcal{Q}_p(V_0)$, we let $H_Eu$ be the harmonic extension with respect to $f\mapsto\sum_{i\in I}E(f\circ\psi_i)$. For $1\le i\le M$, define $H^E_i:\R^{V_0}\to\R^{V_0}$ by 
\[
H_i^E u:=(H_Eu)\circ\psi_i|_{V_0}. 
\]
This map preserves $\R_{i,+}^{V_0}$.
\end{enumerate}
\end{definition}

In preparation for the proof, we first study some properties of $H_i^E$.

\begin{lemma}\label{lem:strict-group-current}
Let $V$ be a nonempty finite set. Let $E\in\cQ_p(V)$ and $h_1\le h_2$. Suppose that 
\[
  Z:=\{x\in V:h_1(x)=h_2(x)\}
\]
is nonempty. Then
\begin{equation} \label{eq:current-weak}
 J_E(h_1,x)\ge J_E(h_2,x)\quad \hbox{ for every }x\in Z.
\end{equation}
If further $Z\neq V$, then 
\begin{equation}\label{eq:current-strict}
J_E(h_1,x)>J_E(h_2,x)\quad\hbox{ for some }x\in Z.
\end{equation}
\end{lemma}

\begin{proof}
Fix $x\in Z$. Then for any $\varepsilon>0$, 
\[
h_2\wedge (h_1+\varepsilon\mathbf1_x)=h_1\ \hbox{ and }\ h_2\vee (h_1+\varepsilon\mathbf1_x)=h_2+\varepsilon\mathbf1_x. 
\] 
As the map $(s,t)\mapsto(s\wedge t,s\vee t)$ is non-expansive from $\ell^p$ to $\ell^p$ in the sense of \eqref{nonexpansive}, 
\[
E(h_1)+E(h_2+\varepsilon\mathbf1_x)\le E(h_2)+E(h_1+\varepsilon\mathbf1_x),
\] 
by the generalized contraction property (RF5). 
This implies 
\[
\begin{aligned}
J_E(h_2,x)
&=\lim_{\varepsilon\downarrow0}\frac{E(h_2+\varepsilon\mathbf1_x)-E(h_2)}{p\varepsilon}\\
&\le\lim_{\varepsilon\downarrow0}\frac{E(h_1+\varepsilon\mathbf1_x)-E(h_1)}{p\varepsilon}
=J_E(h_1,x),
\end{aligned}
\]
which is \eqref{eq:current-weak}. 

Next, we prove \eqref{eq:current-strict} by contradiction. Assume $Z\neq V$ and 
\begin{equation}\label{e:lemmaB1.3}
 J_E(h_1,x)=J_E(h_2,x)\quad \hbox{ for every }x\in Z.
\end{equation}
We note that 
\begin{equation}\label{e:lemmaB1.4}
	h_1\le h_1+t\mathbf1_{V\setminus Z}\le h_2\ \hbox{ for }0\le t\le t_0\hbox{ with equality on }Z,
\end{equation}
where $t_0:=\min_{x\in V\setminus Z}(h_2(x)-h_1(x))>0$. 
Then, \eqref{eq:current-weak}, \eqref{e:lemmaB1.3}  and \eqref{e:lemmaB1.4} together imply 
\[
 J_E(h_1+t\mathbf1_{V\setminus Z},x)=J_E(h_1,x)
\ \hbox{ for every }x\in Z.
\]
The current identities \eqref{eq:current-euler} now imply
\[
 \frac{d}{dt}E(h_1+t\mathbf1_{V\setminus Z})
 =p\sum_{x\in V\setminus Z}
   J_E(h_1+t\mathbf1_{V\setminus Z},x)
 =-p\sum_{x\in Z}J_E(h_1,x).
\]
Thus the derivative is constant on an interval of positive length, contrary
to strict convexity of $E$ modulo constants, since $\mathbf1_{V\setminus Z}$ is nonconstant (see also \cite[Theorem~3.7]{KajinoShimizu2026}).
\end{proof}

\begin{lemma}\label{lem:cell-propagation}
Let $u,v\in\R^{V_0}$ and $E\in\mathcal{Q}_p(V_0)$. Let $h_u,h_v$ be the harmonic extensions of $u,v$ with respect to $f\mapsto\sum_{w\in I^2}E(f\circ\psi_w)$, respectively. Suppose $u\geq v$. Then (a) and (b) hold.
\begin{enumerate}
\item[(a)] For each $w\in I^2$ with $\psi_w(V_0)\cap V_0=\emptyset$, either $h_u>h_v$ or $h_u=h_v$ on $\psi_w(V_0)$.

\item[(b)] For every $1\leq i\leq M$ such that $u(q_i)>v(q_i)$, $h_u|_{\psi_{ii}(V_0)}>h_v|_{\psi_{ii}(V_0)}$. 
\end{enumerate}
\end{lemma}
\begin{proof}
Let $x\in V_2\setminus V_0$ and suppose that $h_u(x)=h_v(x)$. For every $\theta\in I^2$ such that $x\in\psi_\theta(V_0)$, \eqref{eq:current-weak} and the fact $h_u\geq h_v$ give
\[
J_E(h_u\circ\psi_\theta,\psi_\theta^{-1}x)\leq J_E(h_v\circ\psi_\theta,\psi_\theta^{-1}x).
\]
On the other hand, the Kirchhoff equations for both harmonic extensions imply  
\[
\begin{aligned}
\sum_{\theta\in I^2:\,x\in \psi_\theta(V_0)}J_E(h_u\circ\psi_\theta,\psi_\theta^{-1}x)&=0,\\
\sum_{\theta\in I^2:\,x\in \psi_\theta(V_0)}J_E(h_v\circ\psi_\theta,\psi_\theta^{-1}x)&=0.
\end{aligned}
\]
Therefore, for every $x\in V_2\setminus V_0$ with $h_u(x)=h_v(x)$ and every $\theta\in I^2$ such that $x\in\psi_\theta(V_0)$, 
\begin{equation}\label{lemmaB2.1}
J_E(h_u\circ\psi_\theta,\psi_\theta^{-1}x)=J_E(h_v\circ\psi_\theta,\psi_\theta^{-1}x). 
\end{equation}

(a). Now, we fix $w\in I^2$ with $\psi_w(V_0)\cap V_0=\emptyset$, and let $Z=\{x\in\psi_w(V_0):\,h_u(x)=h_v(x)\}$. To prove (a), it suffices to show that either $Z=\emptyset$ or $Z=\psi_w(V_0)$. We prove this by contradiction. Suppose $Z\neq\emptyset$ and $Z\neq\psi_w(V_0)$. Then $J_E(h_u\circ\psi_w,\psi_w^{-1}x)<J_E(h_v\circ\psi_w,\psi_w^{-1}x)$ for some $x\in Z$ by \eqref{eq:current-strict}, which contradicts \eqref{lemmaB2.1}. 

(b). We prove (b) by contradiction. Suppose $Z=\{x\in\psi_{ii}(V_0):\,h_u(x)=h_v(x)\}\neq\emptyset$. Since $h_u(q_i)>h_v(q_i)$, $Z\neq\psi_{ii}(V_0)$. Then, by \eqref{eq:current-strict}, there is $x\in Z$ so that 
\[
J_E(h_u\circ\psi_{ii},\psi_{ii}^{-1}x)<J_E(h_v\circ\psi_{ii},\psi_{ii}^{-1}x). 
\]
Since $\psi_{ii}(V_0)\cap V_0=\{q_i\}$ by \cite[Lemma~3.1.4]{Kohl2022}, we have $Z\subset V_2\setminus V_0$, contradicting \eqref{lemmaB2.1}.
\end{proof}

\begin{lemma}\label{lem:transitive-graph}
For distinct $w,v\in I^n$, $n\geq 1$, we say $w\sim v$ if $\psi_w(V_0)\cap\psi_v(V_0)\neq\emptyset$. 
 	
(a). For each $1\le i\le M$, there is $\mathcal{I}_i\subset I\setminus\{i\}$ so that $\{1,2,\cdots,M\}\setminus\{i\}\subset\mathcal{I}_i$ and both the graphs $(\mathcal{I}_i,\sim)$ and $(\mathcal{I}_i\cup\{i\},\sim)$ are connected.

(b). Let 
\begin{equation}
	\mathcal I=
	\{jj':j\in I\setminus\{1,\ldots,M\},\ j'\in I\}
	\cup
	\{jj':1\le j\le M,\ j'\in \mathcal{I}_j\}.
\end{equation}
Then, $(\mathcal{I},\sim)$ is connected. Moreover, for each $i\in\{1,2,\cdots,M\}$, there is $w\in \mathcal{I}$ so that $w\sim ii$.
\end{lemma}
\begin{proof}
(a). Choose a connected subgraph $\mathcal{G}'$ containing $\{1,2,\cdots,M\}$ with as few edges as possible. Then, $\mathcal{G}'$ is a tree, and every vertex of degree one in $\mathcal{G}'$ belongs to $\{1,2,\cdots,M\}$. By symmetry, we may assume that $i$ has degree $1$ in $\mathcal{G}'$. It suffices to take $\mathcal{I}_i$ to be the set of vertices of $\mathcal{G}'$ after deleting $i$.

(b). It is an easy consequence of (a) and basic properties of nested fractals.
\end{proof}

\begin{proposition}\label{prop:strong-order}
Let $E\in\cQ_p(V_0)$, $1\le i\le M$, and
$u,v\in\mathbb R_{i,+}^{V_0}$. Then 
\[
u\leq v\hbox{ and }u\neq v\ \Longrightarrow\ H_i^{\cR(E)}u(q_j)<H_i^{\cR(E)}v(q_j)\hbox{ for }j\neq i.
\]
In particular, for every nonzero $u\in\mathbb R_{i,+}^{V_0}$, $H_i^{\cR(E)}u(q_j)>0$ for all $j\ne i$.
\end{proposition}
\begin{proof}
Choose $j_0\ne i$ such that $u(q_{j_0})<v(q_{j_0})$. Let $h_u,h_v$ be the harmonic extensions of $u,v$ with respect to $f\mapsto\sum_{w\in I^2}E(f\circ\psi_w)$, respectively. Then, $h_u|_{\psi_{j_0j_0}(V_0)}<h_v|_{\psi_{j_0j_0}(V_0)}$ by Lemma~\ref{lem:cell-propagation}. By Lemma \ref{lem:transitive-graph}(b), there are $w\in \mathcal{I}$ and $x\in \psi_w(V_0)$ such that $h_u(x)<h_v(x)$. 

By \cite[Lemma~3.1.4]{Kohl2022}, it holds that $V_0\cap \bigcup_{w\in\mathcal{I}}\psi_w(V_0)=\emptyset$. Hence, by the connectedness of $\mathcal{I}$ and Lemma \ref{lem:cell-propagation}, $h_u<h_v$ on $\bigcup_{w\in\mathcal{I}}\psi_w(V_0)$. Finally, by Lemma \ref{lem:transitive-graph}(a), $ij\in \mathcal{I}$ for each $j\in\{1,2,\cdots,M\}\setminus\{i\}$. Hence 
\[
\begin{aligned}
H_i^{\cR(E)}u(q_j)
&=h_u(\psi_iq_j)=h_u(\psi_i\psi_jq_j)
<h_v(\psi_i\psi_jq_j)=h_v(\psi_iq_j)
=H_i^{\cR(E)}v(q_j).
\end{aligned}
\]  
\end{proof}

\begin{lemma}\label{lem:harmonic-continuity}
The map
\[
 (E,u)\longmapsto H_Eu
\]
from $[aE_0,bE_0]\times\mathbb R^{V_0}$ to $\mathbb R^{V_1}$ is jointly
continuous for every $0<a\le b<\infty$.
\end{lemma}
\begin{proof}
Let $E_n\to E$ in form norm and $u_n\to u$. By the maximum principle,
the functions $H_{E_n}u_n$ are uniformly bounded.

It suffices to prove that any subsequential limit of $H_{E_n}u_n$ is $H_Eu$. Let $n_j$ be a subsequence so that $H_{E_{n_j}}u_{n_j}$ converges to some limit function $h$ as $j\to\infty$. Then 
\[
\begin{aligned}
\sum_{i\in I}E(h\circ\psi_i|_{V_0})
&=\lim_{j\to\infty}\sum_{i\in I}E_{n_j}(H_{E_{n_j}}u_{n_j}\circ\psi_i|_{V_0})\\
&\le\lim_{j\to\infty}\sum_{i\in I}E_{n_j}(g_{n_j}\circ\psi_i|_{V_0})\\
&=\sum_{i\in I}E(H_Eu\circ\psi_i|_{V_0}),
\end{aligned}
\]
where for $n\geq 1$
\[
g_n(x)=
\begin{cases}
	u_n(x)&\hbox{ if }x\in V_0,\\
	H_Eu(x)&\hbox{ if }x\in V_1\setminus V_0. 
\end{cases}
\]
Thus $h=H_Eu$ by the uniqueness of harmonic extension.  Every convergent subsequence has the same limit, which proves joint continuity.
\end{proof}

\begin{corollary}\label{coroB7}
For $1\le i\le M$ and $0<a\le b<\infty$, $\{H_i^E:E\in[aE_0,bE_0]\}$ is a compact family of continuous maps from $\R^{V_0}_{i,+}$ to $\R^{V_0}_{i,+}$ with respect to uniform convergence on $\Delta_i$. 
\end{corollary}

Combining Proposition \ref{prop:strong-order} and Corollary \ref{coroB7}, we can prove the following convergence theorem. 

\begin{proposition}\label{prop:B8}
Let $G\in\mathcal{Q}_p(V_0)$ and let $G_k=\cR^k(G)$ for $k=0,1,2,\cdots$. Fix $i\in\{1,2,\cdots,M\}$, and define $\mathfrak n:\R^{V_0}_{i,+}\setminus\{0\}\to\Delta_i$ by 
\[
\mathfrak n(u)=\frac{u}{\sum_{j\ne i}u(q_j)}.
\]
Then, there is $v\in\Delta_i$ with $v(q_j)>0$ for $j\ne i$ such that 
\[
\sup_{u\in\R^{V_0}_{i,+}\setminus\{0\}}
\big\|\mathfrak n\big(H^{G_1}_i\cdots H^{G_n}_iu\big)-v\big\|_{\ell^1}\longrightarrow0. 
\]
\end{proposition}

We state the proof in a more general setting.
 
\begin{lemma}\label{lem:backward-products}
Let $\Delta=\{x\in\mathbb R_+^{M-1}:\sum_jx_j=1\}$ and
$\mathfrak n(x)=x/(\sum_jx_j)$ for $x\ne0$.  Let $\mathscr A$ be a compact family,
in the norm of uniform convergence on $\Delta$, of continuous,
degree-one homogeneous maps
$\Theta:\mathbb R_+^{M-1}\to\mathbb R_+^{M-1}$ satisfying
\begin{equation}\label{eq:strong-map}
 0\le x\le y,\quad x\ne y
 \quad\Longrightarrow\quad \Theta x<\Theta y.
\end{equation}
For every sequence $\Theta_n\in\mathscr A,n\geq 0$, there is $v\in \operatorname{int}(\Delta):=\{x\in\Delta:x_j>0\hbox{ for each }j\}$ so that
\begin{equation}\label{eq:backward-synchronization}
 \sup_{x\in\Delta}
 \bigl\|\mathfrak n(\Theta_0\Theta_1\cdots\Theta_{n-1}x)-v\bigr\|_{\ell^1}
 \longrightarrow0
 \qquad\hbox{ as }n\to\infty.
\end{equation}
\end{lemma}
\begin{proof} 
If $M=2$, then $\Delta$ is a
singleton and the assertion is trivial.
Assume $M\ge3$. 

For strictly positive $x,y\in\R^{M-1}_+$, define the Hilbert projective metric by 
\[
 d_{\mathrm{HP}}(x,y)=\log\max_j\frac{x_j}{y_j}+\log\max_j\frac{y_j}{x_j}.
\]
It is easy to check that $d_{\mathrm{HP}}(x,y)=d_{\mathrm{HP}}(\mathfrak n(x),\mathfrak n(y))$ for every strictly positive $x,y\in\R^{M-1}_+$, and $d_{\mathrm{HP}}$ is a metric on $\operatorname{int}(\Delta)$.  

We claim that for every $\Theta\in\mathscr{A}$ and strictly positive $x,y$,
\begin{equation}\label{eq:hilbert-nonexpansive}
d_{\mathrm{HP}}(\Theta x,\Theta y)\le d_{\mathrm{HP}}(x,y)\qquad \hbox{ with equality if and only if $d_{\mathrm{HP}}(x,y)=0$.}
\end{equation}
For every strictly positive $x,y\in\R^{M-1}_+$, let $s=\max_j\frac{x_j}{y_j}$ and $t=\max_j\frac{y_j}{x_j}$. Then $d_{\mathrm{HP}}(x,y)=\log(st)$. By homogeneity and \eqref{eq:strong-map}, $\Theta x\leq s\Theta y$ and $\Theta y\leq t\Theta x$. Moreover, if $st>1$, it is impossible that both $x=sy$ and $y=tx$ hold. Hence either $\Theta x<s\Theta y$ or $\Theta y<t\Theta x$ by \eqref{eq:strong-map}. The claim follows. 

For simplicity, we write $\widetilde\Theta=\mathfrak n\circ\Theta$ for $\Theta\in\mathscr A$. By \eqref{eq:strong-map} and the compactness of $\mathscr{A}\times\Delta$, 
\[
\Gamma:=\big\{\widetilde\Theta x:\,\Theta\in \mathscr{A}\hbox{ and }x\in\Delta\big\}
\]
is a compact subset of $\hbox{int}(\Delta)$.  By compactness and \eqref{eq:hilbert-nonexpansive}, for every $\varepsilon>0$, there is
$\delta>0$, independent of $x,y\in\Gamma$ and $\Theta\in\mathscr{A}$, such that
\begin{equation}\label{eq:uniform-hilbert-drop}
d_{\mathrm{HP}}(x,y)\ge \varepsilon
 \quad\Longrightarrow\quad
 d_{\mathrm{HP}}(\Theta x,\Theta y)
 \le d_{\mathrm{HP}}(x,y)-\delta.
\end{equation}

Now, for a sequence $\Theta_0,\Theta_1,\cdots\in\mathscr{A}$, let 
\[
\Gamma_0=\Delta,\qquad
\Gamma_n=\widetilde\Theta_0\widetilde\Theta_1\cdots\widetilde\Theta_{n-1}(\Delta)\subset\Gamma
\quad(n\ge1). 
\]
Then, 
\[
\Gamma_1\supset\Gamma_2\supset\Gamma_3\supset\cdots
\]
and their $d_{\mathrm{HP}}$-diameters tend to zero. Their intersection is therefore a singleton $\{v\}\subset\Gamma$, which proves \eqref{eq:backward-synchronization}.
\end{proof}\bigskip 

\begin{proof}[Proof of Proposition~\ref{prop:B8}]
Choose $0<a\le b<\infty$ so that $aE_0\le G\le bE_0$.
By Lemma~\ref{lem:compact-annulus}, $G_k\in[aE_0,bE_0]$ for all $k\ge0$.
Proposition~\ref{prop:strong-order} and Corollary~\ref{coroB7} imply that
\[
\mathscr A=\{H_i^{\cR(E)}:E\in[aE_0,bE_0]\}
\]
is a compact family satisfying \eqref{eq:strong-map}.
Apply Lemma~\ref{lem:backward-products} with $\Delta=\Delta_i$ and $\Theta_n=H_i^{G_{n+1}}$.
\end{proof}

We finalize the proof of Theorem \ref{thm:peirone}. 

\begin{lemma}\label{lem:trace-contact}
Let $E,G\in \mathcal{Q}_p(V_0)$ and $u\in\R^{V_0}$. 
Assume $G\le E$ and
$\cR(G)(u)=\cR(E)(u)$.  Then
$H_Eu=H_Gu=:h$, and
\begin{equation}\label{B.10}
 G(h\circ\psi_i|_{V_0})=E(h\circ\psi_i|_{V_0})\quad\hbox{ for }i\in I. 
\end{equation}
\end{lemma}
\begin{proof}
Put $h_E=H_Eu$ and $h_G=H_Gu$.  The identity
\begin{align*}
0&=\cR(E)(u)-\cR(G)(u)\\
&=\rho\sum_{i\in I}
\Bigl[E(h_E\circ\psi_i|_{V_0})
-G(h_G\circ\psi_i|_{V_0})\Bigr]\\
&=\rho\sum_{i\in I}
   \Bigl[E(h_E\circ\psi_i|_{V_0})
         -G(h_E\circ\psi_i|_{V_0})\Bigr]+
   \rho\left[
   \sum_{i\in I}G(h_E\circ\psi_i|_{V_0})
   -\sum_{i\in I}G(h_G\circ\psi_i|_{V_0})
   \right]
\end{align*}
expresses zero as a sum of nonnegative terms. In particular, $\sum_{i\in I}G(h_E\circ\psi_i|_{V_0})
=\sum_{i\in I}G(h_G\circ\psi_i|_{V_0})$, so $h_E=h_G$ by the uniqueness of harmonic extension. \eqref{B.10} follows from the fact that each term in
the first sum also vanishes.
\end{proof}

For forms $E$ and $G$, define
\begin{equation}\label{eq:comparison-ratio}
	\Lambda(G/E):=\max_{u\in\cS}\frac{G(u)}{E(u)}.
\end{equation}
The ratio is continuous in $(E,G)$ on a common interval $[aE_0,bE_0]$. 

\begin{proposition}\label{prop:no-persistent-contact}
If $G_\ell=\cR^\ell(G_0)$ with $G_0\in\mathcal{Q}_p(V_0)$, no constant $\theta>1$
can satisfy either of the following identities:
\begin{align}
 \Lambda(G_{\ell+1}/G_\ell)&=\theta
 \qquad\text{for every }\ell\geq0,\label{eq:persistent-forward}\\
 \Lambda(G_\ell/G_{\ell+1})&=\theta
 \qquad\text{for every }\ell\geq0.\label{eq:persistent-backward}
\end{align}
\end{proposition}

\begin{proof}
We prove the proposition by contradiction. In the following, we assume that
\eqref{eq:persistent-backward} holds. The other case is treated in the same way. \smallskip  

As an immediate consequence of Lemma \ref{lem:trace-contact} and \eqref{eq:persistent-backward}, the following claim holds. 

\noindent\textbf{Claim 1.} Let $\ell\geq 0$, $i\in\{1,2,\cdots,M\}$ and $u_{\ell+1}\in\R^{V_0}$. If $G_{\ell+1}(u_{\ell+1})=\theta G_{\ell+2}(u_{\ell+1})$, then 
\[
H_i^{G_\ell}u_{\ell+1}=H_i^{G_{\ell+1}}u_{\ell+1}=:u_\ell
\]
and
\[
G_\ell(u_\ell)=\theta G_{\ell+1}(u_\ell).
\] 
\begin{proof}[Proof of Claim 1]
By \eqref{eq:persistent-backward} and the assumption,
\[
G_\ell\le\theta G_{\ell+1},\qquad
\cR(G_\ell)(u_{\ell+1})=\cR(\theta G_{\ell+1})(u_{\ell+1}).
\]
Applying Lemma~\ref{lem:trace-contact} to $G_\ell$ and $\theta G_{\ell+1}$, we deduce the claim.
\end{proof} 

For each $n\geq 1$, we choose $u^{(n)}_n\in\R^{V_0}\setminus \R{\bf 1}$ so that $G_n(u^{(n)}_n)=\theta G_{n+1}(u^{(n)}_n)$. By applying Claim 1 repeatedly, we can see the following claim. 

\noindent\textbf{Claim 2.} Let $n\geq 1$ and $i\in\{1,2,\cdots, M\}$. Then for every $0\leq \ell\leq n-1$, 
\[
H_i^{G_{\ell}}H_i^{G_{\ell+1}}\cdots H_i^{G_{n-1}}u^{(n)}_n=H_i^{G_{\ell+1}}H_i^{G_{\ell+2}}\cdots H_i^{G_n}u^{(n)}_n=:u^{(n)}_\ell
\] 
and 
\[
G_\ell(u^{(n)}_\ell)=\theta G_{\ell+1}(u^{(n)}_\ell).\medskip 
\]

By translating and renormalizing $u^{(n)}_n$ and passing to a subsequence, we may assume that there are $i\in\{1,2,\cdots,M\}$ and $n_j$ so that $u^{(n_j)}_{n_j}\in\Delta_i$ for every $j\geq 1$. For $\ell\geq 1$, by Proposition \ref{prop:B8}, there is $v_\ell\in \Delta_i$ so that 
\[
\sup_{u\in\R^{V_0}_{i,+}\setminus\{0\}}\|\mathfrak n\big(H^{G_\ell}_i\cdots H^{G_n}_iu\big)-v_\ell\|_{\ell^1}\longrightarrow0. 
\]
Then, for any $\ell\geq 1$, Claim 2 implies 
\[
\lim_{j\to\infty}\mathfrak n(u_\ell^{(n_j)})=v_\ell=v_{\ell+1}.
\]
Hence, $v_1=v_2=\cdots$ and we denote their common value by $v$. Again, by Claim 2 and homogeneity,
\[
G_\ell(\mathfrak n(u_\ell^{(n_j)}))
=\theta G_{\ell+1}(\mathfrak n(u_\ell^{(n_j)})).
\]
Taking $j\to\infty$ gives $G_\ell(v)=\theta G_{\ell+1}(v)$ for every $\ell\ge1$.
On the other hand, choose $a>0$ so that $G_0\ge aE_0$. Then $G_\ell\ge aE_0$ for all $\ell\ge0$, so
\[
0<aE_0(v)\le G_\ell(v)=\theta^{1-\ell}G_1(v)\longrightarrow0,
\]
a contradiction.
\end{proof}

\begin{proof}[Proof of Theorem~\ref{thm:peirone}]
Fix $E\in\cQ_p(V_0)$. Choose $0<a\le b<\infty$
so that $aE_0\le E\le bE_0$. Then $aE_0\le \cR^n(E)\le bE_0$ for any $n\geq 0$.

Define
\begin{equation}\label{eq:successive-ratios}
 \Lambda_n^+=\Lambda\big(\cR^{n+1}(E)/\cR^n(E)\big),
 \qquad
 \Lambda_n^-=\Lambda\big(\cR^n(E)/\cR^{n+1}(E)\big).
\end{equation}
By the order preservation and homogeneity properties of $\cR$, both $\Lambda_n^+$ and $\Lambda_n^-$ are nonincreasing in $n\geq 0$. They are at least one.  Indeed, if $\Lambda_n^+=c<1$ for some $n\geq 0$, then $\cR^{n+k}(E)\le c^k\cR^n(E)$ for any $k\geq 0$, contradicting the common lower bound; if
$\Lambda_n^-=c<1$ for some $n\geq 0$, then $\cR^{n+k}(E)\ge c^{-k}\cR^n(E)$ for any $k\geq 0$, contradicting the common
upper bound.  Thus $\Lambda_n^+\downarrow\Lambda^+\ge1$ and
$\Lambda_n^-\downarrow\Lambda^-\ge1$.

We choose a convergent subsequence $\cR^{n_k}(E),k\geq 1$ and let $G=\lim\limits_{k\to\infty}\cR^{n_k}(E)$. Then, for each $\ell\geq 0$,  $\cR^{n_k+\ell}E\to G_\ell:=\cR^\ell(G)$. By the continuity of \eqref{eq:comparison-ratio}, 
\begin{equation}\label{e:G}
 \Lambda(G_{\ell+1}/G_\ell)=\Lambda^+
 \ \hbox{ and }\ 
 \Lambda(G_\ell/G_{\ell+1})=\Lambda^-
\quad\hbox{ for }\ell\geq 0.
\end{equation}
By Proposition~\ref{prop:no-persistent-contact}, 
\begin{equation}\label{eq:ratios-to-one}
\Lambda^+=1\ \hbox{ and }\
\Lambda^-=1.
\end{equation}
Hence, \eqref{e:G} implies $\cR(G)=G$. 

Finally, we notice that $\cR$ is non-expansive with respect to the Thompson metric
\begin{equation}\label{eq:thompson-metric}
d_T(D,D'):=\log\max\{\Lambda(D/D'),\Lambda(D'/D)\}\quad\hbox{ for }D,D'\in\mathcal{Q}_p(V_0). 
\end{equation}
Hence, for any $n\geq 0$ and $n_k\leq n$,
\[
 d_T(\cR^n(E),G)
 =d_T(\cR^{n-n_k}\cR^{n_k}(E),\cR^{n-n_k}(G))
 \le d_T(\cR^{n_k}(E),G).
\]
As the right-hand side tends to $0$ by the convergence $\cR^{n_k}(E)\to G$, $\cR^n(E)$ converges to $G$ as $n\to\infty$.
\end{proof}

\begingroup
\renewcommand{\addcontentsline}[3]{}
\section*{\texorpdfstring{{Declaration of use of AI tools}}{Declaration of use of AI tools}}
{The authors used GPT-6 Astra to help organize and refine their proof of the $p$-version of Peirone's convergence theorem in Appendix~\ref{sec:peirone-proof}. They also used the tool to polish the English and proofread the manuscript. The authors take full responsibility for the manuscript.}\par
\endgroup

\enlargethispage{3\baselineskip}


\begin{thebibliography}{99}
\bibitem{AndresBarlowDeuschelHambly2013}
S.~Andres, M.~T. Barlow, J.-D. Deuschel and B.~M. Hambly,
\emph{Invariance principle for the random conductance model},
Probab. Theory Related Fields \textbf{156} (2013), no.~3--4, 535--580.

\bibitem{AndresDeuschelSlowik2015}
S.~Andres, J.-D. Deuschel and M.~Slowik,
\emph{Invariance principle for the random conductance model in a degenerate ergodic environment},
Ann. Probab. \textbf{43} (2015), no.~4, 1866--1891.

\bibitem{BandtKeller1991}
C.~Bandt and K.~Keller,
\emph{Self-similar sets 2. A simple approach to the topological structure of fractals},
Math. Nachr. \textbf{154} (1991), 27--39.

\bibitem{Barlow1998}
M.~T. Barlow,
\emph{Diffusions on fractals},
in: Lectures on Probability Theory and Statistics (Saint-Flour, 1995),
Lecture Notes in Math. 1690, Springer, Berlin, 1998, pp.~1--121.

\bibitem{BarlowBass1989}
M.~T. Barlow and R.~F. Bass,
\emph{The construction of Brownian motion on the Sierpi\'nski carpet},
Ann. Inst. H. Poincar\'e Probab. Statist. \textbf{25} (1989), no.~3, 225--257.

\bibitem{BarlowBass1999}
{M.~T. Barlow and R.~F. Bass,
\emph{Brownian motion and harmonic analysis on Sierpi\'nski carpets},
Canad. J. Math. \textbf{51} (1999), no.~4, 673--744.}

\bibitem{BBK}
M.~T. Barlow, R.~F. Bass and T.~Kumagai, 
\emph{Stability of parabolic Harnack inequalities on metric measure spaces},
J. Math. Soc. Japan 58 (2006), no. 2, 485--519.

\bibitem{BarlowDeuschel2010}
M.~T. Barlow and J.-D. Deuschel,
\emph{Invariance principle for the random conductance model with unbounded conductances},
Ann. Probab. \textbf{38} (2010), no.~1, 234--276.

\bibitem{BGK} 
M.~T. Barlow, A.~Grigor'yan and T.~Kumagai, \emph{On the equivalence of parabolic Harnack inequalities and heat kernel estimates}, J. Math. Soc. Japan 64 (2012), no. 4, 1091--1146.

\bibitem{BarlowPerkins1988}
M.~T. Barlow and E.~A. Perkins,
\emph{Brownian motion on the Sierpi\'nski gasket},
Probab. Theory Related Fields \textbf{79} (1988), no.~4, 543--623.

\bibitem{BS}
P. Bella and M. Schäffner, \emph{Quenched invariance principle for random walks among random degenerate conductances.} Ann. Probab. 48 (2020), no. 1, 296--316.

\bibitem{BergerBiskup2007}
N.~Berger and M.~Biskup,
\emph{Quenched invariance principle for simple random walk on percolation clusters},
Probab. Theory Related Fields \textbf{137} (2007), no.~1--2, 83--120.

\bibitem{Bogachev2007}
V.~I. Bogachev,
\emph{Measure Theory, Vol. I},
Springer-Verlag, Berlin, 2007.

\bibitem{BoucheronLugosiMassart2013}
S.~Boucheron, G.~Lugosi and P.~Massart,
\emph{Concentration Inequalities: A Nonasymptotic Theory of Independence},
Oxford University Press, Oxford, 2013.

\bibitem{CaoGuQiu2022}
S.~Cao, Q.~Gu and H.~Qiu,
\emph{$p$-energies on p.c.f. self-similar sets},
Adv. Math. \textbf{405} (2022), Paper No.~108517.

\bibitem{CF}
Z.-Q. Chen and M. Fukushima, \emph{Symmetric Markov Processes, Time Change, and Boundary Theory}, London Mathematical Society Monographs
Series, vol.~35, Princeton University Press, Princeton, NJ, 2012.


\bibitem{CroydonHamblyKumagai2017}
D.~A. Croydon, B.~M. Hambly and T.~Kumagai,
\emph{Time-changes of stochastic processes associated with resistance forms},
Electron. J. Probab. \textbf{22} (2017), Paper No.~82, 41 pp.

\bibitem{EfronStein1981}
B.~Efron and C.~Stein,
\emph{The jackknife estimate of variance},
Ann. Statist. \textbf{9} (1981), no.~3, 586--596.

\bibitem{EthierKurtz1986}
S.~N. Ethier and T.~G. Kurtz,
\emph{Markov Processes: Characterization and Convergence},
Wiley, New York, 1986.

\bibitem{FHK}
	P. J. Fitzsimmons, B. M. Hambly and T. Kumagai,
	\emph{Transition density estimates for Brownian motion on affine nested fractals},
	Comm. Math. Phys. \textbf{165} (1994), no. 3, 595--620.

\bibitem{FukushimaOshimaTakeda2011}
M.~Fukushima, Y.~Oshima and M.~Takeda,
\emph{Dirichlet Forms and Symmetric Markov Processes},
2nd revised and extended ed.,
De Gruyter Studies in Mathematics 19,
De Gruyter, Berlin, 2011.

\bibitem{Goldstein1987}
S.~Goldstein,
\emph{Random walks and diffusions on fractals},
in: Percolation Theory and Ergodic Theory of Infinite Particle Systems,
IMA Vol. Math. Appl. 8, Springer, New York, 1987, pp.~121--129.

\bibitem{HamblyMetz1998}
B.~M. Hambly and V.~Metz,
\emph{The homogenization problem for the Vicsek set},
Stochastic Process. Appl. \textbf{76} (1998), no.~2, 167--190.

\bibitem{HermanPeironeStrichartz2004}
P.~E. Herman, R.~Peirone and R.~S. Strichartz,
\emph{$p$-energy and $p$-harmonic functions on Sierpi\'nski gasket type fractals},
Potential Anal. \textbf{20} (2004), 125--148.

\bibitem{Hutchinson1981}
J.~E. Hutchinson,
\emph{Fractals and self similarity},
Indiana Univ. Math. J. \textbf{30} (1981), no.~5, 713--747.

\bibitem{KajinoShimizuSurvey}
N.~Kajino and R.~Shimizu,
\emph{$p$-energy forms on fractals: recent progress},
arXiv:2501.09002v1 [math.FA], 15 January 2025.

\bibitem{KajinoShimizu2026}
N.~Kajino and R.~Shimizu,
\emph{Contraction properties and differentiability of $p$-energy forms with applications to nonlinear potential theory on self-similar sets},
arXiv:2404.13668v5 [math.FA], 10 September 2026.

\bibitem{Kesten1986}
H.~Kesten,
\emph{Subdiffusive behavior of random walk on a random cluster},
Ann. Inst. H. Poincar\'e Probab. Statist. \textbf{22} (1986), no.~4, 425--487.

\bibitem{Kigami1993}
J.~Kigami, \emph{Harmonic calculus on p.c.f. self-similar sets,}
Trans. Amer. Math. Soc. {\bf 335} (1993), no. 2, 721--755.

\bibitem{Kigami2001}
J.~Kigami,
\emph{Analysis on Fractals},
Cambridge Tracts in Mathematics 143, Cambridge University Press, Cambridge, 2001.

\bibitem{Kigami2020}
J.~Kigami, \emph{Geometry and analysis of metric spaces via weighted partitions.} Lecture Notes in Math., 2265, Springer, Cham, 2020, viii+162 pp.


\bibitem{Kigami2023}
J.~Kigami, \emph{Conductive homogeneity of compact metric spaces and construction of p-energy},  Mem. Eur. Math. Soc., 5
European Mathematical Society (EMS), Berlin, 2023, viii+129 pp.


\bibitem{Kohl2022}
S.~Kohl,
\emph{Inhomogeneous Fractals as a Martin Boundary},
doctoral dissertation, University of Stuttgart, 2022.

\bibitem{Kozlov1985}
S.~M. Kozlov,
\emph{The method of averaging and walks in inhomogeneous environments},
Russian Math. Surveys \textbf{40} (1985), no.~2, 73--145.

\bibitem{Kozlov1993}
S.~M. Kozlov,
\emph{Harmonization and homogenization on fractals},
Comm. Math. Phys. \textbf{153} (1993), no.~2, 339--357.

\bibitem{KumagaiKusuoka1996}
T.~Kumagai and S.~Kusuoka,
\emph{Homogenization on nested fractals},
Probab. Theory Related Fields \textbf{104} (1996), no.~3, 375--398.

\bibitem{Kunnemann1983}
R.~K\"unnemann,
\emph{The diffusion limit for reversible jump processes on $\mathbb Z^d$ with ergodic random bond conductivities},
Comm. Math. Phys. \textbf{90} (1983), no.~1, 27--68.

\bibitem{Kusuoka1987}
S.~Kusuoka,
\emph{A diffusion process on a fractal},
in: Probabilistic Methods in Mathematical Physics (Katata/Kyoto, 1985),
Academic Press, Boston, MA, 1987, pp.~251--274.


\bibitem{KusuokaZhou}
S. Kusuoka and X.~Y.~Zhou, \emph{Dirichlet forms on fractals: Poincaré constant and resistance,} Probab. Theory Related Fields {\bf93} (1992), no. 2, 169--196.


\bibitem{KuwaeShioya2003}
K.~Kuwae and T.~Shioya,
\emph{Convergence of spectral structures: a functional analytic theory
and its applications to spectral geometry},
Comm. Anal. Geom. \textbf{11} (2003), no.~4, 599--673.

\bibitem{Lindstrom1990}
T.~Lindstr{\o}m,
\emph{Brownian Motion on Nested Fractals},
Mem. Amer. Math. Soc. \textbf{83} (1990), no.~420, iv+128 pp.

\bibitem{MathieuPiatnitski2007}
P.~Mathieu and A.~Piatnitski,
\emph{Quenched invariance principles for random walks on percolation clusters},
Proc. R. Soc. Lond. Ser. A Math. Phys. Eng. Sci. \textbf{463} (2007), 2287--2307.




\bibitem{Mosco1969}
U.~Mosco,
\emph{Convergence of convex sets and of solutions of variational inequalities},
Adv. Math. \textbf{3} (1969), no.~4, 510--585.

\bibitem{MuruganShimizu2025}
{M.~Murugan and R.~Shimizu,
\emph{First-order Sobolev spaces, self-similar energies and energy measures on the Sierpi\'nski carpet},
Comm. Pure Appl. Math. \textbf{78} (2025), no.~9, 1523--1608.}

\bibitem{Peirone2000}
R.~Peirone,
\emph{Convergence and uniqueness problems for Dirichlet forms on fractals},
Boll. Unione Mat. Ital. (8) \textbf{3-B} (2000), no.~2, 431--460.

\bibitem{Sabot1997}
C.~Sabot,
\emph{Existence and uniqueness of diffusions on finitely ramified self-similar fractals},
Ann. Sci. \'{E}cole Norm. Sup. (4) \textbf{30} (1997), no.~5, 605--673.



\bibitem{Shimizu2024}
{R.~Shimizu,
\emph{Construction of $p$-energy and associated energy measures on Sierpi\'nski carpets},
Trans. Amer. Math. Soc. \textbf{377} (2024), no.~2, 951--1032.}

\bibitem{SidoraviciusSznitman2004}
V.~Sidoravicius and A.-S. Sznitman,
\emph{Quenched invariance principles for walks on clusters of percolation or among random conductances},
Probab. Theory Related Fields \textbf{129} (2004), no.~2, 219--244.
\end{thebibliography}
\end{document}